\documentclass[11pt,twoside]{article}
\usepackage{amsmath,amsfonts,amssymb,cases,mathdots,color,colordvi,url,tikz}
\usepackage{hyperref}
\usepackage{bm}
\usepackage{bbm}
\usepackage{pgfpages}
\usepackage{tabularx}
\usepackage{graphics}%\textwidth 16.cm
\usepackage{subfigure}
\usepackage{epsfig}
\usepackage{algorithm}
\usepackage{algorithmic}
\usepackage{physics}
\usepackage{epstopdf}
\usepackage{parskip}
\usepackage{footmisc}
\usepackage{tablefootnote}
\usepackage[justification=centering]{caption}

\pdfoutput=1
\newtheorem{proposition}{Proposition}%[section]
\newtheorem{theorem}[proposition]{Theorem}
\newtheorem{lemma}[proposition]{Lemma}

\newtheorem{definition}[proposition]{Definition}

\newtheorem{remark}{Remark}
\newtheorem{assumption}{Assumption}
\newtheorem{example}{Example}

\newcommand{\be}{\begin{equation}}
	\newcommand{\ee}{\end{equation}}
\newcommand{\ba}{\begin{eqnarray}}
	\newcommand{\ea}{\end{eqnarray}}
\newcommand{\bas}{\begin{eqnarray*}}
	\newcommand{\eas}{\end{eqnarray*}}

\graphicspath{ {./Figures/} }
\def\R{{\mathcal R}}
\def\C{{\mathcal C}}
\def\F{{\mathcal F}}
\def\I{{\mathcal I}}
\def\cL{{\mathcal L}}
\def\M{{\mathcal M}}
\def\N{{\mathcal N}}

\def\P{{\mathcal P}}
\def\V{{\mathcal V}}

\def\bfp{{\bf p}}
\def\bfq{{\bf q}}

\def\bfs{{\bf s}}

\def\bfa{{\bf a}}
\def\bfb{{\bf b}}

\def\bfd{{\bf d}}

\def\bfv{{\bf v}}
\def\bfx{{\bf x}}
\def\bfy{{\bf y}}
\def\bfz{{\bf z}}
\def\bfw{{\bf w}}
\def\bfu{{\bf u}}

\def\Prox{\mbox{Prox}}
\def\OG{{\overline{\Gamma}}}

\def\wz{{\widetilde{\zeta}}}

\title{iNALM: An inexact Newton Augmented Lagrangian Method for Zero-One Composite Optimization}
\author{Penghe Zhang\thanks{School of Mathematics and Statistics, Beijing Jiaotong University, Beijing 100044, PR China, E-mail: {19118011@bjtu.edu.cn} },
\ \ Naihua Xiu\thanks{School of Mathematics and Statistics, Beijing Jiaotong University, Beijing 100044, PR China, E-mail: {nhxiu@bjtu.edu.cn} } \ \ and \ \
Hou-Duo Qi\thanks{Department of Applied Mathematics, The Hong Kong Polytechnic University, Hung Hom, Hong Kong, E-mail: {houduo.qi@polyu.ac.hk} }
	}

\date{}

\begin{document}
	
	\maketitle	

\begin{abstract}
	Zero-One Composite Optimization (0/1-COP) is a prototype of nonsmooth, nonconvex optimization problems and it has attracted much attention recently.
	The augmented Lagrangian Method (ALM) has stood out as a leading methodology for such problems. 
	The main purpose of this paper is to extend the classical theory of ALM from
	smooth problems to 0/1-COP.
	We propose, for the first time, second-order optimality conditions for 0/1-COP.
	In particular, under a second-order sufficient condition (SOSC), we prove R-linear convergence rate of the proposed ALM. 
	In order to identify the subspace used in SOSC, we employ the proximal operator of
	the 0/1-loss function, leading to an active-set identification technique.
	Built around this identification process, we design practical stopping criteria for 
	any algorithm to be used for the subproblem of ALM. 
	We justify that Newton's method is an ideal candidate for the subproblem and
	it enjoys both global and local quadratic convergence.
	Those considerations result in an inexact Newton ALM (iNALM).
	The method of iNALM is unique in the sense that it is active-set based, it is inexact (hence more practical), and SOSC 
	%instead of widely assumed 	Kurdyka-$\L$ojasiewicz (K$\L$) properties 
	plays an important role in its R-linear convergence analysis. 
	The numerical results on both simulated and real datasets show the fast running speed and high accuracy of iNALM when compared with several leading solvers.
	
	\noindent{\bf \textbf{Keywords}:}
	zero-one composite optimization problem, 
	inexact augmented Lagrangian method,
	second-order optimality conditions,
	P-stationary point, 
	convergence rate, 
	subspace Newton method.
%	\subclass{90C26 \and 90C30 \and 90C90}
\end{abstract}

\section{Introduction} 

We are concerned with the zero-one Composite Optimization Problem (0/1-COP):
\begin{equation} \label{COP}
	\min_{\bfx \in \mathbb{R}^n} f(\bfx) + \lambda h( A \bfx + \bfb) 
\end{equation}
where $f:\mathbb{R}^n \to \mathbb{R}$ is a smooth function, $\lambda>0$ is a penalty parameter,  and $A \in \mathbb{R}^{m\times n}, \bfb\in \mathbb{R}^m$
are given. Moreover, the function $h: \mathbb{R}^m \to \mathbb{R}$ counts the number of positive entries of $\bfu$:
\[
h(\bfu) := \| \bfu_+\|_0,
\]
where $\| \bfz\|_0$ is the $\ell_0$ quasi-norm of $\bfz$ counting its nonzero elements, and $\bfu_+ := \max\{ \bfu, 0\}$ (component-wise maximum).
The function $h(\cdot)$ is known as the $0/1$-loss function and  
the second part of the objective function in (\ref{COP}) is a composition between the $0/1$ loss and a linear operator.
Problem (\ref{COP}) arises from many applications including support vector machines (SVM) and multi-label classifications (MLC) \cite{CV95}, the one-bit compressed sensing \cite{BB08,5955138}, the maximum rank correlation \cite{han1987non}, and the problem of area under curves \cite{ma2005regularized}.  However, optimization related to the $0/1$-loss function is discontinuous and NP-hard, see \cite{ben2003difficulty}.
We refer to a recent paper \cite{zhou2021quadratic} for a more detailed discussion
on this type of problems from a combinatorial perspective.

The main purpose of this paper is to tackle (\ref{COP}) under the framework of the augmented Lagrangian method (ALM) with its subproblems being solved by a Newton method.
For the first time for a nonconvex, nonsmooth optimization problem, we are able to establish the global as well as local R-linear convergence of the
Newton Augmented Lagrangian method (NALM) under a second-order sufficient condition.
We will report promising numerical results on the two applications SVM and MLC mentioned 
above.
In the following, we first review some pertinent references that strongly motivated our research. We then describe our approach and main contributions.

\subsection{Literature review}
%There exists a large body of research on a variety of ALMs for smooth constrained convex optimization, see 
%\cite{hestenes1969multiplier, powell1969method, rockafellar1976augmented, bertsekas2014constrained, cui2019r}.
%There is also growing literature on nonsmooth, nonconvex optimization problems,
%see \cite{bolte2018nonconvex, zeng2021moreau} and the references therein. 
There exists a large body of research on a variety of ALMs for smooth constrained convex optimization, see 
\cite{andreani2008augmented2,bertsekas2014constrained,cui2019r,grapiglia2021complexity,alm_ma,powell1969method,rockafellar1976augmented}.
There is also growing literature on nonsmooth, nonconvex optimization problems,
see \cite{bolte2018nonconvex,zeng2021moreau} and the references therein. 
Our review here is to discuss only those which have close connections to or 
motivated our framework. 
We group them into four categories with a question in mind whether some of the nice 
properties of those reviewed methods may be extended to the problem (\ref{COP}).

{\bf (a) ALM under Second-Order Sufficient Condition (SOSC)}.
This is the classical approach to study the convergence properties of ALM for
smooth constrained optimization problems and is well explained, for example,
in \cite[Chps.~12 and 17]{nocedal2006numerical}. 
It has become a standard that this approach often requires certain regularity assumptions such as linear independence constraint qualification at the interested solution.
The implication of SOSC on convergence also extends to
Newton's method for generalized equations
\cite{fischer2002local} and sequential quadratic programming (SQP) methods
\cite{fernandez2012local,izmailov2012stabilized}.
A key element of SOSC in such extension is that the perturbed 
KKT conditions enjoy an upper Lipschitzian property. 
A question to us is whether we may extend this to the nonsmooth case like problem (\ref{COP}).
This forms one of our main contributions in this paper that an appropriate SOSC is proposed at a stationary point of (\ref{COP}) and the convergence of ALM is established 
under this SOSC.

{\bf (b) Newton ALM for convex optimization}.
This represents a major advance on ALM for structural convex optimization with its subproblems being solved by a highly efficient Newton's method, see \cite{9053722,li2018highly,li2020asymptotically,lin2021augmented} for a few examples. 
An essential observation in this framework is that each Newton equation has a structural 
sparse pattern so that the equation can be solved in reduced dimensions (often in very low-dimensional space). 
The resulting algorithm is carefully designed to ensure global and fast local convergence.
Given $\ell_0$-norm promoting sparsity in the solutions, we ask whether 
an efficient Newton ALM can be developed for the problem (\ref{COP}). 
An important issue to be resolved is the subspace tracking so that Newton's equation may be properly defined in a low-dimensional space.
We will design a novel technique based on the proximal operator of the function $h(\cdot)$ to identify the active set leading to the correct subspace being identified.

{\bf (c) Proximal ALM for nonconvex optimization}.
There is a significant progress in this part recently.
The paper \cite{bolte2018nonconvex}, which includes (\ref{COP}) as a special case, provides us with a deep understanding of how ALM would behave 
for nonconvex, nonsmooth optimization problems.
It recognizes that for ALM to converge there needs a certain regularity assumption such
as the full row-rank of the matrix $A$ in (\ref{COP}),
see also \cite[Example 7]{li2015global} and \cite[Remark 2]{bot2019proximal}.
In addition to the two standard steps (Primal step and Multiplier step) in ALM,
the generic algorithm (called ALBUM) in \cite{bolte2018nonconvex} also 
introduces a third step (Adaptive step) by the means of a Lyapunov function.
The Lagrangian algorithmic map in its primal step needs to satisfy certain 
conditions to ensure convergence.
For example, a proximal alternating direction method of multipliers (ADMM) can be used.
One such method is \cite{bot2019proximal}, which also addressed the boundedness of 
the primal-dual iterates provided that $f(\cdot)$ is coercive (i.e., $f(\bfx) \rightarrow \infty$ as $\| \bfx\| \rightarrow \infty$). 
If $h(\cdot)$ is assumed weakly convex, the Moreau envelope ALM in \cite{zeng2021moreau}
can be applied. Unfortunately, the $\ell_0$-norm is not weakly convex.
The innovative techniques introduced in those papers 
(e.g., Lyapunov function and Lagrangian arithmetic map \cite{bolte2018nonconvex},
Moreau envelop \cite{zeng2021moreau} and boundedness analysis \cite{bot2019proximal})
are motivating us to ask whether they can be applied to develop a Newton ALM for
problem (\ref{COP}) with both global and local R-linear convergence.
Though sounding too demanding, it is exactly what we are going to achieve in this paper.

{\bf (d) Hard thresholding for sparse optimization}.
Hard thresholding has emerged as an important technique for subspace pursuit in
sparse optimization, see e.g., \cite{beck2013sparsity,blumensath2008gradient}:
\be \label{SparseOpt}
\min_{\bfx \in \mathbb{R}^n} \ f(\bfx) \quad \mbox{s.t.} \quad 
\bfx\in \C := \left\{   \bfx \; | \ \| \bfx\|_0 \le s  \right\},
\ee 
where $f: \mathbb{R}^n \mapsto \mathbb{R}$ is smooth and $s>0$ is a given integer to control
the sparsity level in $\bfx$.
It is observed in \cite{zhou2021global} that the iterative hard-thresholding update
\[
\bfx^{k+1} \in \P_{\C}(  \bfx^k - \eta \nabla f(\bfx^k) )
\]
also implicitly updates the active indices of the sparse vector $\bfx^{k+1}$. 
Here, $\P_{\C}$ is the orthogonal projection operator onto $\C$ and $\eta>0$ is step-length. 
Once the active set is identified, a subspace Newton method can be developed for (\ref{SparseOpt}) with a globalization strategy.
This idea was then extended to problem (\ref{COP}) in \cite{zhou2021quadratic}, which can only establish its local
quadratic convergence under the assumption of $f$ being strongly convex and $A$ being of
full row rank.
Relevant extension to sparse optimization (\ref{SparseOpt}) with extra equality constraints has been done in \cite{zhao2021lagrange} via the Lagrange function.
Inspired by those developments, we ask whether a thresholding strategy can be developed
in the framework of ALM so that Newton's method can be used in solving its subproblems.
We will see that the corresponding technique in the ALM context is related to the proximal operator of the function $h(\cdot)$.

%%%%%%%%%%%%%%%%%%%%%%%%%%%%%%%%%%%%%%%%%%%%%%%%
\subsection{Our approach and main contributions}

It has become roughly clear what properties we would like our ALM to enjoy:
(i) Its subproblems are to be solved by Newton's method in a subspace;
(ii) The generated iterates sequence is bounded; and
(iii) It enjoys global as well as local $R$-linear convergence under a suitable
SOSC. We outline the approach how we arrive at such a method.

Firstly, we adopt the strategy used in \cite{bolte2018nonconvex,bot2019proximal}
to reformulate (\ref{COP}) as a constrained optimization problem:
\be \label{COP-Constrained}
\min_{\bfx, \bfu} \ f(\bfx) + \lambda h(\bfu) , \quad
\mbox{s.t.} \quad A\bfx + \bfb = \bfu.
\ee 
Let $ (\bfx^*, \bfu^*, \bfy^*)$ be a Karush-Kuhn-Tucker (KKT) point of (\ref{COP-Constrained}) with
$\bfy^*$ being the Lagrange multiplier.
Our key observation is that there exists a corresponding KKT point for the
following smooth optimization problem:
\be \label{Smooth-P}
\min_{\bfx, \bfu} \ f(\bfx), \quad \mbox{s.t.} \quad \bfu_{\I_-^*} \le 0,
\quad A \bfx + \bfb = \bfu,
\ee 
where $\I^*_- := \left\{ i \; | \ u^*_i \le 0 , \ i=1, \ldots, m\right\}$
is the index set of non-positive components of $\bfu^*$. 
Through this correspondence, we are able to propose SOSC for the problem (\ref{COP-Constrained}) as (\ref{Smooth-P}) is a smooth problem.
This is where the current research departs from the previous research \cite{zhang2021global}, where the strong convexity of $f(\bfx)$ is assumed.
The correspondence of KKT points between (\ref{COP-Constrained}) and (\ref{Smooth-P})
is established by studying $P$-stationary points of (\ref{COP-Constrained}),
whose detail is deferred to Sect.~\ref{Section-SOSC}.

Secondly, motivated by \cite{bolte2018nonconvex}, we use the Lyapunov function
$\V_{\rho, \mu}$ as our merit function to measure improvement of ALM:
\begin{eqnarray*}
	\V_{\rho, \mu} (\bfx, \bfu, \bfy, \bfv) &:=&
	\cL_{\rho} (\bfx, \bfu, \bfy) + \frac {\mu}2 \| \bfx - \bfv\|^2
	%	&=& \underbrace{f(\bfx) + \langle \bfy, A\bfx + \bfb - \bfu \rangle + \frac{\rho}2 
	%	\|A\bfx + \bfb - \bfu \|^2 +  \frac {\mu}2 \| \bfx - \bfv\|^2}_{:= g(\bfx, \bfu, \bfy, \bfv)} + \lambda h(\bfu) ,
\end{eqnarray*}
where $\mu >0 $ is the associated Lyapunov constant and $\cL_\rho(\bfx, \bfu, \bfy)$ is the augmented Lagrange function of 
(\ref{COP-Constrained}) with $\rho>0$ being a penalty parameter:
\[
\cL_{\rho} (\bfx, \bfu, \bfy) := 
f(\bfx) + \langle \bfy, A\bfx + \bfb - \bfu \rangle + \frac{\rho}2 
\|A\bfx + \bfb - \bfu \|^2  + \lambda h(\bfu),
\]
and $\bfv \in \mathbb{R}^n$ is called Lyapunov regularization vector.
In the particular choice of $\bfv = \bfx^k$, which is the current iterate of $\bfx$,
the Lyapunov function is simply the augmented Lagrangian added by a Euclidean proximal 
term on the $\bfx$-part.
The subproblem at the $k$th iteration of our ALM is to compute an approximate solution of
the  following problem:
\begin{align} \label{Subproblem}
	&(\bfx^{k+1}, \bfu^{k+1}) \approx
	\arg\min_{\bfx, \bfu} \ \V_{\rho, \mu} (\bfx, \bfu, \bfy^k, \bfx^k) \\
	&= \arg\min_{\bfx, \bfu} \; \underbrace{f(\bfx) + \langle \bfy^k, A\bfx + \bfb - \bfu \rangle + \frac{\rho}2 
		\|A\bfx + \bfb - \bfu \|^2 +  \frac {\mu}2 \| \bfx - \bfx^k\|^2}_{=:\; g_k(\bfx, \bfu)} + \lambda h(\bfu) \nonumber
\end{align}
and we update the Lagrange multiplier $\bfy^k$ accordingly.
Subproblem (\ref{Subproblem}) can be approximately solved by any gradient descent or alternating 
minimization method.
Our proposal is to add a Newton's step hoping to achieve fast convergence for the inner step. 

Thirdly, we recall the aim is to establish convergence of ALM under a SOSC associated with the problem (\ref{Smooth-P}).
An important step is to identify the active set of those indices that $u^*_i =0$, $i=1, \ldots, m$. 
We make use of the proximal operator of $h(\bfu)$ to achieve this purpose: Let
\[
\Gamma_k := \left\{
i \ \left | \ 0 \in \left[
\Prox_{\alpha \lambda h(\cdot)} \Big(
\bfu^k - \alpha \nabla_{\bfu} g_k(\bfx^k, \bfu^k)
\Big)
\right]_i, \ i=1, \ldots, m \right.
\right\},
\]
where $\alpha>0$ is a steplength,  $\Prox_{\psi(\cdot)}(\bfu)$ is the proximal operator of a given function $\psi: 
\mathbb{R}^m \to \mathbb{R}$,
and for a given set $\Omega$, $[\Omega]_i$ denotes the set of
its $i$th elements.
Since $h(\bfu)$ is nonsmooth, its proximal operator may be a set, see its detailed 
description in Sect.~\ref{Section-SOSC}.
Consequently, we will solve subproblem (\ref{Subproblem}) over the subspace $\bfu_{\Gamma_k} =0$.

However, the introduction of the active-set identification strategy to the standard
ALM framework also brings a challenging task how to ensure global as well as local convergence of ALM. 
As we will see, the analysis is quite technical. 
In order to reduce the technical 
complexity, we will assume the objective function $f(\bfx)$ twice continuously differentiable with Lipschitzian properties. 
We list some of the main results below.

\begin{itemize}
	\item[(i)] We characterize the KKT point of (\ref{COP-Constrained}) in terms of P-stationary points (Prop.~\ref{Prop-P}) and
	prove that any local minimizer of (\ref{COP-Constrained}) is a P-stationary point
	(Prop.~\ref{Prop-First-Order-Optimality}). The latter result removes the previously required assumption that $A$ has full row rank used in \cite{bolte2018nonconvex,bot2019proximal,zhou2021quadratic}.
	
	\item[(ii)] We directly derive the second-order optimality condition of 0/1-COP from the reformulation \eqref{Smooth-P}. If SOSC is satisfied at a P-stationary point, then the stationary point must be a local minimum and a quadratic growth condition holds at this point (Thm.~\ref{Thm-SOSC}). 
	Moreover, a local minimizer must satisfy a second-order necessary condition
	(a weak version of SOSC), see Thm.~\ref{Thm-SONC}. 
	Those results extend the classical second-order optimality conditions for smooth problems in \cite{nocedal2006numerical} to (\ref{COP-Constrained}).
	
	\item[(iii)] We develop an inexact Augmented Lagrangian Method (iNALM), whose subproblem is solved by a subspace Newton method.
	The subspace is identified by means of the proximal operator of $h(\bfu)$. 
	The Newton method is globally convergent for $f$ being weakly convex
	(Thm.~\ref{Thm-Global})
	and 
	has a local quadratic convergence rate if a strict complementarity condition holds
	(Thm.~\ref{Thm-Quadratic}).
	Newton's method is terminated as soon as a set of stopping criteria is met.
	The criteria, motivated by the upper Lipschitzian property of the KKT solution map of (\ref{COP-Constrained}), are carefully designed so that
	the resulting iNALM is globally convergent to a P-stationary point of (\ref{COP-Constrained}) under standard assumptions, see Thm.~\ref{Thm-Stationarity}.
	In particular, when $f$ is coercive, the generated KKT sequence
	$\{(\bfx^k, \bfu^k, \bfy^k)\}$ is bounded (Lemma~\ref{Lemma-Boundedness}).
	Furthermore, if SOSC is satisfied at a P-stationary point, then both 
	the Lyapunov sequence $\{ \V_{\rho, \beta} (\bfx^k, \bfu^k, \bfy^k, \bfx^k) \}$
	($\beta$ is properly chosen and different from $\mu$ )
	and the KKT sequence enjoy a R-linear convergence rate, see Thms.~\ref{Thm-Lyapunov} and \ref{Thm-R-Iterates}. 	
	%Those results are in contrast to that of \cite{bot2019proximal}, where K$\L$ properties of $f(\cdot)$ and $h(\cdot)$ are assumed.
\end{itemize}

We also note that the inexact ALM is more practical for implementation than its exact
counterpart. 
%%%%%%%%%%%%%%%%%%%%%%%%%%%%%%%%%%%%%%%%%%%%%%%%%
\subsection{Organization}
The paper is organized as follows.
In the next section, we describe the P-stationarity and second-order conditions
for problem (\ref{COP-Constrained}) and study their consequences.
In Section~\ref{Section-Newton}, we propose a subspace Newton's method aimed for
the subproblem arising from ALM and prove both its global and local quadratic convergence under suitable conditions.
In Section~\ref{Section-NALM}, we present our main iNALM algorithm, 
proving that it is well-defined and it enjoys nice convergence properties as listed above.
Extensive numerical experiments are reported in Section~\ref{Section-Numerical}, confirming the efficiency of iNALM. 
%We conclude the paper in Section~\ref{Section-Conclusion}.
%\section{Section title}
%\label{sec:1}
%Text with citations \cite{RefB} and \cite{RefJ}.
%\subsection{Subsection title}
%\label{sec:2}
%as required. Don't forget to give each section
%and subsection a unique label (see Sect.~\ref{sec:1}).
%\paragraph{Paragraph headings} Use paragraph headings as needed.
%\begin{equation}
%\end{equation}

%%%%%%%%%%%%%%%%%%%%%%%%%%%%%%%%%%%%%%%%%%%%%%%%%%%%%%%%%%%%%%%%%%%%%%%%%%%%%%%
\section{P-Stationarity and Second-Order Conditions} \label{Section-SOSC}

\subsection{Notation and preliminaries}

We use boldfaced lowercase letters to denote vectors. For example, $\bfx \in \mathbb{R}^n$ is a column vector of size $n$ and $\bfx^\top$ is its transpose.
Let $x_i$ or $[\bfx]_i$ denote the $i$th element of $\bfx$. 
The norm $\| \bfx\|$ denotes the Euclidean norm of $\bfx$ and for
a matrix $A$, $\|A\|$ is the induced norm by the Euclidean norm so that we always have
$
\| A\bfx \| \le \| A \| \| \bfx\|.
$ $\langle \cdot \rangle$ is the inner product in Euclidean space.
For two column vectors $\bfx$ and $\bfy$, we use the matlab notation $[\bfx; \bfy]$ to 
denote the new column vector concatenating $\bfx$ and $\bfy$.
The neighborhood of $\bfx^* \in \mathbb{R}^n$ with radius $\delta > 0$ is denoted by $\mathcal{N}(\bfx^*, \delta) := \{ \bfx \in \mathbb{R}^n \ | \ \| \bfx - \bfx^* \| < \delta \}$, where ``$:=$'' means ``define''.
We let
$
I
$ denote the identity matrix of the appropriate dimension.

Let $[m]$ denote the set of indices $\{1, \ldots, m\}$.
For a subset $\Gamma \subset [m]$,  $\OG$ consists of those indices of $[m]$ not in $\Gamma$ and $|\Gamma|$ denotes the number of elements in $\Gamma$ (cardinality of $\Gamma$).
For $\bfy \in \mathbb{R}^m$, $\bfy_{\Gamma}$ denotes the subvector of
$\bfy$ indexed by $\Gamma$.
Let
$
\F
$ denote the feasible region of (\ref{COP-Constrained}) and
$
\F^*_- : = \left\{
(\bfx, \bfu) \ |\ (\bfx, \bfu) \in \F \ \mbox{and} \
u_i \le 0, \ i \in \I^*_-
\right\}
$ (feasible region of the smooth problem (\ref{Smooth-P})).

Suppose $\psi: \mathbb{R}^n  \to \mathbb{R}$ is lower semi-continuous (lsc).
In this paper, we only consider such $\psi$ that is bounded from below.
For given $\alpha >0$, 
the Moreau envelop and the proximal operator of $\lambda\psi$ are respectively defined by
\begin{eqnarray*}
	\Phi_{\alpha\psi(\cdot)} (\bfx) 
	&:=& \min_{\bfz \in \mathbb{R}^n} \left\{ \alpha \psi(\bfz) 
	+ \frac 12 \| \bfz - \bfx\|^2 \right \} \\
	\Prox_{\alpha \psi(\cdot)}(\bfx) &:=& \arg\min_{\bfz \in \mathbb{R}^n}\left\{
	\alpha \psi(\bfz) 
	+ \frac 12 \| \bfz - \bfx\|^2
	\right\}.
\end{eqnarray*} 
Suppose there is sequence $\bfx^k \rightarrow \bfx^*$ with $\bfz^k \in \Prox_{\alpha\psi(\cdot)} (\bfx^k)$ and $\bfz^k \rightarrow \bfz^*$.
The Theorem of Proximal Behavior \cite[Thm.~1.24]{RockWets98} says
$\bfz^* \in \Prox_{\alpha\psi(\cdot)} (\bfx^*)$.

We say $\psi: \mathbb{R}^n  \to \mathbb{R}$ is $\ell_\psi$-smooth if its gradient is
Lipschitz continuous with modulus $\ell_\psi$:
\[
\|\nabla \psi(\bfx) - \nabla \psi(\bfz)\| \le \ell_\psi \| \bfx - \bfz\|,
\quad \forall \ \bfx, \bfz \in \mathbb{R}^n.
\] 

\begin{lemma}[Descent Lemma (\cite{beck2017first})] \label{Descent-Lemma}
	Let $\psi: \mathbb{R}^n \to (- \infty, \infty]$ be a $\ell_\psi$-smooth function. Then 
	the following holds:
	\[
	\psi(\bfu) \leq \psi(\bfv) + \langle \nabla \psi(\bfv), \bfu - \bfv \rangle + \frac{\ell_\psi}{2} \| \bfu - \bfv \|^2, \quad \forall \ \bfu, \ \bfv \in \mathbb{R}^n.
	\]
	In particular, we have
	\[
	\psi\Big( \bfu - (1/\ell_\psi) \nabla \psi(\bfu) \Big)
	\le \psi(\bfu) - \frac{1}{2 \ell_\psi} \| \nabla \psi(\bfu) \|^2 .
	\]
\end{lemma}

%%%%%%%%%%%%%%%%%%%%%%%%%
\subsection{P-stationary points}

In this part, we define the P-stationary point of (\ref{COP-Constrained}) and study its
relationship with the KKT point of (\ref{COP-Constrained}).
We then prove that any local minimum of (\ref{COP-Constrained}) is also
a P-stationary point.

Define the Lagrange function of (\ref{COP-Constrained}) by
\[
\cL(\bfx, \bfu, \bfy) := f(\bfx) + \lambda h(\bfu) + \langle \bfy, \; A\bfx + \bfb - \bfu \rangle.
\]
We recall that a KKT point $\bfz^* = (\bfx^*, \bfu^*, \bfy^*)$ of (\ref{COP-Constrained})
satisfies
\be \label{KKT}
\left\{
\begin{array}{l}
	0 = \nabla_{\bfx} \cL(\bfz^*) = \nabla f(\bfx^*) + A^\top \bfy^* \\
	0 \in \partial_{\bfu} \cL(\bfz^*) = - \bfy^* + \lambda \partial h(\bfu^*) \\
	0 = A\bfx^* + \bfb - \bfu^*,
\end{array} 
\right .
\ee 
where $\partial h(\cdot)$ is the limiting subdifferential defined for lsc functions \cite[Def.~8.3]{RockWets98},
and it has a very simple structure, see \cite[Eq.~(2.1)]{zhou2021quadratic}:
\be \label{Partial_h}
\partial h(\bfu) =
\left\{
\bfv \in \mathbb{R}^m \; \left|
\ v_i \left\{
\begin{array}{ll}
	\ge 0 , \ & \mbox{if} \ u_i = 0 \\
	=0,     \ & \mbox{if} \ u_i \not=0,
\end{array} , \ \ i \in [m]
\right\}
\right .
\right\}
\ee 
A useful fact is that a complementarity condition holds for any pair $(\bfu, \bfv)$ with
$\bfv \in \partial h(\bfu)$:
\be \label{Complementarity-Condition}
u_i v_i =0, \ i =1, \ldots, m \ \ \mbox{for any} \ \bfv \in \partial h(\bfu).
\ee 
We note that the second condition in (\ref{KKT}) holds if and only if
\[
\alpha \bfy^* \in \alpha \lambda \partial h(\bfu^*) \quad \forall \ \alpha > 0.
\]
A sufficient condition for this to hold is
\begin{eqnarray*}
	\bfu^* &\in& \arg\min_{\bfu} \left\{
	(\alpha \lambda) h(\bfu) + \frac 12 \| \bfu - ( \bfu^* + \alpha \bfy^* ) \|^2 
	\right\} \\
	&=& \Prox_{\alpha \lambda h(\cdot)} ( \bfu^* + \alpha \bfy^* ),
\end{eqnarray*}
where the proximal operator $\Prox_{\alpha \lambda h(\cdot)}$ of $h(\cdot)$ is well-defined because $h(\cdot)$ is bounded from below.
We refer the interested reader to \cite[Chp.~6]{beck2017first} for detailed study of
proximal operators of various functions.
This motivates us to define the P-stationary point of (\ref{COP-Constrained}).

\begin{definition} \label{Def-P}
	(P-stationary point)
	A point $ (\bfx^*, \bfu^*, \bfy^*) \in \mathbb{R}^n \times \mathbb{R}^m \times \mathbb{R}^m $ is called a P-stationary triplet of
	(\ref{COP-Constrained}) if there exists $\alpha >0$ such that
	\begin{equation} \label{P-stat}
		\left\{ \begin{aligned}
			& \nabla f(\bfx^*) + A^\top \bfy^* = 0,\\
			& \bfu^* \in {\rm Prox}_{\alpha\lambda h(\cdot)} ( \bfu^* + \alpha \bfy^* ), \\
			& A\bfx^* + \bfb - \bfu^* = 0.
		\end{aligned}  \right.
	\end{equation}
	Moreover, $(\bfx^*, \bfu^*)$ is called P-stationary point of (\ref{COP-Constrained})
	and $\bfy^*$ is the corresponding P-stationary multiplier.
\end{definition}

Apparently, a P-stationary triplet is also a KKT point of (\ref{COP-Constrained}).
Interestingly, the converse is also true. 
Moreover, we will characterize the P-stationary point in terms of KKT points of
the smooth optimization problem (\ref{Smooth-P}). 
In order to prepare the proofs, 
we first note that the proximal operator of $h(\cdot)$ can be easily 
computed via its definition (see, e.g., \cite[Eq.~(2.6)]{zhou2021quadratic}):
\begin{equation}\label{prox0+}
	\Big[ 	{\rm Prox}_{\alpha\lambda h(\cdot)} ( \bfu )  \Big]_i 
	= \left\{  \begin{aligned}
		& 0, && u_i \in (0,\; \sqrt{2\lambda\alpha}), \\
		&\{ 0, u_i \}, && u_i \in \{ 0,\; \sqrt{2\lambda \alpha} \}, \\
		&u_i, && u_i \in (-\infty,\; 0)\cup (\sqrt{2\lambda\alpha},\; \infty),
	\end{aligned} \right. \quad i \in [m].
\end{equation}
A straightforward application of this formula leads to the following implication:
\begin{equation} \label{Prox-eq}
	\begin{aligned}	
		& \bfu^* \in {\rm Prox}_{\alpha\lambda h(\cdot)} ( \bfu^* + \alpha \bfy^* )   \\
		& \quad \quad \quad \quad \quad \quad \Downarrow   \\	
		&\left\{ \begin{aligned}
			& u^*_i y^*_i = 0 \quad \mbox{(complementarity condition)}, \\
			& \mbox{if} \ y^*_i = 0,\ \mbox{then} \ u^*_i \in (-\infty,0] \cup [\sqrt{2\lambda\alpha}, \infty),   \\
			& \mbox{if} \ u^*_i = 0,\ \mbox{then} \ y^*_i \in [0,\sqrt{2\lambda/\alpha}],
		\end{aligned}  \right. ~\forall i \in [m]. 
	\end{aligned}
\end{equation}
The complementarity condition follows from (\ref{Complementarity-Condition}) and the fact
that $\bfu^* \in {\rm Prox}_{\alpha\lambda h(\cdot)} ( \bfu^* + \alpha \bfy^* )$ is
sufficient for $\bfy^* \in \partial h(\bfu^*)$.
We say such a pair $(\bfu^*, \bfy^*)$ satisfies the strict complementarity condition if
$
u^*_i + y^*_i \not= 0 
$ 
for all $ i \in [m]$.
%We note that the strict complementarity condition does not hold 
%if and only if there exists $i \in [m]$ such that $u^*_i = y^*_i =0$.

We summarize equivalent characterizations of P-stationary point below.

\begin{proposition} \label{Prop-P}
	(Equivalent characterizations of P-stationary point)
	Let $\bfw^* := (\bfx^*, \bfu^*, \bfy^*)$ be a reference point. The following hold.
	\begin{itemize}
		\item[(i)] If $\bfw^*$ is a P-stationary triplet of (\ref{COP-Constrained}), then
		it is also a KKT point. Conversely, if $\bfw^*$ is a KKT point of (\ref{COP-Constrained}), then it is a P-stationary triplet with any $\alpha \in (0, \alpha^*)$, where $\alpha^* := \min\{
		\alpha_u^*, \; \alpha_y^*
		\}$ and
		\begin{equation} \label{alpha_star}
			\begin{aligned}
				&	
				\alpha^*_{u} := \left\{  
				\begin{aligned}
					& +\infty,  && \bfu^*\leq 0, \\
					&\min\left\{ \frac{(u^*_i)^2}{ 2\lambda }:\ u_i^* > 0 \right\}, &&{\rm otherwise},
				\end{aligned} \right. ,\\
				&	\alpha^*_{y} := \left\{  \begin{aligned}
					& +\infty, && \bfy^*\leq 0, \\
					&\min\left\{ \frac{2\lambda}{ (y^*_i)^2 } : \ y^*_i > 0 \right\}, && {\rm otherwise}. 
				\end{aligned} \right.
			\end{aligned}
		\end{equation}
		
		\item[(ii)] $\bfw^*$ is a KKT point of (\ref{COP-Constrained}) if and
		only if $(\bfx^*, \bfu^*, \bfy^*_{\I^*_-}, \bfy^*)$ is a KKT point of the
		smooth problem (\ref{Smooth-P}).
	\end{itemize}
\end{proposition}

{\bf Proof.}
(i) As mentioned above, P-stationary triplet is a KKT point of (\ref{COP-Constrained}). 
We prove the converse part.
It is sufficient to prove,
for a given KKT point $(\bfx^*, \bfu^*, \bfy^*)$, that
$\bfu^* \in \Prox_{\alpha \lambda h(\cdot)}(\bfu^* + \alpha \bfy^*)$ for any
$\alpha \in (0, \alpha^*)$.
For any $i \in [m]$,
we consider three cases that correspond to the three scenarios in
(\ref{prox0+}).
%can be similarly proved by following the part of the proof
%from from Eq.~(3.14) in \cite[Thm.~3.3]{zhou2021quadratic}.

\begin{itemize}
	\item[(C1)] Assume $u_i^* \not=0$. By the complementarity condition (\ref{Complementarity-Condition}), we have $y^*_i =0$.
	If $u_i^* >0$, then by the definition of $\alpha^*$, we have
	\[
	\alpha < \alpha^* \le \frac{(u^*_i)^2}{2 \lambda},
	\]
	which implies $u^*_i > \sqrt{2 \lambda \alpha}$. Hence, in this case
	$u_i^* + \alpha y^*_i = u^*_i \in (-\infty, 0)$ or $(\sqrt{2 \lambda \alpha}, \infty)$. 
	By the third scenario in (\ref{prox0+}), we have
	$
	u^*_i \in [\Prox_{\alpha \lambda h(\cdot)}(\bfu^*+ \alpha \bfy^*)]_i.
	$
	
	\item[(C2)] Assume $y_i^* \not=0$. Once again, by the complementarity condition
	(\ref{Complementarity-Condition}), we have $u^*_i =0$. By the definition of $\alpha^*$, we have
	\[
	\alpha < \alpha^* \le \alpha_y^* \le  \frac{2 \lambda}{(y^*_i)^2},
	\]
	which implies $y^*_i < \sqrt{2 \lambda /\alpha}$. Hence, in this case
	$u_i^* + \alpha y^*_i = \alpha y^*_i \in (0, \sqrt{2 \lambda \alpha})$
	(note $y^*_i \ge 0$) 
	By the first scenario in (\ref{prox0+}), we have
	$
	u^*_i = 0 \in [\Prox_{\alpha \lambda h(\cdot)}(\bfu^*+ \alpha \bfy^*)]_i.
	$
	
	\item[(C3)] Assume $y^*_i =0$ and $u^*_i =0$. This falls within the second scenario 
	in (\ref{prox0+}), which implies 
	$
	u^*_i = 0 \in [\Prox_{\alpha \lambda h(\cdot)}(0)]_i 
	= [\Prox_{\alpha \lambda h(\cdot)}(\bfu^*+ \alpha \bfy^*)]_i.
	$
	
\end{itemize}
Putting the three cases together, we proved
$\bfu^* \in \Prox_{\alpha \lambda h(\cdot)}(\bfu^* + \alpha \bfy^*)$.

(ii) Let $\bfs \in \mathbb{R}^{|\I_-^* |}$ and $\bfy \in \mathbb{R}^m$ be the Lagrange multipliers of the smooth problem (\ref{Smooth-P}) corresponding to its two constraints.
It is easy to see that $\bfs = \bfy_{\I^*_-}$. 
Thus, any KKT point of the
smooth problem (\ref{Smooth-P}) takes the following form $(\bfx, \bfu, \bfy_{\I^*_-}, \bfy)$ satisfying
\begin{equation}   \label{KKT-Smooth}
	\left\{  \begin{aligned}
		& \nabla f(\bfx) + A^\top \bfy = 0,\\
		& \bfy_{\overline{\I}^*_-} = 0, \  \bfy_{\I^*_-} \geq 0, \
		\bfu_{\I^*_-} \leq 0, \ 
		\langle \bfy_{\I^*_-}, \; \bfu_{{\I^*_-}}  \rangle = 0, \\
		& A \bfx + \bfb - \bfu = 0 .
	\end{aligned}  \right.
\end{equation}
By using the structure (\ref{Partial_h}) of $\partial h(\bfu^*)$, 
the second condition in (\ref{KKT-Smooth}) is exactly an equivalent representation of
$\bfy \in \partial h(\bfu)$. Using the fact $\partial h(\bfu) = \lambda \partial h(\bfu)$ as $\partial h(\bfu)$ is a cone, we see that
(\ref{KKT-Smooth}) is an equivalent representation of the KKT point for (\ref{COP-Constrained}). 
The proof is completed by noting that any P-stationary triplet is a KKT point.
\hfill $\Box$  

The importance of the characterization in terms of the smooth problem (\ref{Smooth-P}) in
Prop.~\ref{Prop-P} is reflected in the following result.

\begin{proposition} \label{Prop-First-Order-Optimality}
	(First-order optimality condition)
	Suppose $(\bfx^*, \bfu^*) \in \mathbb{R}^{n} \times \mathbb{R}^{m}$ is a local optimal
	solution of (\ref{COP-Constrained}). 
	The following hold.
	
	\begin{itemize}
		\item[(i)] $(\bfx^*, \bfu^*)$ is also a local minimum of the smooth problem
		(\ref{Smooth-P}).
		
		\item[(ii)] There exists a Lagrange multiplier $\bfy^* \in \mathbb{R}^m$ such that
		$(\bfx^*, \bfu^*, \bfy^*)$ is a KKT point of (\ref{COP-Constrained}).
		Consequently, $(\bfx^*, \bfu^*)$ is a P-stationary point with any $\alpha \in (0, \alpha^*)$, where $\alpha^* = \min\{ \alpha_u^*, \alpha_y^*\}$ is defined in (\ref{alpha_star}).
		
	\end{itemize} 
\end{proposition}

{\bf Proof.}
(i) Since $(\bfx^*, \bfu^*)$ is the local minimizer of (\ref{COP-Constrained}), 
then there exists a neighbourhood of it, denoted by 
$\mathcal{N}((\bfx^*, \bfu^*), \delta)$ with $\delta>0$ such that for any feasible
point $(\bfx, \bfu)$ of (\ref{COP-Constrained})
\begin{align} 
	f(\bfx) + \lambda \| \bfu_+ \|_0 \geq f(\bfx^*) + \lambda \| \bfu^*_+ \|_0,  \ \forall \ (\bfx, \bfu) \in \mathcal{N}((\bfx^*, \bfu^*), \delta). \label{loc_cop}
\end{align} 
Now we restrict $\bfu$ to $\bfu_{\I^*_-} \le 0$ so that those feasible points
$(\bfx, \bfu)$ are also feasible to the problem (\ref{Smooth-P}).
Those restricted points also form a neighbourhood of $(\bfx^*, \bfu^*)$ for the
smooth problem (\ref{Smooth-P}).
For such a restricted point, we have  $\| \bfu_+ \|_0 \leq \| \bfu^*_+ \|_0$.
Consequently, $f(\bfx) \ge f(\bfx^*)$ 
from (\ref{loc_cop}). 
This means that $(\bfx^*, \bfu^*)$ is a local minimizer of (\ref{Smooth-P}).

(ii) We note that the objective function in (\ref{Smooth-P}) is smooth and constraints are linear.
Therefore, there must exist $\bfy^* \in \mathbb{R}^m$ such that
$(\bfx^*, \bfu^*, \bfy^*_{\I^*_-}, \bfy^*)$ is a KKT point of (\ref{Smooth-P}).
It follows from Prop.~\ref{Prop-P}(ii) that $(\bfx^*, \bfu^*, \bfy^*)$ is a
KKT point of (\ref{COP-Constrained}) and  Prop.~\ref{Prop-P}(i) implies $(\bfx^*, \bfu^*)$ is a P-stationary point with any $\alpha \in (0, \alpha^*)$.  \hfill $\Box$

\begin{remark} \label{Remark-P}
	The first-order optimality result in Prop.~\ref{Prop-First-Order-Optimality} generalizes the result \cite[Lemma~3.1]{zhou2021quadratic}, where it requires a submatrix $A$ having
	a full row rank.
	In general, when $h(\cdot)$ is lsc, a first-order optimality condition would require
	$A$ be of full row rank, see \cite[Def.~2.1 and Prop.~3.1]{bolte2018nonconvex}.
	Our characterization here makes use of the special structure of $h(\bfu)$, which allows
	us to drop the full row rank assumption on $A$.	
	It is also important to note that the P-stationarity characterization of a KKT point 
	establishes its link to the proximal operator of $h(\bfu)$, 
	which in turn provides us with a tool to identify the active set required to define second-order necessary and sufficient optimality conditions.
\end{remark}

%%%%%%%%%%%%%%%%%%%%%%%%%%%%%%%%%%%%%%%%%%%%
\subsection{Second-Order Conditions}

In this subsection, we will study the second-order conditions for the
smooth problem (\ref{Smooth-P}) and prove that those conditions also work for
(\ref{COP-Constrained}).

Suppose $(\bfx^*, \bfu^*, \bfy^*)$ be a KKT point of (\ref{COP-Constrained}).
In the proof of Prop.~\ref{Prop-First-Order-Optimality}(ii), we
proved that $(\bfx^*, \bfy^*, \bfy^*_{\I^*_-}, \bfy^*)$ is also a
KKT point of (\ref{Smooth-P}). 
We define the active sets of (\ref{Smooth-P}) at $(\bfx^*, \bfu^*)$ as follows:
\[
\I^*_0 := \left\{
i \ | \ u^*_i =0, \ y^*_i = 0, \ i\in [m] 
\right\} \quad \mbox{and} \quad
\I^*_+ := \left\{
i \ | \ u^*_i =0, \ y^*_i > 0, \ i\in [m] 
\right\}.
\]
By its definition (see, \cite[Eq.~(12.53)]{nocedal2006numerical},
the critical cone of (\ref{Smooth-P}) at $(\bfx^*, \bfu^*)$ is described as follow:
\be \label{Critical-Cs}
\C^*_s (\bfx^*, \bfu^*, \bfy^*) 
= \left\{
[\bfd^x; \bfd^u] \in \mathbb{Re}^{n+m} \ \Big| \ 
A \bfd^x = \bfd^u, \ \ \bfd^u_{\I^*_+} = 0, \ \
\bfd^u_{\I^*_0} \le 0
\right\}.
\ee 
We project this critical cone to its first components to get the following cone:
\be \label{Critical-C}
\C^* := \C^* (\bfx^*, \bfu^*, \bfy^*)
= \left\{
\bfd^x \in \mathbb{Re}^{n} \ \Big| \ 
A_{\I^*_+} \bfd^x  = 0, \ \
A_{\I^*_0} \bfd^x  \le 0 
\right\}.
\ee 
We will define the second-order conditions on this cone and 
it allows us to extend the classical results on second-order optimality conditions from
smooth problems \cite[Chp.~12]{nocedal2006numerical} to the nonsmooth problem (\ref{COP-Constrained}).
We first state the second-order necessary condition.

\begin{theorem} \label{Thm-SONC}
	(Second-order necessary condition)
	Let $(\bfx^*, \bfu^*)$ be a local solution of (\ref{COP-Constrained}) and let $\bfy^*$ be its KKT multiplier. Then it holds
	\[
	\bfd^\top \nabla^2 f(\bfx^*) \bfd \ge 0, \ \ \forall \ \bfd \in \C^*.
	\]
\end{theorem}

{\bf Proof.}
We recall from Prop.~\ref{Prop-First-Order-Optimality}(i) that $(\bfx^*, \bfu^*)$ is also
a local solution of the smooth problem (\ref{Smooth-P}). 
Since (\ref{Smooth-P}) has only linear constraints, which form one kind of constraint
qualification \cite[Sect.~12.6]{nocedal2006numerical},
the second-order necessary condition theorem of \cite[Thm.~12.5]{nocedal2006numerical}
applies here:
\[
[\bfd^x; \bfd^u]^\top \nabla^2_{\zeta \zeta} \cL^\# (\bfx^*, \bfu^*, \bfy^*_{\I_-^*},
\bfy^*  ) [\bfd^x; \bfd^u] \ge 0 \quad \forall \ [\bfd^x; \bfd^u] \in \C^*_s(\bfx^*, \bfu^*, \bfy^*),
\]
where $\zeta:= (\bfx, \bfu)$ and  $\cL^\#(\cdot)$ is the Lagrange function of (\ref{Smooth-P}):
\be \label{Smooth-Lagrange}
\cL^\# (\bfx, \bfy, \bfs, \bfy) 
= f(\bfx) + \langle \bfs,\; \bfu_{\I^*_-} \rangle
+ \langle \bfy, \; A\bfx + \bfb - \bfu \rangle.
\ee
Simple calculation of the Hessian of $\cL^\#$ gives the equivalent representation of
the second-order necessary condition for any $\bfd^x \in \C^*$
\be \label{Hessian-L}
[\bfd^x; \bfd^u]^\top \nabla^2_{\zeta\zeta} \cL^\# (\bfx^*, \bfu^*, \bfy^*_{\I_-^*}, \bfy^*  ) [\bfd^x; \bfd^u]  =
(\bfd^x)^\top \nabla^2 f(\bfx^*) \bfd^x  \ge 0.
\ee 
This proves our result. \hfill $\Box$

We now state a result of the second-order sufficient condition for (\ref{COP-Constrained}).

\begin{theorem} \label{Thm-SOSC}
	(Second-order sufficient condition)
	Suppose $(\bfx^*, \bfu^*, \bfy^*)$ is a {\rm P}-stationary triplet of (\ref{COP-Constrained}).
	We also suppose the second-order sufficient condition holds:
	\be \label{SOSC}
	\bfd^\top \nabla^2 f(\bfx^*) \bfd > 0 \quad \forall \ \bfd \in \C^*\setminus\{0\}.
	\ee 
	Then $(\bfx^*, \bfu^*)$ is a strict local solution of (\ref{COP-Constrained}).
	Furthermore, there exist constants $\delta^*>0$ and $c^*>0$ such that the
	following quadratic growth condition holds
	\[
	f(\bfx) + \lambda \| \bfu_+\|_0 \ge f(\bfx^*) + \lambda \| \bfu^*_+\|_0 +
	c^* \| (\bfx, \bfu) - (\bfx^*, \bfu^*) \|^2, \
	\forall \ (\bfx, \bfu) \in \N((\bfx^*, \bfu^*), \delta^*) \cap \F . 
	\]
\end{theorem}

{\bf Proof.}
From Prop.~\ref{Prop-First-Order-Optimality}(ii), we know that
$(\bfx^*, \bfu^*, \bfy^*_{\I^*_-}, \bfy^*)$ is  a KKT point of the smooth problem
(\ref{Smooth-P}). 
The first identity in (\ref{Hessian-L}) and the SOSC (\ref{SOSC}) imply 
\[
[\bfd^x; \bfd^u]^\top \nabla^2_{\zeta\zeta} \cL^\# (\bfx^*, \bfu^*, \bfy^*_{\I_-^*}, \bfy^*  ) [\bfd^x; \bfd^u] > 0 , \quad \forall \ 
[\bfd^x; \bfd^u] \in \C^*_s(\bfx^*, \bfu^*, \bfy^*) \setminus\{(0, 0)\}.
\]
This is the second-order sufficient condition for the smooth problem (\ref{Smooth-P})
at $(\bfx^*, \bfu^*)$.
It follows from \cite[Thm.~12.6]{nocedal2006numerical} that there exist
$c^*>0$ and $\delta^*$ such that
\[
f(\bfx) \ge f(\bfx^*) +  c^* \| (\bfx, \bfu) - (\bfx^*, \bfu^*) \|^2,
\quad \forall \ (\bfx, \bfu) \in \N((\bfx^*, \bfu^*), \delta^*) \cap \F^*_-.
\]
%We note that $c^*$ does not depend on $\lambda$.
Without loss of any generality (if necessary, we may reduce $\delta^*$),
we may assume
\[
\{ i \in [m] : u_i^* > 0 \} \subseteq \{ i \in [m] : u_i > 0 \}
\quad \mbox{and} \quad 
\{ i \in [m] : u_i^* < 0 \} \subseteq \{ i \in [m] : u_i < 0 \}
\]
and
\[
c^*\delta^* < \frac{\lambda} 2, \quad 
| f(\bfx) - f(\bfx^*) | \le \frac{\lambda}2, \ \ \forall \ 
(\bfx, \bfu) \in \N((\bfx^*, \bfu^*), \delta^*).
\]
Hence, it holds $\| \bfu_+\|_0 \ge \| \bfu^*_+\|_0$ and for any $\N((\bfx^*, \bfu^*), \delta^*) \cap \F^*_-$, we have
\be \label{Quadratic-growth-1}
f(\bfx) + \lambda \| \bfu_+\|_0 \ge f(\bfx^*) + \lambda \| \bfu^*_+\|_0
+  c^* \| (\bfx, \bfu) - (\bfx^*, \bfu^*) \|^2.
\ee 
Now we consider those feasible points $(\bfx, \bfu)$ of (\ref{COP-Constrained}) that are outside of $\F_-^*$. 
Since they belong to $\F \cap \N((\bfx^*, \bfu^*), \delta^*)$, there exists
an index $i_0$ such that $u_{i_0} >0$ and $u^*_{i_0} =0$. This implies
\[
\| \bfu_+\|_0 \ge \| \bfu^*_+\|_0 + 1 , \quad \forall\
(\bfx, \bfu) \in \N((\bfx^*, \bfu^*), \delta^*) \cap \F \setminus \F^*_-.
\] 
We then have for any $(\bfx, \bfu) \in \N((\bfx^*, \bfu^*), \delta^*) \cap \F \setminus \F^*_-$, 
\begin{eqnarray*}
	f(\bfx) + \lambda \| \bfu_+ \|_0 &\ge&  f(\bfx) + \lambda \| \bfu^*_+ \|_0 + \lambda \\
	&=& f(\bfx^*) + f(\bfx) - f(\bfx^*) + \lambda \| \bfu^*_+ \|_0 + \lambda \\ 
	&\ge& f(\bfx^*) - \frac{\lambda}2 + \lambda \| \bfu^*_+ \|_0 + \lambda \\
	&=& f(\bfx^*)  + \lambda \| \bfu^*_+ \|_0 + \frac{\lambda}2\\
	&\ge& f(\bfx^*)  + \lambda \| \bfu^*_+ \|_0 + c^*\delta^* \\
	&\ge& f(\bfx^*)  + \lambda \| \bfu^*_+ \|_0 +  c^* \| (\bfx, \bfu) - (\bfx^*, \bfu^*) \|^2,
\end{eqnarray*} 
which combining (\ref{Quadratic-growth-1}) yields the claimed bound, and thus $(\bfx^*, \bfu^*)$ is isolated.\hfill $\Box$

Thm.~\ref{Thm-SOSC} weakens the (local) convexity assumption in \cite{zhou2021quadratic} and (locally) strong convexity assumption in \cite{zhang2021global}  of $f(\cdot)$ to
the nonconvex case, where only the second-order sufficient condition is required.

%%%%%%%%%%%%%%%%%%%%%%%%%%%%%%%%%%%%%%%%%%%%%%%%%%%%%%%%%
\section{Subspace Newton's Method for Subproblems} \label{Section-Newton}

The main purpose of this part is to study a Newton-type method for the subproblem 
(\ref{Subproblem}) 
\[
(\bfx^{k+1}, \bfu^{k+1}) \approx \arg\min_{\bfx, \bfu} \ g_k(\bfx, \bfu) + \lambda h(\bfu)
\]
and to establish its global and local quadratic convergence.

%%%%%%%%%%%%%%%%%%%%%%%%%%%%%%%%%%%%%%%%%%%%%%%%%%%%
\subsection{Newton's method}

Let us drop the index $k$ in the subproblem and solve, without loss of any generality, 
the following problem:
\be \label{Problem-G}
\min_{\bfx, \bfu} \ G(\bfx, \bfu) := g(\bfx, \bfu) + \lambda \| \bfu_+\|_0 ,
\ee 
where for given $\bfp \in \mathbb{R}^m$ and $\bfq\in \mathbb{R}^n$,
\[
g(\bfx, \bfu) := f(\bfx) + \langle \bfp,\; A\bfx - \bfu + \bfb \rangle
+ \frac{\rho}2 \| A\bfx - \bfu + \bfb \|^2 + \frac{\mu}2 \| \bfx - \bfq\|^2.
\]
Problem (\ref{Problem-G}) in general has no closed-form solution and has to be
solved iteratively. We assume the current iterate is $(\bfx^j, \bfu^j)$. 
We would like to update it via Newton's method. 
The difficulty is that we have a nonsmooth function $h(\bfu)$ in $G(\bfx, \bfu)$.
We explain our strategy below.

Let $\Gamma_j$ be the index set corresponding to the zero components of $\bfu^j$. 
We approximate $g(\bfx, \bfu)$ by its second-order Taylor expansion at
$(\bfx^j, \bfu^j)$ restricted to the
subspace $\bfu_{\Gamma_j}=0$:
\begin{eqnarray*}
	g(\bfx, \bfu) |_{\bfu_{\Gamma_j} = 0} 
	&\approx g(\bfx^j, \bfu^j) + \langle \nabla_{\bfx} g(\bfx^j, \bfu^j), \; \bfx - \bfx^j \rangle
	+ \langle (\nabla_{\bfu} g(\bfx^j, \bfu^j)  )_{\OG_j}, \;
	(\bfu - \bfu^j)_{\OG_j} \rangle \\
	& + \frac 12 \left[
	\begin{array}{c}
		\bfx - \bfx^j \\
		(\bfu - \bfu^j)_{\OG_j}
	\end{array} 
	\right]^\top
	\underbrace{\left[
		\begin{array}{cc}
			\nabla^2 f(\bfx^j) + \rho A^\top A + \mu I, & \ \rho (A_{\OG_j})^\top \\ [1ex]
			\rho A_{\OG_j}, & \ \rho I
		\end{array}
		\right]}_{=: H^j}
	\left[
	\begin{array}{c}
		\bfx - \bfx^j \\
		(\bfu - \bfu^j)_{\OG_j}
	\end{array} 
	\right]
\end{eqnarray*} 
Assuming the matrix $(\nabla^2 f(\bfx^j) + \mu I)$ is positive definite, then
matrix $H^j$ is also positive definite by the theorem of Schur complement.
Minimizing the quadratic approximation is equivalent to solving the Newton equation:
\[
H^j \left[
\begin{array}{c}
	\bfx - \bfx^j \\ [0.5ex]
	(\bfu - \bfu^j)_{\OG_j}
\end{array} 
\right] = - \left[
\begin{array}{c}
	\nabla_{\bfx} g(\bfx^j, \bfu^j) \\ [0.5ex]
	(\nabla_{\bfu} g(\bfx^j, \bfu^j)  )_{\OG_j}
\end{array} 
\right].
\]
To accommodate the use of first-order methods \cite{beck2017first} to problem
(\ref{Problem-G}) when the Newton equation is not solvable or the Newton direction is
not a descent direction, 
we apply a gradient step at $(\bfx^j, \bfu^j)$ along the
subspace $\bfu_{\Gamma_i} =0$ to generate a new point denoted by
$(\bfx^{j+1/2}, \bfu^{j+1/2})$ and we then 
try the Newton step at this new point.
The gradient step serves as a globalization strategy (i.e., global convergence)
of the Newton method.
The difficulty in this setting is that the subspace keeps changing each iteration and
the gradient used is truncated information. 
Hence, the classical convergence of gradient methods cannot be directly used here.

To simplify the convergence analysis and to make the Newton-step well defined so that we
can analyze it, we assume that the function $f(\bfx)$ is sufficiently smooth. 
In particular, we assume the following.

\begin{assumption} \label{Assumption-f}
	\begin{itemize}
		\item[(i)] The function $f$ is twice continuously differentiable and it is
		$\sigma_f$-weakly convex, i.e., $f(\bfx) + (\sigma_f/2) \| \bfx\|^2$ is convex.
		
		\item[(ii)] The gradient function $\nabla f(\bfx)$ is Lipschitz continuous with modulus $\ell_f$:
		\[
		\| \nabla f(\bfx) - \nabla f(\bfy)\| \le \ell_f \| \bfx-\bfy \|, \quad \forall \
		\bfx, \bfy \in \mathbb{R}^n .
		\]
		
		\item[(iii)] The Hessian function $\nabla^2 f(\bfx)$ is Lipschitz continuous with modulus $L_f$:
		\[
		\| \nabla^2 f(\bfx) - \nabla^2 f(\bfy)\| \le L_f \| \bfx-\bfy \|, \quad \forall \
		\bfx, \bfy \in \mathbb{R}^n .
		\]
		
	\end{itemize}
\end{assumption}

\begin{remark} \label{Remark-g}
	A consequence of this assumption in terms of $g$ is that (i) $g$ is strongly convex
	provided that $\mu > \sigma_f$ (denote $\sigma_g := \mu- \sigma_f$), and 
	(ii) the gradient function $\nabla g(\bfx, \bfu)$ and the Hessian 
	$\nabla^2 g(\bfx, \bfu)$ are both Lipschitz continuous with their modulus respectively 
	denoted by $\ell_g$ and $L_g$.
	Although the Lipschitzian properties if $f$ are defined over the space $\mathbb{R}^n$,
	they can be restricted to a bounded region because we will prove that the generated
	sequence $\{\bfx^k\}$ will remain bounded.
\end{remark}

The globalized subspace Newton method is described in Alg.~\ref{Alg-GSN}.

\begin{algorithm}[htbp] 
	\caption{GSNewton: Gradient-Subspace-Newton method} \label{Alg-GSN}
	\begin{algorithmic}
		
		\STATE{Initialization: Take $\alpha \in (0,1/\ell_g)$ and $t \in (0,2/\ell_g)$, and then start with $(\bfx^0, \bfu^0) \in \mathbb{R}^n \times \mathbb{R}^m$.}
		\FOR{$j=0,1,\cdots$}
		
		\STATE{\textbf{1. Identification step:} } 
		Denote $\zeta^j := (\bfx^j, \bfu^j)$, $\bfy^j := - \nabla_{\bfu} g(\zeta^j)$ and define the active set  $\Gamma_j$ by
		\be \label{Gamma_j}
		\Gamma_j := \left\{
		i \in [m] \ \Big| \ [\bfu^j + \alpha \bfy^j]_i \in [0, \sqrt{2\alpha \lambda} )
		\right\}.
		\ee 
		\STATE{\textbf{2. Gradient step (the half-step): } }	
		Compute $\zeta^{j+1/2} := (\bfx^{j+1/2},\; \bfu^{j+1/2} ))$ by
		\be \label{Gradient-Step}
		\left\{
		\begin{array}{l}
			u^{j+1/2}_i = \left\{
			\begin{array}{ll}
				0, & \ \mbox{if} \ i \in \Gamma_j \\ 
				\left[ \bfu^j + \alpha \bfy^j \right]_i, & \ \mbox{if} \ i \in \OG_j, 
			\end{array} 
			\right . \quad i \in [m] \\ [3ex]
			%\quad \mbox{and} \quad
			\bfx^{j+1/2} = \bfx^j - t \nabla_{\bfx} g (\bfx^j, \bfu^{j+1/2}).
		\end{array} 
		\right. 
		\ee 
		
		%	     
		%	     
		%		\begin{align} 
		%			&u^{j+1/2}_{\Gamma_{j}^c} = (u^{j} + \alpha z^{j})_{\Gamma_{j}^c},~	u^{j+1/2}_{\Gamma_{j}} = 0, \label{u-step} \\
		%			&w^{j+1/2} = w^j -  t\nabla_w g(w^j,u^{j+1/2}) . \label{w-step}	
		%		\end{align}
		\STATE{\textbf{3. Newton step:} }  Compute the Newton step $\wz^{j+1} := (\widetilde{\bfx}^{j+1}, \widetilde{\bfu}^{j+1})$, which satisfies the following
		truncated Newton equation in $(\bfx, \bfu)$:
		\be \label{Newton-Eq}
		H^{j+1/2}
		\left[
		\begin{array}{c}
			\bfx - \bfx^{j+1/2} \\ [0.5ex]
			(\bfu - \bfu^{j+1/2})_{\OG_j}
		\end{array} 
		\right] = - \left[
		\begin{array}{c}
			\nabla_{\bfx} g(\zeta^{j+1/2}) \\ [0.5ex]
			(\nabla_{\bfu} g(\zeta^{j+1/2})  )_{\OG_j}
		\end{array} 
		\right] \quad \mbox{and}
		\quad
		\bfu_{\Gamma_j} = 0 .
		\ee 
		where
		\[
		H^{j+1/2} := \left[
		\begin{array}{cc}
			\nabla^2 f(\bfx^{j+1/2}) + \rho A^\top A + \mu I, & \ \rho (A_{\OG_j})^\top \\ [1ex]
			\rho A_{\OG_j}, & \ \rho I
		\end{array}
		\right] 
		\]
		
		\STATE{\textbf{4. Update step:} }	Update $\zeta^j$ either by the Newton step or
		the gradient step as follows:
		\be \label{Newton-Condition}
		\zeta^{j+1} = \left\{
		\begin{array}{ll}
			\wz^{j+1} , & \ \mbox{if} \ G ( \zeta^{j+1/2} ) - G ( \widetilde{\zeta}^{j+1}) \geq (\sigma_g/4) \| \widetilde{\zeta}^{j+1} - \zeta^{j+1/2} \|^2 \\ [1ex]
			\zeta^{j+1/2} , & \ \mbox{otherwise}
		\end{array} 
		\right .
		\ee 
		
		\ENDFOR
	\end{algorithmic}
\end{algorithm}

We make some comments on Alg.~\ref{Alg-GSN}.

\begin{itemize}
	\item[(i)] The active set $\Gamma_j$ in (\ref{Gamma_j}) is dynamically updated.
	It involves the gradient information at $\bfu^j$. 
	Hence, it does not in general coincide with the true active set of indices where $u^j_i =0$. However, it does coincide with the true active set at $\bfu^{j+1/2}$
	due to (\ref{Gradient-Step}):
	\be \label{Proximal-Update}
	\Gamma_j = \left\{
	i \in [m] \; \Big| \ [ \bfu^{j+1/2}]_i = 0
	\right\}
	\quad \mbox{and} \quad
	\bfu^{j+1/2} \in \Prox_{\alpha\lambda h(\cdot)}( \bfu^{j} + \alpha \bfy^{j} ).
	\ee 
	The second relationship above follows (\ref{Gamma_j}) and (\ref{prox0+}).
	
	\item[(ii)] Because $\Gamma_j$ represents the true active set at the half step point $\bfu^{j+1/2}$, the Newton matrix $H^{j+1/2}$ in (\ref{Newton-Eq}) 
	is the true truncated Hessian along the subspace $\bfu_{\Gamma_j} =0$ at
	$\zeta^{j+1/2}$. This is why we are able to establish quadratic convergence 
	under reasonable conditions. 
	By applying the Schur-complement theorem to the matrix $H^{j+1/2}$, solving
	(\ref{Newton-Eq}) is equivalent to 
	solving the following type of linear equation in $\bfd_{\bfx}$ with a given $\bfb_{\bfx} \in \mathbb{R}^n$
	\be \label{Reduced-Newton-Eq}
	\Big( \nabla^2 f(\bfx) + \mu I + \rho (A_{\Gamma})^\top A_{\Gamma} \Big) \bfd_{\bfx} =
	\bfb_{\bfx}.
	\ee
	Suppose $(\nabla^2 f(\bfx) + \mu I)$ is nonsingular, the matrix in the equation is
	a low-rank correction of $(\nabla^2 f(\bfx) + \mu I)$.
	We may apply Sherman-Morrison-Woodbury formula \cite[A.28]{nocedal2006numerical} to solve this equation.
	In general, the computational complexity is approximately $O(n^2 \max\{n, |\Gamma| \}$). 
	This is fine for small-sized problems. 
	When $n$ is	large, the computational cost is too high.
	In this case, we would 
	recommend the use of the gradient step in Alg.~\ref{Alg-GSN}  
	and only use the Newton step in the final stage of the algorithm.
	Fortunately, for many real applications, such as SVM,  their
	functions $f(\bfx)$ are separable and have small non-overlapping block structures.  This implies that
	the Hessian matrix 
	$\nabla^2 f (\bfx)$ 
	is of diagonal blocks and is invertible. The worst-case
	computational complexity can be reduced to 
	$O(|\Gamma|^2 \max\{ n, |\Gamma| \})$.
	For SVM, $\Gamma_j$ coincides with the indices of support vectors that take
	a relatively small portion of the total samples. Hence,
	$|\Gamma|$  can be on a small scale
	and computation of the Newton direction can be very cheap.
	
\end{itemize}

%%%%%%%%%%%%%%%%%%%%%%%%%%%%%%%%%%%%%%%%%%%%%%%%%
\subsection{Global convergence}

Since we are not assuming any convexity, the generated sequence can only converge to a stationary point.
Suppose $(\widehat{\bfx}, \widehat{\bfu})$ is a ${\rm P}$-stationary point with ${\rm P}$-stationary multiplier $\widehat{\bfy}$ and $\alpha >0$ for problem (\ref{Problem-G})
so that 
\be \label{P^k-stat}
\nabla_\bfx g(\widehat{\bfx},\widehat{\bfu}) = 0, \quad
\widehat{\bfy} = -\nabla_\bfu g(\widehat{\bfx},\widehat{\bfu}), \quad
\widehat{\bfu} \in {\rm Prox_{\alpha\lambda h(\cdot )}} ( \widehat{\bfu} + \alpha \widehat{\bfy} ).
\ee 
The global convergence below states that (i) the objective sequence is non-increasing, (ii) the whole sequence of iterates converges to a stationary point, and
(iii) the 0/1-function $h(\cdot)$ is
continuous at the limit along the generated iterates. 

\begin{theorem}[Global Convergence of GSNewton] \label{Thm-Global} 
	Suppose Assumption~\ref{Assumption-f}~(i)-(ii)  holds, 
	and $\mu > \sigma_f$. 
	Let $\sigma_g = \mu - \sigma_f$.
	Let $\ell_g$ be the modulus assumed in Remark~\ref{Remark-g}. 
	Let $\{ \zeta^{j} \}_{j \in \mathbb{N}}$ and $\{\bfy^j\}_{j \in \mathbb{N}}$ be 
	respectively the iterate sequence and the multiplier (i.e., gradient) sequence generated by 
	Alg.~\ref{Alg-GSN}.
	The following hold.
	
	\begin{itemize}
		
		\item[(i)] The objective value sequence $\{ G(\zeta^j) \}_{j \in \mathbb{N}}$
		is non-increasing and  satisfies
		\begin{align} \label{G-Decreasing}
			G ( \zeta^{j} ) - G ( \zeta^{j+1}) 
			\ge \tau \| \zeta^{j+1/2} - \zeta^j \|^2 
			+ \frac{\sigma_g}4 \| \zeta^{j+1} - \zeta^{j+1/2} \|^2,
		\end{align}
		where $\tau := \min\{1/(2\alpha) - \ell_g/2, \; \ell_g/2 - 1/t  \} > 0$ 
		Consequently, it holds 
		\begin{equation}\label{msn_j+1-j}
			\lim_{j \to \infty} \| \zeta^{j+1/2} - \zeta^{j} \| = 0
			\quad \mbox{and} \quad
			\lim_{j \to \infty} \| \zeta^{j+1} - \zeta^{j} \| = 0.
		\end{equation}
		
		\item[(ii)] The sequence $\{\zeta^{j}\}_{j \in \mathbb{N}}$ and $\{\bfy^{j}\}_{j \in \mathbb{N}}$ converge respectively to ${\rm P}$-stationary point $\widehat{\zeta}:=(\widehat{\bfx}, \widehat{\bfu})$ and the corresponding ${\rm P}$-stationary multiplier $\widehat{\bfy}=-\nabla_u g(\widehat{\zeta})$ of problem
		(\ref{Subproblem}).  
		
		\item[(iii)]  We have
		\be \label{lim-G}
		\lim_{j \to \infty} \| \bfu^{j+1}_+ \|_0 = \lim_{j \to \infty} \| \bfu^{j+1/2}_+ \|_0 = \| \widehat{\bfu}_+ \|_0  \quad \mbox{and} \quad 
		\lim_{j \to \infty} G (\zeta^{j}) = G (\widehat{\zeta})
		\ee 		
	\end{itemize}
\end{theorem}
%%%%%%%%%%%%%%%%%%%%%%%%%%%%%%%%%%%%%%%%%%%%%%%%%%%%%%%%%%%%%%

{\bf Proof}
(i) It follows from (\ref{Proximal-Update}) and the definition of proximal operator that
\begin{align}\notag
	\frac{1}{2} \| \bfu^{j+1/2} - (\bfu^{j} + \alpha \bfy^j) \|^2 + \alpha \lambda \| \bfu^{j+1/2}_+ \|_0 \leq \frac{1}{2} \| \bfu^{j} - (\bfu^{j} + \alpha \bfy^j) \|^2 + \alpha \lambda \| \bfu^{j}_+ \|_0, 
\end{align}
and thus we have
\begin{align} \notag
	\lambda \| \bfu^{j+1/2}_+ \|_0 - \lambda \| \bfu^j_+ \|_0  \leq \langle \bfy^j, \bfu^{j+1/2} - \bfu^j \rangle -\frac{1}{2\alpha} \| \bfu^{j+1/2} - \bfu^j \|^2.
\end{align}
By Lipschitz continuity of $\nabla_{\bfu} g(\bfx^j, \cdot)$ and the Descent Lemma, the following estimation holds
\begin{align} \notag
	g(\bfx^{j}, \bfu^{j+1/2}) - g(\bfx^j, \bfu^j) \leq \langle -\bfy^j, \bfu^{j+1/2} - \bfu^j \rangle + \frac{\ell_g}{2} \| \bfu^{j+1/2} - \bfu^j \|^2.
\end{align} 
Adding the above two inequalities yields
\begin{align} \label{u_descent}
	G(\zeta^j) - G(\bfx^{j}, \bfu^{j+1/2}) 
	\ge ( \frac{1}{2\alpha} - \frac{\ell_g}{2} ) \| \bfu^{j+1/2} -\bfu^j \|^2
	\ge \tau \| \bfu^{j+1/2} -\bfu^j \|^2.
\end{align}
By Lipschitz continuity of $\nabla_\bfx g(\cdot, \bfu^{j+1/2})$ and the Descent Lemma, we obtain
\begin{align} \notag
	g(\bfx^j, \bfu^{j+1/2}) \geq & g(\bfx^{j+1/2}, \bfu^{j+1/2}) - \langle \nabla_{\bfx} g(\bfx^j, \bfu^{j+1/2}), \bfx^{j+1/2} - \bfx^j \rangle - \frac{\ell_g}{2} \| \bfx^{j+1/2} - \bfx^j \|^2 \notag \\
	\mathop{=}\limits^{(\ref{Gradient-Step})}& g(\bfx^{j+1/2},\bfu^{j+1/2}) + ({\ell_g}/{2} - {1}/{t}) \| \bfx^{j+1/2} - \bfx^j \|^2  \notag \\ 
	= \ & g(\bfx^{j+1/2},\bfu^{j+1/2}) + \tau \| \bfx^{j+1/2} - \bfx^j \|^2. \notag
\end{align}
Adding $\| \bfu^{j+1/2}_+ \|_0$ on both sides of above inequality implies 
\begin{align} 
	&G(\bfx^j, \bfu^{j+1/2}) - G( \zeta^{j+1/2} ) \mathop{\geq} \tau \| \bfx^{j+1/2} - \bfx^j \|^2.  \label{w_descent}
\end{align}
%Adding above inequality and (\ref{u_descent}), we can obtain
%\begin{align} \notag
%	G(\zeta^j) - G( \zeta^{j+1/2} ) \geq \tau_w \| w^{j+1/2} - w^j \|^2 + \tau_u \| u^{j+1/2} - u^j \|^2.
%\end{align}
If Newton step is accepted, we then have from (\ref{Newton-Condition})
\begin{align}\label{newt_des}
	G ( \zeta^{j+1/2} ) - G ( \zeta^{j+1}) \geq (\sigma_g/4) \| \zeta^{j+1} - \zeta^{j+1/2} \|^2.
\end{align}
Adding (\ref{u_descent}), (\ref{w_descent}) and (\ref{newt_des}) together, 
we obtain (\ref{G-Decreasing}).
If Newton step is not accepted, then (\ref{G-Decreasing}) reduces to the case of the gradient step due to $\zeta^{j+1} = \zeta^{j+1/2}$.

The strongly convexity of $g$ indicates that it is bounded from below. 
Due to $\| (\cdot)_+ \|_0 \geq 0$, $G$ is also bounded from below. 
The non-increasing property of  $\{ G(\zeta^j) \}_{j \in \mathbb{N}}$ 
in (\ref{G-Decreasing}) implies that there exists a real number $\widehat{G}$ such that $\lim_{j \to \infty} G(\zeta^j) =  \widehat{G}$. 
By taking limit on both sides of (\ref{G-Decreasing}), we get
$\lim_{j \to \infty} \| \zeta^{j+1/2} - \zeta^j \| = 0$ and $\lim_{j \to \infty} \| \zeta^{j+1} - \zeta^{j+1/2} \| = 0$. Hence, $\lim_{j \to \infty} \| \zeta^{j+1} - \zeta^{j} \| = 0$ also holds.

(ii)
The Assumption~\ref{Assumption-f} implies that the function $g(\cdot)$ is coercive,
and so is $G(\cdot)$.
We established in (i) that  $\{ G(\zeta^j) \}_{j \in \mathbb{N}}$ is  non-increasing. Hence, the sequence $\{ \zeta^j \}_{j \in \mathbb{N}}$ is bounded by the coerciveness of $G(\cdot)$. 
We establish (ii) in two steps. 
First, we prove every accumulation point is a P-stationary point of (\ref{Problem-G}).
We then prove there are only finite such stationary points. 
The convergence of the whole sequence to a stationary point will follow.

{\bf Step 1:} Suppose that $\widehat{\zeta}$ is an accumulation point of $\{ \zeta^j \}_{j \in \mathbb{N}}$, then there exists an index set $J \subseteq \mathbb{N}$ such that
\begin{align} \label{sub_lim}
	\lim_{j \in J, j \to \infty} \zeta^j = \widehat{\zeta} = [\widehat{\bfx}, \widehat{\bfu} ] \quad \mbox{and} \quad
	\lim_{j \in J, j \to \infty} \bfy^j = - \nabla_\bfu g(\widehat{\zeta}) = \widehat{\bfy}. 
\end{align}
From (\ref{sub_lim}) and (\ref{msn_j+1-j}), we have
\begin{align}
	&\lim_{j \in J, j \to \infty} \zeta^{j+1/2} = \lim_{j \in J, j \to \infty} [(\zeta^{j+1/2} - \zeta^j) + \zeta^j] = \widehat{\zeta}, \label{1/2converge} \\ 
	&\lim_{j \in J, j \to \infty} \bfu^{j} + \alpha \bfy^j = \widehat{\bfu} + \alpha \widehat{\bfy}.  \label{lim_uz}
\end{align}
Since \eqref{1/2converge}, \eqref{lim_uz} and \eqref{Proximal-Update} hold, using the proximal behavior theorem \cite[Theorem 1.25]{RockWets98}), we have
\begin{align} \notag
	\widehat{\bfu} \in {\rm Prox}_{\alpha\lambda h(\cdot)}( \widehat{\bfu} + \alpha \widehat{\bfy} ).
\end{align} 
The gradient step (\ref{Gradient-Step}) can be rewritten as 
$\nabla_\bfx g(\bfx^j, \bfu^{j+1/2}) = -(\bfx^{j+1/2} - \bfx^j)/t$.
Since (\ref{sub_lim}) and (\ref{1/2converge}) hold, and $g$ is continuously differentiable, then taking limit as $j \to \infty$ for $j \in J$ on both sides of this equality yields 
$ %\begin{align} \notag
\nabla_\bfx g(\widehat{\bfx},\widehat{\bfu}) = 0.
$ %\end{align} 
Therefore, $(\widehat{\bfx}, \widehat{\bfu})$ together with $\widehat{\bfy}$ satisfies (\ref{P^k-stat}) for $\alpha \in (0, 1/\ell_g)$, and hence it is a P-stationary 
point. 

{\bf Step 2:} We now prove that there are finite proximal triplets.
To proceed, we denote $\widehat{\I}_-:= \{ i \in [n]\; | \; \widehat{u}_i \leq 0\}$ 
and consider the following convex programming with linear inequality constraints.
\begin{equation}\label{sub-conprog}
	\min\limits_{\bfx \in \mathbb{R}^n, \bfu \in \mathbb{R}^m} g(\bfx, \bfu),
	\quad \mbox{s.t.} \quad \bfu_{\widehat{\I}_-} \leq 0.	
\end{equation}
Since $g$ is strongly convex, (\ref{sub-conprog}) has the unique global minimizer and its KKT condition can be represented as
\begin{equation} \label{sub_kkt}
	\left\{  \begin{aligned}
		& \nabla_\bfx g(\bfx, \bfu) = 0, \\
		& \bfy_{\widehat{\I}_-} = -[\nabla_\bfu g(\bfx,\bfu)]_{\widehat{\I}_-},~ [\nabla_\bfu g(\bfx, \bfu)]_i = 0,  \ i \not\in \widehat{\I}_- \\
		& \bfy_{\widehat{\I}_-} \geq 0 ,~ \bfu_{\widehat{\I}_-} \leq 0,~\langle \bfy_{\widehat{\I}_-}, \bfu_{\widehat{\I}_-}  \rangle = 0,
	\end{aligned}  \right.
\end{equation}
where $\bfy_{\widehat{\I}_-} \in \mathbb{R}^{| \widehat{\I}_- |}$ is the Lagrange multiplier. From (\ref{P^k-stat}) and (\ref{prox0+}), we know that a ${\rm P}$-stationary triplet satisfies (\ref{sub_kkt}). 
Since (\ref{sub-conprog}) is a strongly convex problem, $(\widehat{\bfx}, \widehat{\bfu})$ is the unique global minimizer of (\ref{sub-conprog}). Since there are only a finite number of choices of $\widehat{\I}_-$ (no more than $2^n$), the number of ${\rm P}$-stationary points of (\ref{Problem-G}) is finite. 

Finally, because we have proved that each accumulation point of $\{ \zeta^{j} \}_{j \in \mathbb{N}}$ is a ${\rm P}$-stationary point, the number of accumulation points is finite, and thus each accumulation point of $\{ \zeta^j \}_{j \in \mathbb{N}}$ is isolated. Considering that (\ref{msn_j+1-j}) holds, it follows from \cite[Prop.~7]{kanzow1999qp} that the whole sequence $\{ \zeta^j \}_{j \in \mathbb{N}}$ must converge to $\widehat{\zeta}$. Moreover, we have $ \lim_{j \to \infty} \bfy^{j} = \lim_{j \to \infty} -\nabla_\bfu g(\zeta^j) = -\nabla_\bfu g(\widehat{\zeta})$, which means that $\{ \bfy^j \}_{j \in \mathbb{N}}$ converges to the corresponding ${\rm P}$-stationary multiplier.

(iii) From (\ref{Proximal-Update}) and the definition of the proximal operator, we have
\begin{align}
	\frac{1}{2} \| \bfu^{j+1/2} - \bfu^{j} + \alpha \bfy^j \|^2 + \alpha \lambda \| \bfu^{j+1/2}_+ \|_0 \leq \frac{1}{2} \| \widehat{\bfu} - \bfu^{j} + \alpha \bfy^j \|^2 + \alpha \lambda \| \widehat{\bfu}_+ \|_0 . \notag
\end{align}
and then taking the upper limit in the above inequality, we obtain
\begin{align}
	\limsup_{j \to \infty} \| \bfu^{j+1/2}_+ \|_0 \leq \| \widehat{\bfu}_+ \|_0. \label{limsup_j1/2}
\end{align}
This formula together with the lower semi-continuity (lsc) of $\|(\cdot)_+\|_0$ implies that $\lim_{j \to \infty} \| \bfu^{j+1/2}_+ \|_0 = \| \widehat{\bfu}_+ \|_0$.
If Newton  step is accepted, it follows from
(\ref{newt_des}) that   
\begin{align} 
	g(\zeta^{j+1}) + \lambda \| \bfu^{j+1}_+ \|_0 + \frac{\sigma_g}{4} \| \zeta^{j+1} - \zeta^{j+1/2} \|^2 \leq  g(\zeta^{j+1/2}) + \lambda \| \bfu^{j+1/2}_+ \|_0. \notag
\end{align}
Due to (\ref{msn_j+1-j}), passing the upper limit in the above inequality as $j \to \infty$, we can obtain
\begin{align}
	\limsup_{j \to \infty} \| \bfu^{j+1}_+ \|_0 \leq \limsup_{j \to \infty} \| \bfu^{j+1/2}_+ \|_0 \mathop{=}\limits^{(\ref{limsup_j1/2})} \| \widehat{\bfu}_+ \|_0. \notag
\end{align}
Then combining the lsc property of $ \|(\cdot)_+\|_0 $ implies $\lim_{j \to \infty}\| \bfu^{j+1}_+ \|_0 =  \| \widehat{\bfu}_+ \|_0$.
\hfill $\Box$

%%%%%%%%%%%%%%%%%%%%%%%%%%%%%%%%%%%%%%%%%%%%%%%%%%%%%%%%%%%%%%
\subsection{Quadratic convergence}

Having established the global convergence of Alg.~\ref{Alg-GSN}, we would like to
investigate the conditions that guarantee a local fast convergence rate. 
Let us introduce some notation that would make the Newton equation more conducive to
analysis. 
Given an index set $T \subseteq [m]$, we denote 
\begin{eqnarray*}
	&&	\zeta_{:T}:= [\bfx ; \bfu_T], \ \nabla_T g(\zeta):= (\nabla_\bfu g(\zeta))_T, \ \nabla_{:T} g (\zeta):= [ \nabla_\bfx g ( \zeta ); \ (\nabla_\bfu g (\zeta) )_T ],
	\\ [1ex]
	&& \nabla^2_{:T} g(\zeta) = 
	\left[ \begin{array}{cc}
		\nabla^2_\bfx g (\zeta) & (\nabla^2_{\bfx,\bfu} g(\zeta))_{:T} \\ [1ex]
		(\nabla^2_{\bfu, \bfx} g (\zeta))_{T:} &\nabla^2_\bfu g (\zeta)
	\end{array} \right]. \notag
\end{eqnarray*}
With those notation, the Newton equation (\ref{Newton-Eq})
with its solution $\widetilde{\zeta}$ can be written as
\be \label{Newton-Eq-2}
\nabla^2_{:\OG_j} g(\zeta^{j+1/2} ) (  \widetilde{\zeta} - \zeta^{j+1/2} )_{: \OG_j}
= - \nabla_{:\OG_j} g(\zeta^{j+1/2}).
\ee
Suppose $(\widehat{\bfx}, \widehat{\bfu})$ is a P-stationary point with the multiplier
$\widehat{\bfy}$ that are obtained by Alg.~\ref{Alg-GSN}.
It follows from (\ref{Prox-eq}) that the complementarity condition holds: 
$\widehat{u}_i \widehat{y}_i =0$ for $i \in [m]$.
We say the {\em strict
	complementarity} condition holds at $(\widehat{\bfu}, \widehat{\bfy})$ if
$
\widehat{u}_i + \widehat{y}_i  \not= 0,  \forall \ i \in [m].
$
Under this condition, the identification step in Alg.~\ref{Alg-GSN} can eventually correctly identify the active set at $(\widehat{\bfx}, \widehat{\bfu})$. 
Consequently, the Newton step will be eventually accepted, leading to quadratic convergence. We state this key technical result as a lemma.

\begin{lemma} \label{quadra_lemma}
	Suppose Assumption~\ref{Assumption-f} hold.	
	Let $\{ \zeta^{j} \}_{j \in \mathbb{N}}$ be a sequence generated by Alg.~\ref{Alg-GSN} converging to a ${\rm P}$-stationary point 
	$\widehat{\zeta} = (\widehat{\bfx}, \widehat{\bfu}) $ with ${\rm P}$-stationary multiplier $\widehat{\bfy}$. 
	We further assume that the strict complementarity condition holds at $(\widehat{\bfu}, \widehat{\bfy})$, i.e., $\widehat{u}_i + \widehat{y}_i \not=0$ for all $i \in [m]$.
	Then there exists a sufficiently large index $\widehat{j} \in \mathbb{N}$, 
	such that the following hold:
	
	\begin{itemize}
		\item[(i)] The index set $\Gamma_j$ correctly identifies the active set at
		$(\widehat{\bfx}, \widehat{\bfu})$:
		\begin{equation} \label{set_iden}
			\Gamma_{j} = \widehat{\I}_0 := \left\{ i \in [m] \; | \ 
			[\widehat{\bfu}]_i = 0 \right\} , \quad \forall \ j \ge \widehat{j}.
		\end{equation}
		
		\item[(ii)]  The Newton step is always accepted at each iterate $({\bfx}^j, {\bfu}^j)$ for all $j \ge \widehat{j}$.
		
	\end{itemize} 
\end{lemma}

{\bf Proof.}
(i) Since $\widehat{\zeta}$ is a ${\rm P}$-stationary point with ${\rm P}$-stationary multiplier $\widehat{z}$ of problem (\ref{Problem-G}), it follows from 
(\ref{Prox-eq}) that under the strict complementarity condition we have
\be \label{uzval}
\widehat{u}_i = 0 \quad \mbox{if and only if} \quad
\widehat{y}_i \in (0, \sqrt{2\lambda/\alpha}].
\ee 
%\textcolor{blue}
{At $\widehat{\zeta}$, we define two index sets
	\begin{align} \notag
		\widehat{\Gamma}_1:= \{ i \in [m] \ \Big| \ \widehat{u}_i = 0, 0< \widehat{y}_i < \sqrt{2\lambda/\alpha} \},~~\widehat{\Gamma}_2:= \{ i \in [m] \ \Big| \ \widehat{u}_i = 0,  \widehat{y}_i = \sqrt{2\lambda/\alpha} \},
	\end{align}
	and then the relationship $\widehat{\I}_0 = \widehat{\Gamma}_1 \cup \widehat{\Gamma}_2$ can be obtained, under the strict complementarity condition, using the definition of $\Prox_{\lambda \alpha h(\cdot)} (\cdot)$ in  (\ref{prox0+}) and the fact
	$\widehat{\bfu} \in \Prox_{\alpha \lambda h(\cdot)}( \widehat{\bfu} + \alpha \widehat{\bfy}  )$.}

%\textcolor{blue}
{We now prove ${\Gamma}_j \supseteq \widehat{\Gamma}_1$ and ${\Gamma}_j \supseteq \widehat{\Gamma}_2$ whenever $j$ is sufficiently large, and thus ${\Gamma}_j \supseteq \widehat{\I}_0$. Firstly, note that $\widehat{\Gamma}_1 = \{ i \in [m] \ \Big| \ 0 < \widehat{u}_i + \alpha \widehat{y}_i < \sqrt{2\lambda\alpha} \}$, then from the definition of $\Gamma_j$ and $\lim_{j \to \infty} (\bfu^j,\bfy^j) = (\widehat{\bfu},\widehat{\bfy})$, ${\Gamma}_j \supseteq \widehat{\Gamma}_1$ holds when $j$ is large enough. After that, let us prove ${\Gamma}_j \supseteq \widehat{\Gamma}_2$.
	Suppose there exists an infinite set $J$ such that $\Gamma_{j} \not\supseteq \widehat{\Gamma}_2$ for any $j \in J$. 
	Without loss of generality, there exists a fixed index $\widehat{i}$ (if necessary, we may draw a subsequence from $J$) such that
	$\widehat{i} \in \widehat{\Gamma}_2$ but not in $\Gamma_j$. That is
	$\widehat{y}_{\widehat{i}} \not=0$. Since $\widehat{i} \not\in \Gamma_j$,
	by the update rule in (\ref{Gradient-Step}), we know
	\[
	{y}_{\widehat{i}}^j = \frac{ u^{j+1/2}_{\widehat{i}} - u^{j}_{\widehat{i}} }{\alpha}.
	\]
	Taking limit on both sides on $j \in J$ and using the convergence of $\zeta^j$ and $\zeta^{j+1/2}$ to $\widehat{\zeta}$, we get
	$
	\widehat{y}_{\widehat{i}} = 0,
	$
	which contradicts the fact $\widehat{y}_{\widehat{i}} \not=0$. 
	Hence, ${\Gamma}_j \supseteq \widehat{\Gamma}_2$ for sufficiently large $j$.}

We now prove the converse part $\Gamma_{j} \subseteq \widehat{\I}_0$ when $j$ is large enough. 
Let $\widehat{\delta} := \{ |\widehat{u}_i|/2 \; |\; \widehat{u}_i \neq 0 \}$. 
Since $\lim_{j \to \infty} \bfu^{j+1/2} = \widehat{\bfu}$, 
$\bfu^{j+1/2} \in \mathcal{N}( \widehat{\bfu}, \widehat{\delta} )$ holds for all $j$ 
sufficiently large.
Taking any index $i_0 \not\in \widehat{\Gamma}$, we have
\begin{align} \notag
	| u^{j+1/2}_{i_0} | \geq | \widehat{u}_{i_0} | - | u^{j+1/2}_{i_0} - \widehat{u}_{i_0} |  > | \widehat{u}_{i_0} | - \widehat{\delta} > 0,
\end{align}
which means that $u^{j+1/2}_{i_0} \neq 0$.
Note that from (\ref{Proximal-Update}), $\Gamma_j \supseteq \{ i \in [n]: u^{j+1/2}_i = 0 \} $ holds. This means $i_0 \not\in \Gamma_j$, which proves
$\Gamma_{j} \subseteq \widehat{\I}_0$ and establishes (i).

(ii) From now on, we only consider those iterates $j \ge \widehat{j}$ so that $\Gamma_j$ identifies $\widehat{\I}_0$ and we focus on the Newton iterate $\widetilde{\zeta}^{j+1}$. Our aim is to show $\zeta^{j+1} = \widetilde{\zeta}^{j+1}$
in (\ref{Newton-Condition}).

We emphasize once again the Newton step is computed on the restricted subspace $\bfu_{\Gamma_j}=0$. Given $\widehat{\bfu}_{\Gamma_j}=0$ (because $\Gamma_j = \widehat{\I}_0$), we have the following chain of inequalities:
%\begin{align}
\begin{eqnarray}
	&& \| \widetilde{\zeta}^{j+1} - \widehat{\zeta} \| 
	= \| (\widetilde{\zeta}^{j+1} - \widehat{\zeta})_{:\OG_j} \| \notag \\
	& \mathop{=}\limits^{(\ref{Newton-Eq-2})} & \| (\zeta^{j+1/2} - \widehat{\zeta})_{:\OG_j} - ( \nabla^2_{:\OG_j}  g(\zeta^{j+1/2}) )^{-1} \nabla_{:\OG_j}  g(\zeta^{j+1/2})\| \notag \\
	&\leq & \frac{1}{\sigma_g} \| \nabla^2_{:\OG_j}  g(\zeta^{j+1/2})(\zeta^{j+1/2} - \widehat{\zeta})_{:\OG_j} - \nabla_{:\OG_j}  g(\zeta^{j+1/2})\| \notag \\
	& \mathop{=}\limits^{(\ref{P^k-stat})} 	&
	\frac{1}{\sigma_g} \| \nabla^2_{:\OG_j}  g(\zeta^{j+1/2})(\zeta^{j+1/2} - \widehat{\zeta})_{:\OG_j} - (\nabla_{:\OG_j}  g(\zeta^{j+1/2}) - \nabla_{:\OG_j}  g(\widehat{\zeta}) )\| \notag \\
	& = & \frac{1}{\sigma_g} \| \int_{0}^{1} ( \nabla^2_{:\OG_j}  g(\zeta^{j+1/2}) - \nabla^2_{:\OG_j}  g( \widehat{\zeta} + s(\zeta^{j+1/2} - \widehat{\zeta})) ) ( \zeta^{j+1/2} - \widehat{\zeta} ) d s \| \notag \\
	&\leq &  \frac{1}{\sigma_g} \int_{0}^{1}  \| \nabla^2_{:\OG_j}  g(\zeta^{j+1/2}) - \nabla^2_{:\OG_j}  g( \widehat{\zeta} + s (\zeta^{j+1/2} - \widehat{\zeta})) \| \| \zeta^{j+1/2} - \widehat{\zeta} \| d s  \notag \\
	&\leq &  \frac{L_g}{\sigma_g} \| \zeta^{j+1/2} - \widehat{\zeta} \|^2 \int_{0}^{1} (1- s) d s = \frac{L_g}{2\sigma_g} \| \zeta^{j+1/2} - \widehat{\zeta} \|^2, \label{quadra_half}
	%\end{align}
\end{eqnarray} 
where the first inequality follows from strong convexity of $g$.
Given $\zeta^{j+1/2}$ converges to $\widehat{\zeta}$, 
the relation in (\ref{quadra_half}) implies that $\widetilde{\zeta}^{j+1}$ can be
made arbitrarily close to $\widehat{\zeta}$. 
In fact,  taking $\lim_{j \to \infty}  \zeta^{j+1/2} = \widehat{\zeta} $ into account, we have
\begin{equation} \label{lim_zetatilde}
	\lim_{j \to \infty}	\widetilde{\zeta}^{j+1} = \lim_{j \to \infty} \zeta^{j+1/2} =
	\lim_{j \to \infty} \zeta^{j+1} = 
	\widehat{\zeta}.
\end{equation}
Bearing in mind that the strict complementarity
condition holds, we have
\begin{equation} \label{+0equal}
	\begin{aligned}
		& \| \widetilde{\bfu}^{j+1}_+ \|_0 \mathop{=}\limits^{(\ref{Newton-Eq})} \| (\widetilde{\bfu}^{j+1}_{\OG_j})_+ \|_0 
		\mathop{=}\limits^{(\ref{set_iden})}   \| (\widehat{\bfu}_{ \OG_j })_+ \|_0,  \\
		& \| \bfu^{j+1/2}_+ \|_0 \mathop{=}\limits^{(\ref{Gradient-Step})} \| (\bfu^{j+1/2}_{ \OG_{j} })_+ \|_0 \mathop{=}\limits^{(\ref{set_iden})}   \| (\widehat{\bfu}_{ \OG_j })_+ \|_0. 
	\end{aligned}
\end{equation}
%%%%%%%%%%%%%%%%%%%%%%%%%%%% Not Needed
%Moreover, we can give the following estimation
%\begin{eqnarray*}
%	\| \widetilde{\zeta}^{j+1} - \zeta^{j+1/2} \| &=& \| (\widetilde{\zeta}^{j+1} - \zeta^{j+1/2})_{ :\OG_j } \|  \\ 
%	& \mathop{=}\limits^{(\ref{Newton-Eq-2})} & \| ( \nabla^2_{ :{\OG_j }}  g(\zeta^{j+1/2}) )^{-1} \nabla_{ :{\OG_j }}  g(\zeta^{j+1/2}) \| \\
%	&\mathop{=}\limits^{\ref{P^k-stat})} &  \| ( \nabla^2_{ :{\OG_j }}  g(\zeta^{j+1/2}) )^{-1} (\nabla_{ :{\OG_j }}  g(\zeta^{j+1/2}) -  \nabla_{ :{\OG_j }}  g(\widehat{\zeta}) ) \| \\
%	&\leq & \frac{1}{\sigma_g} \| \nabla_{ :{\OG_j }}  g(\zeta^{j+1/2}) -  \nabla_{ :{\OG_j }}  g(\widehat{\zeta}) \| 	\leq  \frac{l_g}{\sigma_g} \| \zeta^{j+1/2} - \widehat{\zeta} \|. 
%\end{eqnarray*}
%where the last two inequalities follows from the strong convexity of $g$ and the Lipschitz continuity of its gradient.
%%%%%%%%%%%%%%%%%%%%%%%%%%%%%%%%%%%%%%%%%%%%%%%%%%%%%%%%%%%%%%
The facts established in (\ref{lim_zetatilde}) and (\ref{+0equal}) means that
there exists a sufficiently large $\widehat{j}$ such that $\forall j \geq \widehat{j}$,  the following hold:
\begin{align}
	&  G ( \widetilde{\zeta}^{j+1}) - G ( \zeta^{j+1/2}) \notag \\
	= & g(\widetilde{\zeta}^{j+1}) - g(\zeta^{j+1/2}) + \lambda \| \widetilde{u}^{j+1}_+ \|_0 - \lambda \| u^{j+1/2}_+ \|_0 \notag \\
	\mathop{=}\limits^{(\ref{+0equal})}&  \langle \nabla g( \zeta^{j+1/2} ) , \widetilde{\zeta}^{j+1} - \zeta^{j+1/2}  \rangle + \frac{1}{2}  ( \widetilde{\zeta}^{j+1} - \zeta^{j+1/2} )^\top \nabla^2 g(\zeta^{j+1/2}) ( \widetilde{\zeta}^{j+1} - \zeta^{j+1/2} ) \notag \\
	&  + o( \| \widetilde{\zeta}^{j+1} - \zeta^{j+1/2} \|^2 ) \notag \\
	{=}&  \langle \nabla_{ :\OG_j } g( \zeta^{j+1/2} ) , (\widetilde{\zeta}^{j+1} - \zeta^{j+1/2})_{ :\OG_j }  \rangle + o( \| \widetilde{\zeta}^{j+1} - \zeta^{j+1/2} \|^2 ) \notag \\
	& + \frac{1}{2}  ( \widetilde{\zeta}^{j+1} - \zeta^{j+1/2} )^\top_{ :{\OG_j }} \nabla^2_{ :{\OG_j }} g(\zeta^{j+1/2}) ( \widetilde{\zeta}^{j+1} - \zeta^{j+1/2} )_{ :\OG_j} \notag \\
	= & - \frac{1}{2}  ( \widetilde{\zeta}^{j+1} - \zeta^{j+1/2} )^\top_{ :{\OG_j }} \nabla^2_{ :{\OG_j }} g(\zeta^{j+1/2}) ( \widetilde{\zeta}^{j+1} - \zeta^{j+1/2} )_{ :\OG_j } + o( \| \widetilde{\zeta}^{j+1} - \zeta^{j+1/2} \|^2 ) \notag \\
	\leq & -\frac{\sigma_g}{2} \| \widetilde{\zeta}^{j+1} - \zeta^{j+1/2} \|^2  + o( \| \widetilde{\zeta}^{j+1} - \zeta^{j+1/2} \|^2 ) \notag \\
	\mathop{\leq} \limits^{(\ref{lim_zetatilde})} & -\frac{\sigma_g}{4} \| \widetilde{\zeta}^{j+1} - \zeta^{j+1/2} \|^2,   
\end{align}
which implies that the Newton step is always accepted in (\ref{Newton-Condition}) when $j \geq \widehat{j}$.
\hfill $\Box$ 
%%%%%%%%%%%%%%%%%%%%%%%%%%%%%%%%%%%%%%%%%%%%%%%%%%%%%%%%%%%%

We have so far established that (i) the active set is correctly identified after a finite
number of iterations, and (ii) the Newton step is accepted. 
In other words, the Newton method is carried out in a (fixed) subspace.
Its quadratic convergence would follow from the classical theory due to the strong
convexity of $g(\cdot)$. 
However, we should note that we have applied a gradient step before the Newton step. 
To put it another way, our algorithm eventually carries out a gradient step, followed by a Newton step, and both are done in a same subspace. This scheme is illustrated as follows:
\[
\zeta^{j} \ \longrightarrow \ \zeta^{j+1/2} \ (\mbox{by Gradient Step}) \ 
\longrightarrow \ 
\zeta^{j+1} \ (\mbox{by Newton Step}).
\]
The next result below establishes the quadratic convergence of the sequence $\{\zeta^j\}$.

\begin{theorem}[Quadratic Convergence] \label{Thm-Quadratic}
	Suppose Assumption \ref{Assumption-f} hold. 
	Let $\{ \zeta^{j} \}_{j \in \mathbb{N}}$ be a sequence generated by Alg.~\ref{Alg-GSN} converging to a ${\rm P}$-stationary point $\widehat{\zeta}$ with ${\rm P}$-stationary multiplier $\widehat{\bfy}$. 
	We further assume that the strict complementarity condition holds at
	$(\widehat{\bfu}, \widehat{\bfy})$.
	Then $\{ \zeta^{j} \}_{j \in \mathbb{N}}$ converges to $\widehat{\zeta}$ quadratically, namely, there exists an index $\widehat{j}$ such that
	\begin{equation} \label{pamn_quadra_con}
		\| \zeta^{j+1} - \widehat{\zeta} \| \leq \eta \| \zeta^{j} - \widehat{\zeta} \|^2, \quad \forall j \geq \widehat{j},
	\end{equation}
	where $\eta := \frac{L_g}{\sigma_g} \max\{ \alpha^2 \ell_g^2 + \eta_1^2,\; ( 1+ \alpha \ell_g )^2 + (1+\eta_1)^2 \}$, $\eta_1 := t\ell_g  ( 1+ \alpha \ell_g)$.
\end{theorem}

{\bf Proof.} Since the active set $\widehat{\I}_0$ is correctly identified for
$j \ge \widehat{j}$, in the proof below, we are going to use the following facts without 
explicitly referring to it: 
$
\bfu^{j+1/2}_{\Gamma_j} =\widehat{\bfu}_{\Gamma_j}= 0
$ and $\widehat{\bfy}_{\OG_j}=0$.
From Lemma \ref{quadra_lemma}, we know that (\ref{quadra_half}) holds for all $j \geq \widehat{j}$. 
To verify (\ref{pamn_quadra_con}), we need to estimate the relationship between 
$\| \zeta^{j+1/2} - \widehat{\zeta} \|$ and $\| \zeta^{j} - \widehat{\zeta} \|$. 
For $\| \bfu^{j+1/2} - \widehat{\bfu} \|$, we give the following estimation for $j \geq \widehat{j}$,
\begin{align}
	&\| \bfu^{j+1/2} - \widehat{\bfu} \| \notag \\
	= &
	\| (\bfu^{j+1/2} - \widehat{\bfu})_{\OG_j} \| 
	=\| ( \bfu^{j} - \widehat{\bfu} )_{\OG_j} + \alpha \bfy^{j}_{\OG_j} \| \notag \\
	= & \| ( \bfu^{j} - \widehat{\bfu} )_{\OG_j} + \alpha (\bfy^{j}_{\OG_j} - \widehat{\bfy}_{\OG_j} ) \| \leq \| \bfu^{j} - \widehat{\bfu} \| + \alpha \| \bfy^{j} - \widehat{\bfy} \| \notag \\
	= & \| \bfu^{j} - \widehat{\bfu} \| + \alpha \| \nabla_\bfu g( \zeta^{j} ) - \nabla_\bfu g( \widehat{\zeta} ) \| \leq  \alpha \ell_g\| \bfx^{j} - \widehat{\bfx} \| + ( 1+ \alpha \ell_g ) \| \bfu^{j} - \widehat{\bfu} \|, \label{uj+1/2-hat}
\end{align}
where the last inequality used the Lipschitz continuity of $\nabla g(\zeta)$. We also have the following inequalities
%\begin{align}
%	(\sigma_g + \theta_{\min}(Q)) \| w^{j+1/2} - w^j \| \leq & \| ( \nabla_w g(w^{j+1/2}, u^{j+1/2}) + Q w^{j+1/2} ) - ( \nabla_w g(w^{j}, u^{j+1/2}) + Q w^{j} ) \| \notag \\
%	\mathop{=}\limits^{(\ref{P^k-stat},\ref{w-var})}& \| \nabla_w g(w^{j}, u^{j+1/2}) - \nabla_w g(\widehat{w}, \widehat{u}) \| \leq l_g ( \| w^{j} - \widehat{w} \| + \| u^{j+1/2} - \widehat{u} \| ), \notag
%\end{align}
%where the first inequality follows from the strongly convexity of $g$ in Assumption \ref{g_c2+sc} and $Q$ is positive definite, and the last inequality follows from the Lipschitz continuity in Assumption \ref{lc_g}. This directly leads to
%\begin{align} 
%	\| w^{j+1/2} - w^j \| \leq \frac{l_g}{\sigma_g + \theta_{\min}(Q)} ( \| w^{j} - \widehat{w} \| + \| u^{j+1/2} - \widehat{u} \| ) \label{wj+1/2-wj}
%\end{align}
%Then for $\| w^{j+1/2} - \widehat{w} \|$, we have the following estimation
\begin{align} \label{wj+1/2-hat}
	\| \bfx^{j+1/2} - \widehat{\bfx} \| \mathop{=}\limits^{(\ref{Gradient-Step})}&  \| \bfx^{j} - \widehat{\bfx} - t \nabla_\bfx g(\bfx^j,\; \bfu^{j+1/2}) \| \notag \\
	\mathop{=}\limits^{(\ref{P^k-stat})}& \| \bfx^{j} - \widehat{\bfx} - t (\nabla_\bfx g(\bfx^j,\bfu^{j+1/2}) - \nabla_\bfx g(\widehat{\bfx},\widehat{\bfu}) )\| \notag \\
	\leq\ & \| \bfx^j - \widehat{\bfx} \| + t \ell_g ( \| \bfx^j - \widehat{\bfx} \| + \| \bfu^{j+1/2} - \widehat{\bfu} \| ) \notag \\
	= \ & (t\ell_g + 1) \| \bfx^j - \widehat{\bfx} \| + t \ell_g \| \bfu^{j+1/2} - \widehat{\bfu} \| \notag \\
	\mathop{=}\limits^{(\ref{uj+1/2-hat})} & \eta_1 \| \bfx^j - \widehat{\bfx} \| + (1+\eta_1) \| \bfu^j - \widehat{\bfu} \|,
	%	\mathop{=}\limits^{(\ref{sub_kkt})}& \| w^{j} - \widehat{w} \| + (1/t) \| \nabla_w g(w^j,u^{j+1/2}) - \nabla_w g(\widehat{w},\widehat{u}) \| \notag \\
	%	\leq & \| w^{j} - \widehat{w} \| + (l_g/t) ( \| w^{j} - \widehat{w} \| + \| u^{j+1/2} - \widehat{u} \| ) \notag \\
	%	\mathop{\leq}\limits^{(\ref{uj+1/2-hat})} & \eta_1 \| w^{j} - \widehat{w} \| + 1+\eta_1 \| u^{j} - \widehat{u} \| \label{wj+1/2-hat}.	
	%	(1+l_g/t+\alpha l_g^2/t)\| w^{j} - \widehat{w} \| + (l_g(1+\alpha l_g)/t) \| u^{j} - \widehat{u} \|  \label{wj+1/2-hat} .
\end{align}
and 
\begin{align} \label{j+1/2-h_j-h}
	\| \zeta^{j+1/2} - \widehat{\zeta} \|^2 =& \| \bfu^{j+1/2} - \widehat{\bfu} \|^2 + \| \bfx^{j+1/2} - \widehat{\bfx} \|^2 \notag \\
	\mathop{\leq}\limits^{(\ref{uj+1/2-hat},\ref{wj+1/2-hat})} & 2(\alpha^2 \ell_g^2 + \eta_1^2) \| \bfx^j - \widehat{\bfx} \|^2 + 2(( 1+ \alpha \ell_g )^2 + (1+\eta_1)^2) \| \bfu^j - \widehat{\bfu} \|^2 \notag \\
	\leq~~ &  2 \max\{ \alpha^2 \ell_g^2 + \eta_1^2, ( 1+ \alpha \ell_g )^2 + (1+\eta_1)^2 \} \| \zeta^j - \widehat{\zeta} \|^2.
\end{align}
Combining (\ref{j+1/2-h_j-h}) with (\ref{quadra_half}), we obtain $\| \zeta^{j+1} - \widehat{\zeta} \| \leq \eta \| \zeta^{j} - \widehat{\zeta} \|^2.$
\hfill $\Box$

%%%%%%%%%%%%%%%%%%%%%%%%%%%%%%%%%%%%%%%%%%%%%%%%%%%%%
\section{iNALM: Inexact Newton Augmented Lagrangian Method} \label{Section-NALM}

In this section, we study the Newton Augmented Lagrangian Method (NALM), which
largely follows the standard framework of ALM and has a primal step and a multiplier step as its core iteration.
A key departure from other ALMs is that the subproblem for the primal step
is solved inexactly by Newton's method in Alg.~\ref{Alg-GSN}.
The inexactness of this primal step is more practical than exact solution methods,
yet presents some challenges in ensuring monotone decrease of the augmented Lagrange function.
To resolve this issue, we will design computable criteria to measure the inexactness and 
show that the inexact NALM is well defined.
This is done in Subsect.\ref{Subsection-Inexact}.
We prove the boundedness of the iterates and their global convergence of NALM in Subsect.\ref{Subsection-Global} under reasonable assumptions. 
Further with SOSC, we establish its R-linear convergence rate both in terms of the Lyapunov function and the iterates in Subsect.\ref{Subsection-R}.

%%%%%%%%%%%%%%%%%%%%%%
\subsection{Inexact NALM} \label{Subsection-Inexact}

We remind ourself that the main proposal is to use the gradient-subspace-Newton method Alg.~\ref{Alg-GSN} to
approximately solve subproblem (\ref{Subproblem}), which is recalled here for 
convenient reference. Given the current iterate 
$(\bfx^k, \bfu^k, \bfy^k)$, solve
\begin{eqnarray} \label{Subproblem-k}
	(\bfx^{k+1}, \bfu^{k+1})  &\approx&
	\arg\min_{\bfx, \bfu} \ \V_{\rho, \mu} (\bfx, \bfu, \bfy^k, \bfx^k) 
	= \arg\min_{\bfx, \bfu}\ g_k(\bfx, \bfu) + \lambda \| \bfu_+\|_0,
\end{eqnarray}
where 
\[
g_k(\bfx, \bfu) := f(\bfx) + \langle \bfy^k, A\bfx + \bfb - \bfu \rangle + \frac{\rho}2 
\|A\bfx + \bfb - \bfu \|^2 +  \frac {\mu}2 \| \bfx - \bfx^k\|^2 .
\]
Under Assumption~\ref{Assumption-f}, $\nabla g_k (\bfx, \bfu)$ are uniformly Lipschitz continuously with a module $\overline{\ell}_g >0 $ for any $k = 1,2,\cdots.$ 
As proved in Thm.~\ref{Thm-Global} under Assumption~\ref{Assumption-f},
Newton's method is able to find a P-stationary triplet 
$(\widehat{\bfx}^{k+1}, \widehat{\bfu}^{k+1}, \widehat{\bfy}^{k+1})$ of
subproblem (\ref{Subproblem-k}) and it satisfies the following stationary condition:
\begin{align}\label{KKT-k}
	\left\{ \begin{aligned} 
		& \nabla_\bfx g_k (\widehat{\bfx}^{k+1},\widehat{\bfu}^{k+1} ) = 0, \\ &\widehat{\bfu}^{k+1} \in \Prox_{\alpha \lambda h(\cdot)} (\widehat{\bfu}^{k+1} + \alpha \widehat{\bfy}^{k+1}), \ \\
		&\widehat{\bfy}^{k+1} = - \nabla_\bfu g_k (\widehat{\bfx}^{k+1}, \widehat{\bfu}^{k+1} ).
	\end{aligned} \right.
\end{align}
Let $\widehat{\Gamma}_{k+1}$ be the active set at $\widehat{\bfu}^{k+1}$.
The second condition in (\ref{KKT-k}) implies the complementarity condition for
$(\widehat{\bfu}^{k+1}, \widehat{\bfy}^{k+1})$:
\be \label{R2}
\widehat{u}^{k+1}_i = 0, \ \ \mbox{for} \ i \in \widehat{\Gamma}_{k+1} \quad
\mbox{and} \quad 
\widehat{y}^{k+1}_i = 0, \ \ \mbox{for} \ i \not\in \widehat{\Gamma}_{k+1}.
\ee 
In terms of the objective values, the second condition in (\ref{KKT-k}) also implies
\be \label{R3}
\frac{\alpha^2}2 \| \widehat{\bfy}^{k+1} \|^2 + \alpha \lambda \| \widehat{\bfu}^{k+1}_+ \|_0 
= \Phi_{\alpha \lambda h(\cdot)} ( \widehat{\bfu}^{k+1} + \alpha \widehat{\bfy}^{k+1}),
\ee 
where $\Phi_{\alpha \lambda h(\cdot)}(\cdot)$ is the Moreau envelop of the function
$(\alpha \lambda) h(\cdot)$.
Given any pair $(\bfx, \bfu)$, we may measure its closeness to 
$(\widehat{\bfx}^{k+1}, \widehat{\bfu}^{k+1})$ by computing the residuals for the first
equation in (\ref{KKT-k}) and those in (\ref{R2}) and (\ref{R3}). 
\be \label{Residual-Function}
\left\{
\begin{array}{l}
	\R_1 (\bfx, \bfu) := \| \nabla_\bfx g_k(\bfx, \bfu) \| \\ [1ex]
	\R_2 (\bfx, \bfu) :=  \| [\bfu_{\Gamma};\; \alpha \nabla_\OG g_k(\bfx, \bfu)] \| \\ [1ex]
	\R_3 (\bfx, \bfu) := \frac{\alpha^2}2 \| \nabla_\bfu g_k(\bfx, \bfu) \|^2 
	+ \alpha \lambda \| {\bfu}_+ \|_0
	- \Phi_{\alpha \lambda h(\cdot)} ( {\bfu} - \alpha \nabla_\bfu g_k(\bfx, \bfu) ),
\end{array} 
\right.
\ee 
where $\Gamma$ is the active set at $\bfu$ and we replaced $\bfy$ with $-\nabla_\bfu g_k(\bfx, \bfu)$ (this is justified by the third equation in (\ref{KKT-k})).
Apparently, $\R_i (\widehat{\bfx}^{k+1}, \widehat{\bfu}^{k+1})=0$ for $i=1,2,3$.
Therefore, we apply the Newton method Alg.~\ref{Alg-GSN} to subproblem
(\ref{Subproblem-k}) and terminate it as soon as the residuals fall below 
certain threshold. This leads to the inexact Newton Augmented Lagrangian Method (iNALM)
described in Alg.~\ref{Alg-iNALM}.

\begin{algorithm}[H]
	\caption{iNALM: inexact Newton Augmented Lagrangian Method} \label{Alg-iNALM}
	\begin{algorithmic}
		\STATE{\textbf{Initialization:}} Given a positive sequence $\{\epsilon_k\}_{k \in \mathbb{N}}$ going to $0$, two constants $c_1>0$ and $c_2>0$,
		and initial point $(\bfx^0, \; \bfu^0, \; \bfy^0)$.
		
		\FOR{$k = 0,1, \cdots$}
		\STATE{\textbf{1. Primal step:}} 
		Starting with $(\bfx^k, \; \bfu^k)$, apply GSNewton Alg.~\ref{Alg-GSN} to
		subproblem (\ref{Subproblem-k}) to generate a sequence
		$\{ \bfx^{k, j}, \; \bfu^{k, j}\}$, $j=0, 1, \ldots, .$ 
		We terminate GSNewton as soon as there exists an index $j_k$ such that
		\be \label{Residual-error}
		\left\{
		\begin{array}{l}
			\R_1( \bfx^{k, j_k}, \bfu^{k, j_k}  ) \le c_1 \| \bfx^{k, j_k} - \bfx^k \|, \\ [1ex]
			\R_2( \bfx^{k, j_k}, \bfu^{k, j_k}  ) \le c_2 \| \bfx^{k, j_k} - \bfx^k \|^2, \\ [1ex]
			\R_3( \bfx^{k, j_k}, \bfu^{k, j_k}  ) \le \epsilon_k.
		\end{array} 
		\right.
		\ee 
		Let
		\[
		\bfx^{k+1} = \bfx^{k, j_k} \quad \mbox{and} \quad
		\bfu^{k+1} = \bfu^{k, j_k}.
		\]
		
		\STATE{\textbf{2. Multiplier step:}}
		\begin{equation} \label{Multiplier-update}
			\bfy^{k+1} = \bfy^k + \rho(A \bfx^{k+1} + \bfb - \bfu^{k+1}).
		\end{equation}
		
		\ENDFOR
	\end{algorithmic}
\end{algorithm}

We note that the residual errors in (\ref{Residual-error}) are controlled by the distance between the current iterate $\bfx^{k, j}$ and the initial iterate $\bfx^k$ for each subproblem. 
As we will prove, $\{\bfx^k\}$ will converge to a stationary point, the error eventually converges to $0$. 
In other words, $(\bfx^{k+1}, \bfu^{k+1})$ is an approximate P-stationary point of subproblem (\ref{Subproblem-k}). The next lemma proves that such approximate point satisfying (\ref{Residual-error}) exists when $\alpha \in (0,1/\overline{\ell}_g)$ and $\bfx^{k} \neq \widehat{\bfx}^{k+1}$, where $\widehat{\bfx}^{k+1}$ has been defined in \eqref{KKT-k}.

\begin{lemma} \label{Lemma-welldefined}
	If Assumption \ref{Assumption-f} ,$\alpha \in (0,1/\overline{\ell}_g)$ and $\bfx^{k} \neq \widehat{\bfx}^{k+1}$, then there exists $j_k$ such that
	the termination criterion (\ref{Residual-error}) is satisfied. 
	Consequently, Alg.~\ref{Alg-iNALM} is well-defined and
	\begin{equation} \label{descent}
		\V_{\rho, \mu}(\bfx^{k+1}, \bfu^{k+1}, \bfy^k, \bfx^k  )
		\le \V_{\rho, \mu}(\bfx^{k}, \bfu^{k}, \bfy^k, \bfx^k  ).
	\end{equation}
\end{lemma}

{\bf Proof.}
Let $(\bfx^{k, j}, \bfu^{k, j}, \bfy^{k,j})$, $j \in \mathbb{N}$ be the sequence generated by GSNewton Alg.~\ref{Alg-GSN} applied to subproblem (\ref{Subproblem-k}).
Let $\Gamma_{k,j}$ be the corresponding active set defined in the
identification step in Alg.~\ref{Alg-GSN}.
The algorithm also generates a half-step sequence $\{\bfx^{k, j+1/2}, \bfu^{k, j+1/2}\}$.
It follows from Thm.~\ref{Thm-Global} that 
the sequence converges to a P-stationary triplet, denoted by
$(\widehat{\bfx}^{k+1}, \widehat{\bfu}^{k+1}, \widehat{\bfy}^{k+1})$, which satisfies
(\ref{KKT-k}).
%\be \label{gk-stat}
%\nabla_\bfx g_k(\widehat{\bfx},\widehat{\bfu}) = 0, \quad
%\widehat{\bfu} \in {\rm Prox_{\alpha\lambda h(\cdot)}} ( \widehat{\bfu} + \alpha \widehat{\bfy} ), \quad
%\widehat{\bfy} = -\nabla_\bfu g_k(\widehat{\bfx},\widehat{\bfu}).
%\ee 
Let $\zeta^{k, j} := (\bfx^{k,j}, \bfu^{k, j})$.
%There are two cases to be considered depending whether $\bfx^k$ is distinctive from $\widehat{\bfx}^{k+1}$.

%{\bf Case I}: 
Since $\bfx^{k, j}$ converges to $\widehat{\bfx}^{k+1}$ and we have assumed $\bfx^k \neq \widehat{\bfx}^{k+1}$, there
exists $\overline{\epsilon} > 0$
such that following inequality holds
\begin{align} 
	c_1 \| \bfx^{k,j} - \bfx^k \| > \overline{\epsilon}
	\quad \mbox{and} \quad
	c_2 \| \bfx^{k,j} - \bfx^k \|^2 > \overline{\epsilon}  \label{mid}
\end{align} 
By (\ref{KKT-k}), we have
\begin{align} \label{limw}
	&\lim_{j \to \infty} \nabla_\bfx g_k(\bfx^{k, j}, \bfu^{k, j}) = \nabla_\bfx g_k(\widehat{\bfx}^{k+1},\widehat{\bfu}^{k+1}) = 0, \\
	&\lim_{j \to \infty}\| [ \bfu^{k,j}_{\Gamma_{k,j}}; \alpha \nabla_{\OG_{k,j}} g_k( \bfx^{k, j}, \bfu^{k, j} ) ] \| 
	\mathop{=}\limits^{(\ref{Gradient-Step})}
	\lim_{j \to \infty} \| \bfu^{k,j+1/2} - \bfu^{k,j} \| = 0. \label{lim2}
\end{align}
It follows from the definition of Moreau envelop and (\ref{Proximal-Update}) that
\begin{align}\label{lim3}
	&	\lim_{j \to \infty} \frac{\alpha^2}{2} \| \nabla_\bfu g_{k} ( \zeta^{k,j} ) \|^2 + \alpha\lambda \| \bfu_+^{k,j}   \|_0 - \Phi_{\alpha\lambda h(\cdot)} ( \bfu^{j} - \alpha \nabla_\bfu g_{k} ( \zeta^{k,j} ) ) ) \notag  \\
	&= \lim_{j \to \infty}   - \frac{1}{2} \| \bfu^{k,j+1/2} - \bfu^{k,j} \|^2 - \alpha \langle \nabla_\bfu g_{k} ( \zeta^{k,j} ), \bfu^{k,j+1/2} - \bfu^{k,j} \rangle \notag \\
	&+ \alpha\lambda (\| \bfu_+^{k,j}   \|_0 -  \| \bfu^{k,j+1/2}_+ \|_0)    \mathop{=}\limits^{(\ref{lim-G})}  0.
\end{align}
From (\ref{mid})--(\ref{lim3}), we know that there exists positive integer 
$j_k$
such that the following inequality holds
\begin{align}
	& \| \nabla_\bfx g_k(\zeta^{k, {j_k}}) \| \leq \overline{\epsilon} \mathop{<}\limits^{(\ref{mid})} c_1 \| \bfx^{k,{j_k}} - \bfx^k \|, \notag \\
	& \| ( \bfu^{j_k}_{\Gamma_{k,{j_k}}}; \alpha \nabla_{\OG_{k,{j_k}}} g_k( \zeta^{k, {j_k}} ) ) \| \leq \overline{\epsilon} \mathop{<}\limits^{(\ref{mid})} c_2 \| \bfx^{k,{j_k}} - \bfx^k \|^2, \notag \\
	& (\alpha^2/2) \| \nabla_\bfu g_{k} ( \zeta^{k,{j_k}} ) \|^2 + \alpha\lambda \| \bfu_+^{k,{j_k}}   \|_0 - \Phi_{\alpha\lambda h(\cdot)} ( \bfu^{k,{j_k}} - \alpha \nabla_\bfu g_{k} ( \zeta^{k,{j_k}} ) ) ) \leq \epsilon_k. \notag
\end{align}
Taking $\zeta^{k+1} = \zeta^{k,{j_k}}$, the above inequalities imply that stopping criteria (\ref{Residual-error}) is met.
%Since $\alpha = 1/\rho$, we have 
%\begin{align}\label{u-alp_hk}
%	\bfu - \alpha \nabla_\bfu g_k(\zeta) = A \bfx + \bfb + \bfy^k/\rho.
%\end{align}
%Let us consider the first half-step iterate in (\ref{Gradient-Step}).
%It follows from (\ref{Proximal-Update}) that 
%\begin{align} 
%	\bfu^{k,1/2} \in& {\rm Prox}_{\alpha\lambda h(\cdot)}( \bfu^{k,0} - \alpha \nabla_\bfu g_k ( \zeta^{k,0} ) ) \mathop{=}\limits^{(\ref{u-alp_hk})}  {\rm Prox}_{\alpha\lambda h(\cdot)}( A \bfx^{k,0} + \bfb + \bfy^k / \rho ) \notag \\
%	= &{\rm Prox}_{\alpha\lambda h(\cdot)}( A\widehat{\bfx} + \bfb + \bfy^k / \rho ). \notag
%\end{align}
%This indicates that $(\bfx^k, \bfu^{k, 1/2})$ is a potential candidate for a P-stationary point of the subproblem (\ref{Subproblem-k}). 
%Indeed, with some technical argument (details omitted), we can prove $\bfx^{k, 1/2} = \bfx^k$ and the Newton step in (\ref{Newton-Condition}) is not met. Hence, the Newton step in this case is not accepted. Therefore, $\zeta^{k, 1}$ is a P-stationary point and the stopping
%criteria (\ref{Residual-error}) is met (the residual is $0$ in this case). In other words, $j_k = 1$.  
%
%Case I and Case II combined shows that the inexact NALM Alg.~\ref{Alg-iNALM} is well defined. Thm.~\ref{Thm-Global} has proved that the objective sequence 
%$\V_{\rho, \mu}(\bfx^{k, j}, \bfu^{k, j}, \bfy^k, \bfx^k)$ is non-increasing for $j \in \mathbb{N}$. This finishes the proof.
\hfill $\Box$

%\begin{remark} 
%%\textcolor{blue}
%{We like to note that the stopping criteria in Alg.~\ref{Alg-iNALM} are also met in finite steps even when $\alpha \neq 1/\rho$. If $\widehat{\bfx}^{k+1} \neq \bfx^k$, then following the same procedure of {\bf Case I} in Thm \ref{Lemma-welldefined}, we can obtain the desired conclusion. If $\widehat{\bfx}^{k+1} = \bfx^k$, then from $\widehat{\bfu}^{k+1} \in \Prox_{\alpha \lambda h(\cdot)} (\widehat{\bfu}^{k+1} + \alpha \widehat{\bfy}^{k+1}), \
%	\widehat{\bfy}^{k+1} = - \nabla_\bfu g_k (\widehat{\bfx}^{k+1}, \widehat{\bfu}^{k+1} ) = \rho( A\widehat{\bfx}^{k+1} + \bfb - \widehat{\bfu}^{k+1} + \bfz^k/\rho )$ in (\ref{KKT-k}), $\widehat{\bfu}^{k+1}$ can be uniquely identified by
%	\begin{align} \notag
%			\widehat{u}^{k+1}_i = \left\{
%			\begin{array}{ll}
%				0, & \ \mbox{if} \left[ A\bfx^k + \bfb + \bfz^k/\rho \right]_i \in [0,\sqrt{2\lambda/(\alpha\rho^2)}], \\ 
%				\left[ A\bfx^k + \bfb + \bfz^k/\rho \right]_i, & \mbox{otherwise}. 
%			\end{array} \right.
%	\end{align}
%We can verify that $(\bfx^{k+1};\bfu^{k+1})$ satisfies 
%the stopping criteria (\ref{Residual-error}).}
%\end{remark}

%%%%%%%%%%%%%%%%%%%%%
%\subsection{Well-posedness} \label{Subsection-Well}

\begin{remark}
	If $\bfx^k = \widehat{\bfx}^{k+1}$, we can actually find a  {\rm P}-stationary point   
	$(\widehat{\bfx}^{k+1}, \widehat{\bfu}^{k+1})$ of
	subproblem (\ref{Subproblem-k}), which satisfy (\ref{descent}) and termination criterion (\ref{Residual-error}). Indeed, the last two formulas in (\ref{KKT-k}) implies that a potential candidate for $\widehat{\bfu}^{k+1}$ must be taken from a finite set
	\begin{align} \notag
		\Theta^k:= \left\{ \bfu \in \mathbb{R}^n : u_i = \left\{\begin{array}{ll}
			0, &~{\rm if}~ [A\bfx^k + \bfb +  \frac{\bfz^k}{\rho}]_i \in [0, \sqrt{2\lambda\alpha}) \\ [1ex]
			0~{\rm or}~[A\bfx^k + \bfb + \frac{\bfz^k}{\rho}]_i,&~{\rm if}~ [A\bfx^k + \bfb +  \frac{\bfz^k}{\rho}]_i \in [ \sqrt{2\lambda\alpha}, \sqrt{\frac{2\lambda}{\alpha\rho^2}}] \\ [1ex]
			[A\bfx^k + \bfb +  \frac{\bfz^k}{\rho}]_i,&~{\rm otherwise}. 
		\end{array} \right.  \right\}.
	\end{align}
	Actually, for any $\bfu \in \Theta^k$, $\R_2 (\widehat{\bfx}^{k+1}, \bfu) = 0$ and $\R_3 (\widehat{\bfx}^{k+1}, \bfu) = 0$ always hold. Thus, we only need to find a $\widehat{\bfu}^{k+1} \in \Theta^k$ such that $\R_1 (\widehat{\bfx}^{k+1}, \widehat{\bfu}^{k+1}) = 0$ and $
	\V_{\rho, \mu}(\widehat{\bfx}^{k+1}, \widehat{\bfu}^{k+1}, \bfy^k, \bfx^k  ) \le \V_{\rho, \mu}(\bfx^{k}, \bfu^{k}, \bfy^k, \bfx^k  ).$	This can be achieved in finite tests since $\Theta^k$ is finite.
\end{remark}

%%%%%%%%%%%%%%%%%%%%%%
\subsection{Global convergence of iNALM} \label{Subsection-Global}

To ensure global convergence of ALM algorithms for nonconvex, nonsmooth composite optimization problems, it is widely accepted (see, e.g., \cite{bolte2018nonconvex,bot2019proximal}) that certain regularity should be assumed. 
In the case of composition with linear operators such as those treated in this paper, 
a full rank assumption is often assumed. 

\begin{assumption} \label{Assumption-Fullrank}
	Suppose the matrix $A$ has full row-rank. Let $\gamma^2 := \lambda_{\min}(AA^\top) >0$
	be the smallest eigenvalue of $AA^\top$.
\end{assumption}

Furthermore, we need to find a proper merit function that measures improvement each iteration. As the Lagrange multiplier $\bfy^k$ is updated in each iteration, the natural 
objective function $\V_{\rho, \mu}$ in (\ref{Subproblem-k}) is not necessarily decreasing although it is decreasing for the sequence $\{(\bfx^{k, j}, \bfu^{k, j})\}$ for each subproblem (see Lemma~\ref{Lemma-welldefined}). 
Instead, we are going to consider an objective sequence defined by
\[
\V_k := \V_{\rho, \beta} (\bfx^k, \bfu^k, \bfy^k, \bfx^{k-1}), \quad k=1,2, \ldots
\]
where $\beta>0$ is properly chosen below.

{\bf Parameters setup:} 
Let $\mu> \sigma_f$ ($\sigma_f$ comes from $\sigma_f$-weak convexity of $f$),  
$\ell_f$ be the Lipschitz modulus of the gradient function of $f$,
and $c_1>0$ used in Alg.~\ref{Alg-iNALM}. 
Define
\begin{align} \label{para_rhobeta}
	& c_3:= \frac{\mu+ \ell_f+c_1}{\gamma}, \ c_4: = \frac{ \mu + c_1}{\gamma}, \ 
	\rho > \max\left\{ \frac{8(c_3^2+c_4^2)}{\mu },  \frac{4\ell_f}{\gamma^2} \right\},~ \beta := \frac{4c_4^2}{\rho}, 
\end{align}
The following fact will be frequently used:
\begin{eqnarray}
	\nabla_\bfx g_k(\bfx^{k+1}, \bfu^{k+1})
	&=& \nabla f(\bfx^{k+1}) + \mu ( \bfx^{k+1} - \bfx^k ) + \rho A^\top( A \bfx^{k+1} + \bfb - \bfu^{k+1} + \bfy^k/\rho ) \nonumber \\ [1ex]
	&=& \nabla f(\bfx^{k+1}) + \mu ( \bfx^{k+1} - \bfx^k ) + A^\top \bfy^{k+1}, 
	\label{gk-gradient} \\
	\nabla_\bfu g_k(\bfx^{k+1}, \bfu^{k+1}) &=& - \bfy^{k+1}, \label{gk-ugradient}
\end{eqnarray} 
where we used the update formula in (\ref{Multiplier-update}).
It is also easy to verify the following 
\begin{equation} \label{auglag_des}
	{\mathcal{L}}_\rho ( \bfx^k, \bfu^k, \bfy^k ) - {\mathcal{L}}_\rho ( \bfx^{k+1}, \bfu^{k+1},\bfy^k ) \geq \frac{\mu}{2} \| \bfx^{k+1} - \bfx^{k} \|^2.
\end{equation}
Having defined $\beta$ in (\ref{para_rhobeta}), we now show that there is a sufficient
decrease in $\V_k$.

\begin{lemma}[Descent Property] \label{Lemma-Vk}
	Suppose that both Assumptions \ref{Assumption-f} and \ref{Assumption-Fullrank}  hold, and the parameters are chosen as in (\ref{para_rhobeta}). 
	Let $\{ (\bfx^k,\bfu^k, \bfy^k) \}_{k\in \mathbb{N}}$ be a sequence generated by iNALM. Then we have
	\begin{equation}  \label{Decreasing-Vk}
		\mathcal{V}_k  - \mathcal{V}_{k+1}  \geq \frac{\mu}4 \| \bfx^{k+1} - \bfx^k \|^2,
		\quad k=0,1,\ldots, .
	\end{equation}
\end{lemma}

{\bf Proof.}
Using the identity in (\ref{gk-gradient}), we have
\begin{align} 	
	&\nabla_\bfx g_k (\bfx^{k+1}, \bfu^{k+1}) 
	- \nabla_\bfx g_{k-1} (\bfx^{k}, \bfu^{k}) \notag \\
	=& \nabla f(\bfx^{k+1}) - \nabla f(\bfx^k) + \mu (\bfx^{k+1} - \bfx^k) - \mu (\bfx^k - \bfx^{k-1}) + A^\top (\bfy^{k+1} - \bfy^k). \notag
\end{align}
This equality, with the stopping criterion on $\R_1$ in (\ref{Residual-error}) and Assumption \ref{Assumption-Fullrank}, implies
\begin{align} 
	& \gamma \|\bfy^{k+1} - \bfy^k \| 
	\leq  \| A^\top ( \bfy^{k+1} - \bfy^k )  \| \notag \\
	\leq &\ (\mu+ \ell_f) \| \bfx^{k+1} - \bfx^k \| + \mu \| \bfx^k - \bfx^{k-1} \|  + \|\nabla_\bfx g_{k-1}(\bfx^k,\bfu^k)\| + 	\|\nabla_\bfx g_{k}(\bfx^{k+1},\bfu^{k+1})\| \notag\\
	\leq &\ (\mu+\ell_f+c_1) \| \bfx^{k+1} - \bfx^k \| + (\mu + c_1) \| \bfx^k - \bfx^{k-1} \|.
\end{align}
Equivalently, we have
\begin{equation} \label{eq3.5}
	\| \bfy^{k+1} - \bfy^k \| \leq c_3 \| \bfx^{k+1} - \bfx^k \| + c_4 \| \bfx^k - \bfx^{k-1} \|.
\end{equation}
Using the fact $(t_1 + t_2)^2 \leq 2 t_1^2 + 2t_2^2 $, the bound in (\ref{eq3.5}) yields
\begin{equation} \label{eq3.6}
	\| \bfy^{k+1} - \bfy^k \|^2 
	\leq 2c_3^2 \| \bfx^{k+1} - \bfx^k \|^2 + 2c_4^2 \| \bfx^k - \bfx^{k-1} \|^2 .
\end{equation}
We get the following estimation:
\begin{align}
	&{\mathcal{L}}_\rho( \bfx^k,\bfu^k,\bfy^k ) 
	- {\mathcal{L}}_\rho( \bfx^{k+1},\bfu^{k+1},\bfy^{k+1} ) \nonumber \\
	=& {\mathcal{L}}_\rho( \bfx^k, \bfu^k, \bfy^k )  
	- {\mathcal{L}}_\rho( \bfx^{k+1}, \bfu^{k+1}, \bfy^k ) 
	+ {\mathcal{L}}_\rho( \bfx^{k+1}, \bfu^{k+1}, \bfy^k ) 
	- {\mathcal{L}}_\rho( \bfx^{k+1}, \bfu^{k+1}, \bfy^{k+1} ) \nonumber \\
	\mathop{\geq}\limits^{(\ref{auglag_des})} & 
	\frac{\mu}{2} \| \bfx^{k+1} - \bfx^k \|^2 + \langle \bfy^k - \bfy^{k+1}, A \bfx^{k+1} + \bfb - \bfu^{k+1} \rangle 	 \notag \\
	\mathop{\geq}\limits^{(\ref{Multiplier-update})} & \frac{\mu}{2} 
	\| \bfx^{k+1} - \bfx^k \|^2 - \frac{1}{\rho} 	\| \bfy^{k+1} - \bfy^k \|^2  \notag \\	\mathop{\geq}\limits^{(\ref{eq3.6})}  & 
	\left( \frac{\mu}{2} - \frac{2c_3^2}{\rho} \right) \| \bfx^{k+1} - \bfx^k \|^2 - \frac{2c_4^2}{\rho}  \| \bfx^k - \bfx^{k-1} \|^2. \label{auglag_diff}
\end{align}
In terms of $\V_k$, (\ref{auglag_diff}) implies
\begin{equation} \notag
	\begin{aligned} 
		\mathcal{V}_k  - \mathcal{V}_{k+1}	=& 
		{\mathcal{L}}_\rho( \bfx^k, \bfu^k, \bfy^k ) 
		- {\mathcal{L}}_\rho( \bfx^{k+1}, \bfu^{k+1}, \bfy^{k+1} ) 
		+ \frac{\beta}{2}\big( \| \bfx^{k} - \bfx^{k-1} \|^2 
		- \| \bfx^{k+1} - \bfx^k \|^2 \big) \notag\\
		\mathop{\geq}\limits^{(\ref{auglag_diff})}& 
		\underbrace{\left( \frac{\mu}{2} - \frac{2c_3^2}{\rho} - \frac{\beta}{2} \right)}_{\ge \mu/4} 
		\| \bfx^{k+1} - \bfx^k \|^2 + 
		\underbrace{\left( \frac{\beta}{2} - \frac{2c_4^2}{\rho} \right)}_{= 0} 
		\| \bfx^k - \bfx^{k-1} \|^2   \notag \\
		\ge & (\mu/4)  \| \bfx^{k+1} - \bfx^k \|^2.
	\end{aligned}
\end{equation}
This establishes the result. 
\hfill $\Box$

We now establish the boundedness of the generated sequence under the assumption that $f$ is coercive. We note that a similar result has been recently proved in \cite{bot2019proximal}, where the subproblem in their ALM is solved exactly.
In contrast, ours is inexact.

\begin{assumption} \label{Assumption-Coercive}
	The function $f$ is coercive, i.e., $f(\bfx) \rightarrow \infty$ as $\| \bfx\| \rightarrow \infty$.	
\end{assumption}

\begin{lemma}[Boundedness] \label{Lemma-Boundedness}
	Suppose that Assumptions \ref{Assumption-f}, \ref{Assumption-Fullrank} and \ref{Assumption-Coercive} hold. Let the parameters  be chosen to satisfy (\ref{para_rhobeta}). 
	Then the sequence $\{ (\bfx^k, \bfu^k, \bfy^k) \}_{k\in \mathbb{N}}$ 
	generated by iNALM is bounded. Moreover, it holds
	\begin{equation} \label{succsive}
		\lim_{k \to \infty} \| \bfx^{k+1} - \bfx^k \| = 0, \
		\lim_{k \to \infty} \| \bfu^{k+1} - \bfu^k \| = 0, \
		\lim_{k \to \infty} \| \bfy^{k+1} - \bfy^k \| = 0.
	\end{equation}
\end{lemma}

{\bf Proof.}
It follows from (\ref{gk-gradient}), Assumption~\ref{Assumption-Fullrank},
and the stopping criterion on $\R_1$ in (\ref{Residual-error}) that
\begin{equation} \notag
	\gamma \left\|\bfy^{k+1}\right\| \leq  \|A^\top \bfy^{k+1}\| \leq \| \nabla f(\bfx^{k+1}) \| + (\mu + c_1 )\| \bfx^{k+1} - \bfx^{k} \|,
\end{equation}
and thus
\begin{equation} \label{z_bound1}
	\left\|\bfy^{k+1}\right\| \leq \frac{1}{\gamma} \| \nabla f(\bfx^{k+1}) \| + \frac{\mu + c_1}{\gamma} \| \bfx^{k+1} - \bfx^{k} \|.
\end{equation}
Then using the fact $(t_1 + t_2)^2 \leq 2t_1^2 + 2t_2^2 $, we can obtain
\begin{equation} \label{z_bound}
	\|\bfy^{k+1}\|^2 \leq \frac{2}{\gamma^2} \| \nabla f(\bfx^{k+1}) \|^2 + 2c_4^2 \| \bfx^{k+1} - \bfx^{k} \|^2 .
\end{equation} 
Since $f$ is $\ell_f$-smooth, it follows from the Descent Lemma~\ref{Descent-Lemma} that
\begin{align}
	f( \bfx^{k+1} - \frac{1}{\ell_f} \nabla f( \bfx^{k+1} ) ) 
	\le  f(\bfx^{k+1}) - \frac{1}{2 \ell_f} \| \nabla f(\bfx^{k+1}) \|^2. \label{f-nabla_f}
\end{align}

Considering the descent property in Lemma \ref{Lemma-Vk}, we have the following inequalities: 
\begin{align}
	& \mathcal{V}_1\geq \mathcal{V}_{k+1}  + \mathcal{V}_k  - \mathcal{V}_{k+1} \geq \mathcal{V}_{k+1} + (\mu/4) \| \bfx^{k+1} - \bfx^k \|^2  \notag \\
	&= f(\bfx^{k+1}) + \lambda \|\bfu^{k+1}_+\|_0 
	+ \frac{\rho}{2} \| A\bfx^{k+1} + \bfb - \bfu^{k+1} 
	+ \frac{1}{\rho} \bfy^{k+1} \|^2 \notag \\
	& \quad
	+ (\mu/4 + {\beta}/{2})\| \bfx^{k+1} - \bfx^k \|^2 - \frac{1}{2\rho}\|\bfy^{k+1}\|^2 \notag \\
	&\mathop{\geq}\limits^{(\ref{z_bound})} \lambda \|\bfu^{k+1}_+\|_0 
	+ \frac{\rho}{2} \| A\bfx^{k+1} + \bfb - \bfu^{k+1} 
	+ \frac{1}{\rho}\bfy^{k+1} \|^2 + \frac{1}{2} f(\bfx^{k+1}) \notag \\
	&\ + \frac{1}{2} f(\bfx^{k+1}) 
	- \frac{1}{\rho\gamma^2} \| \nabla f(\bfx^{k+1}) \|^2 
	+ ( \frac{\mu}4 + \frac{\beta}{2} - \frac{c_4^2}{\rho} ) \| \bfx^{k+1} - \bfx^k \|^2  \notag \\
	& = \lambda \|\bfu^{k+1}_+\|_0 
	+ \frac{\rho}{2} \| A \bfx^{k+1} + \bfb - \bfu^{k+1} 
	+ \frac{1}{\rho}\bfy^{k+1} \|^2  
	+ ( \frac{1}{4\ell_f} 
	- \frac{1}{\rho \gamma^2} ) \| \nabla f(\bfx^{k+1}) \|^2  \notag \\
	&\ + \frac{1}{2} f(\bfx^{k+1})  + \underbrace{\frac{1}{2} f(\bfx^{k+1}) 
		- \frac{1}{4\ell_f} \| \nabla f(\bfx^{k+1}) \|^2}_{\mbox{for use of (\ref{f-nabla_f})}} 
	+ ( \frac{\mu}4 + \frac{\beta}{2} 
	- \frac{c_4^2}{\rho} ) \| \bfx^{k+1} - \bfx^k \|^2  \notag \\
	&\mathop{\geq}\limits^{(\ref{f-nabla_f})} 
	\lambda \|\bfu^{k+1}_+\|_0 
	+ \frac{\rho}{2} \| A \bfx^{k+1} + \bfb - \bfu^{k+1} 
	+ \frac{1}{\rho}\bfy^{k+1} \|^2 
	+ \frac{1}{2} f( \bfx^{k+1} 
	- \frac{1}{\ell_f} \nabla f( \bfx^{k+1} ) ) \notag \\
	& \ + \frac{1}{2} f(\bfx^{k+1}) + ( \frac{1}{4\ell_f} 
	- \frac{1}{\rho \gamma^2} ) \| \nabla f(\bfx^{k+1}) \|^2 
	+ ( \frac{\mu}4 + \frac{\beta}{2} - \frac{c_4^2}{\rho} ) \| \bfx^{k+1} - \bfx^k \|^2  \notag \\
	&\geq \lambda \|\bfu^{k+1}_+\|_0 + \frac{\rho}{2} \| A\bfx^{k+1} + \bfb - \bfu^{k+1} + \frac{1}{\rho}\bfy^{k+1} \|^2 + \frac{1}{2} f(\bfx^{k+1}) + \frac{1}{2}\inf\limits_{\bfx \in \mathbb{R}^n} f(\bfx)  \notag \\
	&\  + ( \frac{\mu}4 + \frac{\beta}{2} - \frac{c_4^2}{\rho} ) 
	\| \bfx^{k+1} - \bfx^k \|^2 + ( \frac{1}{4\ell_f} - \frac{1}{\rho \gamma^2} ) \| \nabla f(\bfx^{k+1}) \|^2. \notag
\end{align}
Consequently,
\begin{align}
	&	\mathcal{V}_1 - \frac{1}{2}\inf\limits_{\bfx \in \mathbb{R}^n} f(\bfx) 
	\geq \lambda \|\bfu^{k+1}_+\|_0 + \frac{\rho}{2} \| A\bfx^{k+1} + \bfb - \bfu^{k+1} + \frac{1}{\rho} \bfy^{k+1} \|^2 + \frac{1}{2} f(\bfx^{k+1}) \notag \\
	& + \underbrace{\left( \frac{\mu}4 + \frac{\beta}{2} - \frac{c_4^2}{\rho} \right)}_{> 0} 
	\| \bfx^{k+1} - \bfx^k \|^2 
	+ \underbrace{\left( \frac{1}{4\ell_f} - \frac{1}{\rho \gamma^2} \right)}_{> 0} \| \nabla f(\bfx^{k+1}) \|^2. \label{upbound1}
\end{align}
By the coerciveness of $f$, the bound in (\ref{upbound1}) means that $\{\bfx^k\}$ is bounded. Consequently, $\{ \|\nabla f(\bfx^k)\| \}$ is bounded.
The fact in (\ref{z_bound1}) implies the boundedness of $\{ \bfy^k \}$.
From (\ref{Multiplier-update}), we have
\begin{equation} \notag
	\|\bfu^{k+1}\| \leq \| A \| \|\bfx^{k+1}\| + \|\bfb\| + \frac{1}{\rho}( \|\bfy^{k+1} \| + \|\bfy^k\| ),
\end{equation}
which means that $\{\bfu^k\}$ is bounded.
Overall, the sequence $\{ (\bfx^k,\bfu^k, \bfy^k) \}$ generated by iNALM is bounded. 

We now prove the consecutive terms of the sequence become arbitrarily close.
We note that
\begin{align}
	\mathcal{V}_k 
	=& f(\bfx^k) + \lambda \| \bfu^k_+ \|_0 
	+ \frac{\rho}{2} \| A \bfx^k + \bfb - \bfu^k +\frac{1}{\rho} \bfy^k \|^2   + \frac{\beta}{2} \| \bfx^k - \bfx^{k-1} \|^2 - \frac{1}{2\rho}\|\bfy^k\|^2 , \notag
\end{align} 
and that  $f$ is bounded from below, and $\{\bfy^k\}$ is bounded. 
We conclude that $\{ \mathcal{V}_k \}$ is also bounded from below.
Since $\{ \V_k\}$ is a non-increasing sequence, it must converge.
Taking limit on both sides of (\ref{Decreasing-Vk}) gives $
\lim_{k \to \infty} \| \bfx^{k+1} - \bfx^k \| = 0.$
Moreover, taking limit in (\ref{eq3.5}), we have 
$\lim_{k \to \infty} \| \bfy^{k+1} - \bfy^k \| = 0$.
From (\ref{Multiplier-update}), we have 
\[ 
\bfu^{k+1} - \bfu^k = \frac{1}{\rho}( \bfy^k - \bfy^{k-1} ) 
- \frac{1}{\rho}( \bfy^{k+1} - \bfy^{k} ) + A( \bfx^{k+1} - \bfx^k ),
\]
which means that
\[
\|\bfu^{k+1} - \bfu^k \| 
\leq \frac{1}{\rho}\| \bfy^k - \bfy^{k-1}\| 
+ \frac{1}{\rho}\| \bfy^{k+1} - \bfy^{k} \| + \|A\| \| \bfx^{k+1} - \bfx^k \|.
\]
Taking limits on both sides yields $\lim_{k \to \infty} \| \bfu^{k+1} - \bfu^k \| = 0$.\hfill $\Box$

We are ready to prove our main result that any accumulated point of the generated sequence must be a P-stationary point of (\ref{COP-Constrained}).

\begin{theorem}[Convergence to Stationarity] \label{Thm-Stationarity}
	Suppose that Assumptions \ref{Assumption-f}, \ref{Assumption-Fullrank}, and \ref{Assumption-Coercive} hold.
	Let the parameters be chosen to satisfy (\ref{para_rhobeta}). 
	Let $\{ (\bfx^k, \bfu^k, \bfy^k) \}$ be the sequence generated by
	iNALM.
	Any accumulated point $(\bfx^*, \bfu^*, \bfy^*)$ is a {\rm P}-stationary triplet
	of (\ref{COP-Constrained}).
\end{theorem}

{\bf Proof.} 
We have proved that the sequence $\{ (\bfx^k, \bfu^k, \bfy^k) \}$ is bounded.
Therefore, there exists a subsequence $\{ (\bfx^{k}, \bfu^{k}, \bfy^{k}) \}_{k\in \mathcal{K}}$ converging to $(\bfx^*, \bfu^*,\bfy^*)$.
For each iteration $k$, we
define a new point $\overline{\bfu}^{k} \in \mathbb{R}^m$ by
\[
\overline{\bfu}^{k}_{\Gamma_{k}} = 0 \quad \mbox{and} \quad
\overline{\bfu}^{k}_{\OG_k} = \bfu^{k}_{\OG_k} + \alpha \nabla_{\OG_k} g_k (\bfx^k, \bfu^k).
\]
By the formula for the proximal operator of $h(\cdot)$ in (\ref{prox0+})
and the definition $\Gamma_k$ in (\ref{Gamma_j}),  we know 
\begin{align} \notag
	\overline{\bfu}^{k} \in {\rm Prox}_{\alpha\lambda h(\cdot)} ( \bfu^{k} + \alpha \bfy^k )
\end{align} 
Stopping criterion (\ref{Residual-error}) implies
\begin{align}
	\| \bfu^k - \overline{\bfu}^k \| =  \| [\bfu^{k}_{\Gamma_{k}}; \alpha \nabla_{\OG_k} g_{k-1} (\bfx^k, \bfu^k ) ] \| \leq c_2 \| \bfx^{k} - \bfx^{k-1} \|^2. \notag 
\end{align}
Since $\lim_{k \to \infty} \| \bfx^k - \bfx^{k-1} \| = 0$ 
from Lemma \ref{Lemma-Boundedness}, we have
$
\lim_{k\in \mathcal{K},k \to \infty}  \overline{\bfu}^{k}  
= \lim_{k\in \mathcal{K},k \to \infty}  {\bfu}^{k}
=\bfu^*. 
$
The Theorem of Proximal Behavior \cite[Theorem 1.25]{RockWets98} implies
\begin{equation} \label{sub_Psta}
	\bfu^* \in {\rm Prox}_{\alpha\lambda h(\cdot)} ( \bfu^{*} + \alpha \bfy^{*} ).
\end{equation}
We also have the following estimation
\begin{align}
	& \| \nabla f(\bfx^k) + \mu ( \bfx^k - \bfx^{k-1} ) + A^\top \bfy^{k} \| \mathop{=}\limits^{(\ref{gk-gradient})} \| \nabla_\bfx g_{k-1}(\bfx^k,\bfu^k) \| \mathop{\leq}\limits^{(\ref{Residual-error})} c_1 \| \bfx^{k} - \bfx^{k-1} \|,  \notag \\
	& A\bfx^k + \bfb - \bfu^k = (\bfy^{k+1} - \bfy^k)/\rho.  \notag
\end{align}	
Taking $k \in \mathcal{K}$, $k \to \infty$ on both sides of the above two equations and using Lemma~\ref{Lemma-Boundedness}, we have	
\begin{equation} \notag \label{sub_Psta1}
	\nabla f(\bfx^*) + A^\top \bfy^* = 0 \quad \mbox{and} \quad
	A\bfx^* + \bfb = \bfu^* .
\end{equation}
Combining with (\ref{sub_Psta}), we can conclude that $ (\bfx^*,\bfu^*,\bfy^*) $ is a P-stationary triplet of (\ref{COP-Constrained}).
\hfill $\Box$

%%%%%%%%%%%%%%%%%%%%%%%
\subsection{R-linear convergence under SOSC} \label{Subsection-R}

In this part, we study what convergence properties that iNALM may enjoy under the additional condition SOSC. We report three results.
One is that the whole sequence $\{ \bfx^k, \bfu^k, \bfy^k\}$ converges.
The other two are about R-linear convergence in terms of the Lyapunov function values $\V_k$ and in terms of the iterate sequence. 
Throughout this section, we assume that the three Assumptions \ref{Assumption-f}, \ref{Assumption-Fullrank} and \ref{Assumption-Coercive} hold, and that
the parameters are chosen to satisfy (\ref{para_rhobeta}).

%%%%%%%%%%%%%%%%%%%%%%%%%%%%%%%%%%%%%%%%%%%%%%%%
\subsubsection{Sequence convergence}

The first result is a simple consequence of Thm.~\ref{Thm-SOSC}.

\begin{theorem}[Sequence Convergence]\label{whole_con}
	Let $\{ (\bfx^k, \bfu^k, \bfy^k) \}$ be the sequence generated iNALM for problem (\ref{COP-Constrained}).
	Let $\{ (\bfx^*, \bfu^*, \bfy^*) \}$ be one of its accumulated points.
	If SOSC (\ref{SOSC}) is satisfied at $\{ (\bfx^*, \bfu^*, \bfy^*) \}$,
	the following statements are true.
	
	(i) The whole sequence  $\{ (\bfx^k, \bfu^k, \bfy^k) \}$ converges to
	$\{ (\bfx^*, \bfu^*, \bfy^*) \}$.
	
	(ii) It holds
	\begin{equation} \label{k+0iden}
		\begin{aligned}	
			&\lim_{k \to \infty} \| \bfu^k_+ \|_0 = \| \bfu^*_+ \|_0, \\  
			&\lim_{k \to \infty} \mathcal{V}_{k} = \V_* := \mathcal{V}_{\rho,\beta}( \bfx^*, \bfu^*,\bfy^*, \bfx^* ) = f(\bfx^*) + \lambda \| \bfu^*_+\|_0,
		\end{aligned}
	\end{equation}
	and
	\begin{equation} \label{u^*ga0}
		\bfu^*_{\Gamma_{k+1}} = 0 \ \ \mbox{for all $k$ sufficiently large}.
	\end{equation} 
\end{theorem}

{\bf Proof.}
(i) According to Thm.~\ref{Thm-Stationarity}, 
each accumulation point of $\{ (\bfx^k, \bfu^k, \bfy^k) \}$ is a {\rm P}-stationary triplet of (\ref{COP-Constrained}) with constant $\alpha$. 
Since $(\bfx^*,\bfu^*,\bfy^*)$ satisfies (\ref{SOSC}), it follows from Thm.~\ref{Thm-SOSC} that $(\bfx^*, \bfu^*)$ is an isolated {\rm P}-stationary point with constant $\alpha$. 
Thus, $(\bfx^*, \bfu^*)$ is also an isolated accumulation point of sequence $\{ (\bfx^k, \bfu^k) \}$. Taking this and (\ref{succsive}) into consideration, \cite[Prop.~7]{kanzow1999qp} 
indicates that $\{ (\bfx^k, \bfu^k) \}$ must converge to $(\bfx^*, \bfu^*)$. 
Furthermore, it follows from (\ref{gk-gradient}) and the stopping criteria (\ref{Residual-error}) that 
\begin{align} \notag
	\| \nabla f(\bfx^{k+1}) + A^\top \bfy^{k+1} \| 
	\le \| \nabla_\bfx g(\bfx^{k+1}, \bfu^{k+1}) \| + \mu \| \bfx^{k+1} - \bfx^k\|
	\leq (\mu + c_1) \| \bfx^{k+1} - \bfx^k \|.
\end{align}
Taking limits on both sides yields
\begin{align}\notag
	\nabla f(\bfx^*) +  \lim_{k \to \infty} A^\top \bfy^{k+1} = 0. 
\end{align}
Since matrix $A$ is full row rank and $\{\bfy^k\}$ is bounded, we arrive at
\[
\lim_{k \rightarrow \infty} \bfy^k = (AA^\top)^{-1}A \nabla f(\bfx^*)
\] 
and this limit must be $\bfy^*$.
This proves (i).

(ii)
By the definition of the Moreau envelope and the stopping criteria (\ref{Residual-error}) on
the residual function $\R_3$, we have the following estimation
\begin{eqnarray*}
	&& \frac{\alpha^2}{2} \| \bfy^{k+1} \|^2 + \alpha\lambda \| \bfu_+^{k+1}   \|_0 \\
	&\leq& \Phi_{\alpha\lambda h(\cdot)} ( \bfu^{k+1} + \alpha \bfy^{k+1} ) ) + \epsilon_k \leq \frac{1}{2} \| \bfu^* - ( \bfu^{k+1} + \alpha \bfy^{k+1} ) \|^2 
	+ \alpha\lambda \| \bfu^*_+ \|_0 + \epsilon_k  \\
	&= & \frac{\alpha^2}{2} \| \bfy^{k+1} \|^2 
	+ \alpha \langle \bfy^{k+1}, \bfu^{k+1} - \bfu^* \rangle 
	+ \frac{1}{2} \| \bfu^{k+1} - \bfu^* \|^2 
	+ \alpha \lambda \| \bfu^*_+ \|_0 + \epsilon_k. 
\end{eqnarray*}
Taking the upper limit on both sides of the above inequality, together with the facts $\lim_{k \to \infty} \bfu^k = \bfu^*$ and $\epsilon_k \rightarrow 0$,  leads to
\begin{equation} \notag
	\limsup_{k \to \infty} \| \bfu^k_+ \|_0 \leq \| \bfu^*_+ \|_0.
\end{equation}
Taking the lower semicontinuity of $\|(\cdot)_+\|_0$ into account, we have
\[
\lim_{k \to \infty} \| \bfu^k_+ \|_0 = \| \bfu^*_+ \|_0.
\]
By the definition of the Lyapunov function, we obtain
\begin{eqnarray*}
	\lim_{k \to \infty} \mathcal{V}_{k} 
	&=& \lim_{k \to \infty}f(\bfx^{k}) + \lambda \| \bfu^k_+ \|_0 
	+ \langle \bfy^k, A\bfx^k + \bfb - \bfu^k \rangle 
	+ \frac{\rho}{2} \| A\bfx^k + \bfb - \bfu^k \|^2  \\
	&&+ \frac{\beta}{2} \| \bfx^k - \bfx^{k-1} \|^2   \\
	&=&f(\bfx^{*})+ \lambda \| \bfu^*_+ \|_0 + \langle \bfy^*, A\bfx^* + \bfb - \bfu^* \rangle + \frac{\rho}{2} \| A\bfx^* + \bfb - \bfu^* \|^2  \notag \\ 
	&=& \mathcal{V}_{\rho,\beta}( \bfx^*, \bfu^*, \bfy^*, \bfx^* )
	= f(\bfx^{*})+ \lambda \| \bfu^*_+ \|_0 \quad (\mbox{using} \ A\bfx^* + \bfb = \bfu^*).
\end{eqnarray*}
Due to $\lim_{k \to \infty} \bfx^k = \bfx^*$, taking limit as $k \to \infty$ in $\R_2$ of  (\ref{Residual-error}) leads to 
$\lim_{k \to \infty} \| \bfu^{k+1} _{ \Gamma_{k+1} } \| = 0$.
Considering that $\| \bfu^*_{ \Gamma_{k+1} } \| \leq \| ( \bfu^{k+1} - \bfu^*)_{ \Gamma_{k+1} } \| + \| \bfu^{k+1}_{ \Gamma_{k+1} } \|\leq \| \bfu^{k+1} - \bfu^* \| + \| \bfu^{k+1}_{ \Gamma_{k+1} } \|$ holds, taking limits on both sides of this inequality leads to $\lim_{k \to \infty} \|  \bfu^*_{ \Gamma_{k+1} } \| = 0$.	Since $\Gamma_{k+1} \subseteq [m]$ and $\bfu^*$ is a fixed vector, the claim in (\ref{u^*ga0}) holds. 
We complete the proof.
\hfill $\Box$

%%%%%%%%%%%%%%%%%%%%%%%%%%%%%%%%%%%%%%%%%%%%%%%%%%%%%%%%%%%%%
\subsubsection{Review of KKT solution map}

Not only does the SOSC imply the whole sequence convergence of
$\{(\bfx^k, \bfu^k, \bfy^k) \}$ to a KKT point $(\bfx^*, \bfu^*, \bfy^*)$ of
(\ref{COP-Constrained}), it also implies an upper Lipschitz continuity of
the KKT point of a slightly perturbed problem of (\ref{COP-Constrained}) near
$(\bfx^*, \bfu^*, \bfy^*)$. 
We show this result in terms of the smooth problem (\ref{Smooth-P}).

It follows from Prop.~\ref{Prop-P}(ii) that
$(\bfx^*, \bfu^*, \bfy^*_{\I^*_-}, \bfy^*)$ is a KKT point of the problem (\ref{Smooth-P})
and the KKT conditions are as follows.
\be \label{KKT-O}
\left\{
\begin{array}{l}
	\frac{\partial_\bfx \cL^\#(\bfx^*, \bfu^*, \bfy^*_{\I^*_-}, \bfy^*) }{\partial \bfx} = 0 , \quad 
	\frac{\partial_\bfu \cL^\#(\bfx^*, \bfu^*, \bfy^*_{\I^*_-}, \bfy^*) }{\partial \bfu} = 0
	\\ [1ex]
	\bfu^*_{\I^*_-} \le 0, \quad \bfy^*_{\I^*_-} \ge 0, \quad 
	\langle \bfu^*_{\I^*_-}\;, \bfy^*_{\I^*_-} \rangle = 0 \\ [1ex]
	A\bfx^* + \bfb - \bfu^* = 0,
\end{array}
\right .
\ee 
where $\cL^\#(\bfx, \bfu, \bfs, \bfy)$ is the Lagrange function of (\ref{Smooth-P}) defined in (\ref{Smooth-Lagrange})
%\[
%  \cL^\#(\bfx, \bfu, \bfs, \bfy) :=
%  f(\bfx) + \langle \bfs, \; \bfu_{\I^*_-} \rangle + \langle \bfy,\; A\bfx - \bfb - \bfu \rangle
%\]
and $(\bfy^*_{\I^*_-}, \bfy^*)$ is the Lagrange multiplier at $(\bfx^*, \bfu^*)$.
Let $\M := \M(\bfx^*, \bfu^*)$ be the set of all Lagrange multipliers associated with
$(\bfx^*, \bfu^*)$.

Now let us consider a canonically perturbed KKT system of (\ref{KKT-O}):
\be \label{KKT-Perturbed}
\left\{
\begin{array}{l}
	\frac{\partial_\bfx \cL^\#(\bfx, \bfu, \bfs, \bfy) }{\partial \bfx} = \bfp_1, \quad 
	\frac{\partial_\bfu \cL^\#(\bfx, \bfu, \bfs, \bfy) }{\partial \bfu} =  \bfp_2
	\\ [1ex]
	\bfu_{\I^*_-} \le \bfp_3, \quad \bfs \ge 0, \quad 
	\langle \bfs, \; \bfu_{\I^*_-} - \bfp_3 \rangle = 0 \\ [1ex]
	A\bfx + \bfb - \bfu = 0,
\end{array}
\right .
\ee 
where $\bfp_i$, $i=1, \ldots, 4$ are given perturbation vectors of compatible dimensions.
Denote $\bfp := [\bfp_1; \bfp_2; \bfp_3; \bfp_4]$ and the solution of (\ref{KKT-Perturbed}) by $(\bfx(\bfp), \bfu(\bfp), \bfs(\bfp), \bfy(\bfp))$
(known as the KKT solution map) to indicate
its dependence on $\bfp$. Denote $\zeta := (\bfx, \bfu)$. It is easy to verify that
\[
\nabla^2_{\zeta, \zeta} \cL^\# (\bfx, \bfu, \bfs, \bfy)
= \left[
\begin{array}{cc}
	\nabla^2 f(\bfx) & 0 \\ 
	0 & 0
\end{array} 
\right].
\]
The (classical) second-order sufficient condition (see, e.g., \cite[Eq. (1.6)]{fernandez2010sharp}) for the problem (\ref{Smooth-P}) at $(\bfx^*, \bfu^*)$ reduces to condition (\ref{SOSC}).
By applying \cite[Lemma 2]{fischer2002local} to (\ref{KKT-Perturbed}) under SOSC (\ref{SOSC}), there exists a constant $\tau_0>0$ such that
\be \label{KKT-ErrorBound}
\| (\bfx(\bfp), \bfu(\bfp))  - (\bfx^*, \bfu^*) \|
+ \mbox{dist} 
\Big( (\bfs(\bfp), \bfy(\bfp) ), \M  \Big) 
\le \tau_0 \| \bfp\|
\ee 
for any $\bfp$ close enough to the zero vector, where $\mbox{dist}(\bfx, \Omega)$ denotes the Euclidean distance between $\bfx$ and a closed set $\Omega$.
We will see the perturbed KKT system (\ref{KKT-Perturbed}) will enter our analysis below.

\begin{lemma}[Primal Error Bound] \label{Lemma-Primal-Errorbound}
	There exists an index $k_\nu$ and a constant $\tau>0$ such that
	\[
	\| \bfx^{k+1} - \bfx^* \| 
	\le \tau \Big(
	\| \bfx^{k+1} - \bfx^k \| 
	+
	\| \bfx^k - \bfx^{k-1} \|
	\Big) \quad \forall \ k \ge k_\nu .
	\]	
\end{lemma}

{\bf Proof.}
We are going to use (\ref{KKT-ErrorBound}) to prove this result.
For this reason, we shall construct a sequence
$\bfp^{k+1} = [\bfp^{k+1}_1; \bfp^{k+1}_2; \bfp^{k+1}_3; \bfp^{k+1}_4]$
satisfying $\bfp^{k+1} \rightarrow 0$, and construct a corresponding sequence
$\bfw^{k+1} := [\bfx(\bfp^{k+1}); \bfu(\bfp^{k+1}); \bfs(\bfp^{k+1}); \bfy(\bfp^{k+1})]$.
Note that we work with the index $(k+1)$ instead of $k$. 
We prove $(\bfp^{k+1}, \bfw^{k+1})$ satisfy the perturbed system (\ref{KKT-Perturbed}). 

{\bf Step 1: Choosing} $\bfp^{k+1}$. Let
\[
\bfx({\bfp^{k+1}}) := \bfx^{k+1}, \quad
\bfu({\bfp^{k+1}}) := \bfu^{k+1}, \quad
\bfy({\bfp^{k+1}}) := \bfy^{k+1},
\]
and $\bfs({\bfp^{k+1}}) := \bfs^{k+1}$ defined by
\[
s^{k+1}_i = \left\{
\begin{array}{ll}
	y^{k+1}_i, & \ \mbox{if} \ i \in \I^*_- \cap \Gamma_{k+1} \\
	0 ,        & \ \mbox{otherwise}.
\end{array} 
\right .
\]
Let
\begin{eqnarray*}
	\bfp_1^{k+1} &:=& \frac{\partial_\bfx \cL^\# (\bfx^{k+1}, \bfu^{k+1} , \bfs^{k+1}, \bfy^{k+1} )}{\partial\bfx} \\ [1ex]
	\bfp_2^{k+1} &:=& \frac{\partial_\bfu \cL^\# (\bfx^{k+1}, \bfu^{k+1} , \bfs^{k+1}, \bfy^{k+1} )}{\partial \bfu} \\ [1ex]
	\bfp_3^{k+1} &:=& \overline{\bfu}^{k+1}_{\I^*_-} \ \ \mbox{with}\ \
	\overline{u}^{k+1}_i  = \left\{  \begin{array}{ll}
		u^{k+1}_i, &\ \mbox{if} \ i \in \I^*_- \cap \Gamma_{k+1} \\
		0 ,        & \ \mbox{otherwise} 
	\end{array}  \right .
	\\ 
	\bfp_4^{k+1} &:=& A \bfx^{k+1} + \bfb - \bfu^{k+1}.
\end{eqnarray*} 
In order for the constructed $\bfw^{k+1}$ and $\bfp^{k+1}$ to satisfy (\ref{KKT-Perturbed}), it is sufficient to prove
(i)
$
\bfu^{k+1}_{\I^*_-} \le \bfp_3^{k+1}
$
and (ii)
$
\bfs^{k+1} \ge 0
$
for all sufficiently large $k$.

{\bf Proof of} (i): 
By the definition of $\I^*_-$, we know that $\bfu^*_{\overline{\I}^*_-   } > 0$ and thus $\bfu^{k+1}_{\overline{\I}^*_-} > 0$ for all sufficiently large $k$. From (\ref{k+0iden}) and $\| (\cdot)_+ \|_0 \in [m]$, we have
\begin{align} \label{u+0eq}
	\| \bfu^{k+1}_+ \|_0 = \| \bfu^*_+ \|_0 \quad \mbox{for all $k$ sufficiently large}.
\end{align}
This combines with the fact $\bfu^*_{\overline{\I}^*_-} > 0$ leads to 
\begin{align} \label{uS^*-leq}
	\bfu^{k+1}_{\I^*_-} \leq 0 \quad \mbox{for all $k$ sufficiently large}.
\end{align}
By the definition of $\bfp_3^{k+1}$, (i) follows from (\ref{uS^*-leq}).

{\bf Proof of} (ii): Suppose (ii) is not true.
Then there exists $i_0 \in \I^*_-\cap \Gamma^{k+1}$ such that $y^{k+1}_{i_0} < 0$.
This, together with (\ref{uS^*-leq}) imply $u^{k+1}_{i_0} + \alpha y^{k+1}_{i_0} < 0$, and thus $i_0 \in \OG_{k+1}$ for all $k$ sufficiently large. This yields a contradiction because of $i_0 \in \Gamma^{k+1}$.

{\bf Step 2: Proving} $\bfp^{k+1} \rightarrow 0$. 
By simple calculation, we have
\[
\bfp^{k+1}_1 = -\nabla f(\bfx^{k+1}) + A^\top \bfy^{k+1}
\stackrel{(\ref{gk-gradient})}{=}
\nabla_\bfx g_k(\bfx^{k+1}, \bfu^{k+1}) - \mu (\bfx^{k+1} - \bfx^k) 
\]
and
\[
\bfp^{k+1}_2 = \left[
\begin{array}{c}
	\bfs^{k+1} - \bfy^{k+1}_{\I^*_-} \\
	- \bfy^{k+1}_{\overline{\I}^*_-}
\end{array} 
\right], \quad \bfp^{k+1}_4 = \frac{\bfy^{k+1} - \bfy^k }{\rho}.
\]
The following bounds on $\bfp^{k+1}_i$, $i=1, \ldots, 4$ are understood for sufficiently large $k$.
It follows from (\ref{Residual-error}) on the residual function $\R_1$ that
\[
\| \bfp^{k+1}_1 \|
\le \mu \| \bfx^{k+1} - \bfx^k \| + \| \nabla_\bfx g_k(\bfx^{k+1}, \bfu^{k+1})\|
\le (\mu + c_1) \| \bfx^{k+1} - \bfx^k \| .
\]
Regarding $\bfp^{k+1}_2$, we have
\[
\| \bfs^{k+1} - \bfy^{k+1}_{\I^*_-} \| = \| \bfy^{k+1}_{\I^*_- \cap \OG_{k+1}  }\|
\le \| \bfy^{k+1}_{\OG_{k+1}}\|.
\]
Since $\bfu^*_{\overline{\I}^*_- } >0$ and $\bfu^*_{\Gamma_{k+1} } =0$ for all $k$ sufficiently large in (\ref{u^*ga0}), we have $\overline{\I}^*_- \subseteq \OG_{k+1}$.
Consequently, we have
\[
\| \bfy^{k+1}_{\overline{\I}^*_- } \| \le \| \bfy^{k+1}_{\OG_{k+1}} \|. 
\]
Therefore, we have for $\bfp^{k+1}_2$
\[
\| \bfp^{k+1}_2 \| \le \| \bfs^{k+1} - \bfy^{k+1}_{\I^*_-} \| + 
\| \bfy^{k+1}_{\overline{\I}^*_- } \| \le 
2  \| \bfy^{k+1}_{\OG_{k+1}} \| \le \frac{2c_2}{\alpha} \| \bfx^{k+1} - \bfx^k \|^2,
\]
where 
the last inequality used the bound on $\R_2$ in (\ref{Residual-error})
and the fact $\bfy^{k+1} = - \nabla_\bfu g_k(\bfx^{k+1}, \bfu^{k+1})$.

For $\bfp^{k+1}_3$, we have
\[
\| \bfp^{k+1}_3\| =
\| \overline{\bfu}^{k+1}_{{\I^*_-}} \| = \| {\bfu}^{k+1}_{{\I^*_-}\cap \Gamma_{k+1}} \|
\le \| {\bfu}^{k+1}_{\Gamma_{k+1}} \| \le 
{ c_2} \| \bfx^{k+1} - \bfx^k \|^2,
\]
where the last inequality used the bound on $\R_2$ in (\ref{Residual-error}).
Finally, it follows from (\ref{eq3.5}) that
\[
\| \bfp^{k+1}_4\| = 
\| (\bfy^{k+1} - \bfy^k)/\rho \| \le (c_3/\rho) \| \bfx^{k+1} - \bfx^k \| 
+ (c_4/\rho) \| \bfx^k - \bfx^{k-1} \|.
\]
Given the fact $\bfx^k \rightarrow \bfx^*$, those bounds on $\bfp^{k+1}_i$, $i=1, \ldots, 4$ above show that $\bfp^{k+1} \rightarrow 0$.

{\bf Step 3: Applying the bound in (\ref{KKT-ErrorBound})}.
There exists a constant $\tau_0>0$ such that
\[
\| \bfx^{k+1} - \bfx^* \| \le \tau_0 \| \bfp^{k+1}\| 
\le \tau_0 \sum_{i=1}^4 \| \bfp^{k+1}_i \|.
\]
Substituting the bounds on each $\| \bfp^{k+1}_i\|$, we get
\begin{eqnarray*}
	&&\| \bfx^{k+1} - \bfx^* \| \\ 
	&\leq& 
	\tau_0 \Big( (\mu + c_1 + \frac{c_3}{\rho}) \| \bfx^{k+1} - \bfx^k \| 
	+ (\frac{2c_2}{\alpha} +c_2) \| \bfx^{k+1} - \bfx^k \|^2 
	+ \frac{c_4}{\rho} \| \bfx^k - \bfx^{k-1} \| \Big) \\
	&\leq & \tau_0 \Big(  (\mu + c_1 + c_3/\rho + c_2) \| \bfx^{k+1} - \bfx^k \| 
	+ (c_4/\rho) \| \bfx^k - \bfx^{k-1} \| \Big)  \\
	&\leq & \tau \Big( \| \bfx^{k+1} - \bfx^k \| + \| \bfx^k - \bfx^{k-1} \| \Big),
\end{eqnarray*}
where the second inequality used the fact $(2/\alpha + 1)\| \bfx^{k+1} - \bfx^k\| <1$ for 
$k$ sufficiently large, and $\tau := \tau_0 (\mu + c_1 + c_3/\rho + c_2 + c_4/\rho)$.
\hfill $\Box$

%%%%%%%%%%%%%%%%%%%%%%%%%%%%%%%%%%%%%%%%%%%%%%%%%%%%%%%%%%%%
\subsubsection{R-linear convergence results} 

The convergence rate estimation in this part will involve some constants. 
Since $\{\bfx^k, \bfu^k, \bfy^k\}$ is bounded, there exists a constant $M>0$ such that
\[
\| (\bfx^k, \bfu^k, \bfy^k  ) \| \le M, \quad k=0, 1,\ldots, .
\]
This bound and the fact $\bfu^*_{\Gamma_{k+1} }=0$ in (\ref{u^*ga0}) imply that for $k$ sufficiently large
\begin{align}
	&| \langle \bfy^{k+1}_{\Gamma_{k+1}},\; (\bfu^{k+1} - \bfu^*)_{\Gamma_{k+1}} \rangle | = | \langle \bfy^{k+1}_{\Gamma_{k+1}}, \bfu^{k+1} _{\Gamma_{k+1}} \rangle | 
	\leq M \| \bfu^{k+1} _{\Gamma_{k+1}} \| \mathop{\leq}\limits^{(\ref{Residual-error})}  c_2 M \| \bfx^{k+1} - \bfx^k \|^2,   \notag\\
	& | \langle \bfy^{k+1}_{\OG_{k+1}},\; (\bfu^{k+1} - \bfu^*)_{\OG_{k+1}} \rangle |  
	\leq \| \bfy^{k+1}_{\OG_{k+1}} \|( \| \bfu^{k+1}  \| + \| \bfu^* \|)  \mathop{\leq}\limits^{(\ref{Residual-error})} (2c_2M/\alpha) 
	\| \bfx^{k+1} - \bfx^k \|^2.  \notag
\end{align}
Adding the above inequalities leads to 
\begin{align}
	| \langle \bfy^{k+1}, \bfu^{k+1} - \bfu^* \rangle | 
	\leq &  c_2 M (1+2/\alpha) \| \bfx^{k+1} - \bfx^k \|^2. \label{prod}
\end{align}

We are ready to state the first convergence rate result in terms of Lyapunov function values.

\begin{theorem}[Convergence rate of Lyapunov function sequence] \label{Thm-Lyapunov}
	There exist a sufficiently large index $k_\nu$ and two positive constants $c_\nu>0$ and $q<1$ such that
	\be \label{Lyapunov-rate}
	\V_k - \V_* \le c_\nu q^k , \quad \forall \ k \ge k_\nu.
	\ee
\end{theorem}

{\bf Proof.}
We will derive an induction bound for $\V_{k+1} - \V_*$. It follows from (\ref{k+0iden}) (\ref{u+0eq}), and the update rule (\ref{Multiplier-update}) that for all $k$ sufficiently large, we have
\begin{eqnarray*}
	&& \V_{k+1} - \V_* \\
	&=& f(\bfx^{k+1}) - f(\bfx^*) 
	+ \langle \bfy^{k+1}, \; A \bfx^{k+1} + \bfb - \bfu^{k+1} \rangle 
	+ \frac{1}{2\rho} \| \bfy^{k+1} - \bfy^k \|^2 + \frac{\beta}2 \| \bfx^{k+1} - \bfx^k \|^2	\\
	&=& f(\bfx^{k+1}) - f(\bfx^*) 
	+ \langle A^\top \bfy^{k+1}, \; \bfx^{k+1} - \bfx^* \rangle 
	+ \langle \bfy^{k+1}, \; \bfu^* - \bfu^{k+1} \rangle \\
	&&  + \frac{1}{2\rho} \| \bfy^{k+1} - \bfy^k \|^2 + \frac{\beta}2 \| \bfx^{k+1} - \bfx^k \|^2 ,
\end{eqnarray*}
where the second equation used the fact $\bfu^* = A\bfx^* + \bfb$.
Using (\ref{prod}), (\ref{gk-gradient}) for $A^\top \bfy^{k+1}$, and (\ref{eq3.6}),
we continue to relax the above bound:
\begin{eqnarray*}
	&& \V_{k+1} - \V_* \\
	&\le& f(\bfx^{k+1}) - f(\bfx^*)
	+ \langle \nabla_\bfx g_k(\bfx^{k+1}, \bfu^{k+1} ) - \nabla f(\bfx^{k+1}) 
	- \mu(\bfx^{k+1} - \bfx^k) , \; \bfx^{k+1} - \bfx^* \rangle  \\
	&& + \; c_2M (1+2/\alpha) \| \bfx^{k+1} - \bfx^k \|^2
	+  \frac{1}{2\rho} \Big(
	2c_3^2 \| \bfx^{k+1} - \bfx^k\|^2 + 2c_4 \| \bfx^{k} - \bfx^{k-1}\|^2
	\Big) \\
	&& + \; \frac{\beta}2 \| \bfx^{k+1} - \bfx^k \|^2 \\
	&\le& \underbrace{f(\bfx^{k+1}) - f(\bfx^*) -  \langle  \nabla f(\bfx^{k+1}) , \; \bfx^{k+1} - \bfx^* \rangle}_{\le (\mu/2) \| \bfx^{k+1} - \bfx^* \|^2  }
	+ \underbrace{\langle \nabla_\bfx g_k(\bfx^{k+1}, \bfu^{k+1} ) , \; \bfx^{k+1} - \bfx^* \rangle}_{\le c_1 \|\bfx^{k+1} - \bfx^k\| \| \bfx^{k+1} - \bfx^*\| } \\
	&& + \ \mu \|\bfx^{k+1} - \bfx^k\| \| \bfx^{k+1} - \bfx^*\| \\
	&& + \
	\underbrace{\Big(   c_2M (1+2/\alpha) + c_3^2/\rho + \beta/2
		\Big)}_{:= c_5} \|\bfx^{k+1} - \bfx^k\|^2 + \frac{c_4}{\rho} \|\bfx^{k} - \bfx^{k-1}\|^2 \\
	&\le& \Big( \mu + c_1/2  \Big) \| \bfx^{k+1} - \bfx^*\|^2
	+ \Big(  
	c_5 + (\mu+c_1)/2
	\Big) \|\bfx^{k+1} - \bfx^k\|^2
	+ \frac{c_4}{\rho} \|\bfx^{k} - \bfx^{k-1}\|^2 ,
\end{eqnarray*}
where the first inequality above used the fact that $f$ is $\mu$-weakly convex (as $\mu > \sigma_f$) and the bound on $\R_1$ in (\ref{Residual-error}),
and the third inequality used the relaxation:
$
\|\bfx^{k+1} - \bfx^k\| \| \bfx^{k+1} - \bfx^*\| 
\le (1/2) ( \|\bfx^{k+1} - \bfx^k\|^2 +  \| \bfx^{k+1} - \bfx^*\|^2 )
$.
Applying the error bound in Lemma \ref{Lemma-Primal-Errorbound} to the term
$\| \bfx^{k+1} - \bfx^*\|$ in the above inequality, we simplify and get
\be \label{lyap_rate}
\V_{k+1} - \V_* 
\le c_6 \Big(    
\|\bfx^{k+1} - \bfx^k\|^2 + \|\bfx^{k} - \bfx^{k-1}\|^2 
\Big),
\ee
where $c_6:=2(2\mu + c_1) \tau^2 + c_5 + (\mu + c_1)/2 + c_4/\rho$.
On the other hand, we follow from the sufficient decrease in Lemma~\ref{Lemma-Vk} that
\[
(\mathcal{V}_{k-1} - \mathcal{V}_{*}) - (\mathcal{V}_{k+1} - \mathcal{V}_{*}) = \mathcal{V}_{k-1} - \mathcal{V}_{k} + \mathcal{V}_{k} - \mathcal{V}_{k+1} \ge
\frac{\mu}4 ( \| \bfx^{k+1} - \bfx^k \|^2 + \| \bfx^k - \bfx^{k-1} \|^2  ) .
\]
This bound, together with (\ref{lyap_rate}), yields the induction decrease
\[
\V_{k+1} -\V_* \le \frac{1}{1 + \mu/(4c_6)} \Big(  \V_{k-1} - \V_*  \Big) 
\quad \mbox{for all $k$ sufficiently large}.
\]
Therefore, there must exist a sufficiently large index $k_\nu$ such that
the bound (\ref{Lyapunov-rate}) holds with 
\[
q := \sqrt{\frac{1}{1 + \mu/(4c_6)} }\quad \mbox{and} \quad
c_\nu := (1/q)^{k_\nu} (\V_0 - \V_*).
\]
The proof is completed. \hfill $\Box$

The R-linear convergence rate on the Lyapunov function sequence also translates to 
the iterate sequences $\{\bfx^k\}$, $\{\bfu^k\}$, and $\{\bfy^k\}$, as proved below.

\begin{theorem}[Sequence Convergence Rate] \label{Thm-R-Iterates}
	There exist a sufficiently large index $k^*$ and three positive constants
	$c_x$, $c_u$ and $c_y$ such that for any $k \geq k^*$, we have
	\begin{equation} \notag
		\| \bfx^k - \bfx^* \| \leq c_x \sqrt{q}^k,\quad
		\| \bfy^k - \bfy^* \| \leq c_y \sqrt{q}^k,\quad
		\| \bfu^k - \bfu^* \| \leq c_u \sqrt{q}^k,
	\end{equation}
	where $q$ is from (\ref{Lyapunov-rate}). 
\end{theorem} 

{\bf Proof.}
For any $k \geq k_{\nu}$, it follows from (\ref{lyap_rate}) that 
\be \label{Lyapu_ub}
\sqrt{\mathcal{V}_{k+1} - \mathcal{V}_*} 
% & \leq  \sqrt{c_6(\| \bfx^{k+1} - \bfx^k \|^2 + \| \bfx^k - \bfx^{k-1} \|^2 )} 
\le  \sqrt{c_6} ( \| \bfx^{k+1} - \bfx^k \| + \| \bfx^k - \bfx^{k-1} \| ).
\ee
Since $\sqrt{(\cdot)}$ is an increasing and concave function, we have the following estimation
\begin{align}
	\delta_k :=& \sqrt{\mathcal{V}_k - \mathcal{V}_*} - \sqrt{\mathcal{V}_{k + 1 } - \mathcal{V}_*} \geq\frac{\mathcal{V}_k - \mathcal{V}_{k+1}}{2\sqrt{\mathcal{V}_k - \mathcal{V}_*}} % \notag\\
	\mathop{\geq}\limits^{(\ref{Decreasing-Vk}, \ref{Lyapu_ub})}  % & 
	\frac{(\mu/4) \| \bfx^{k+1} - \bfx^k \|^2}{ 2\sqrt{c_6}(\| \bfx^k - \bfx^{k-1} \| + \| \bfx^{k - 1} - \bfx^{k-2} \|) } \notag,
\end{align}
and thus
\begin{align}
	\| \bfx^{k+1} - \bfx^k \| 
	&\leq \sqrt{\frac{8\sqrt{c_6}}{\mu} \delta_k 
		( \| \bfx^k - \bfx^{k-1} \| + \| \bfx^{k - 1} - \bfx^{k-2} \| )} \notag \\
	& \leq \frac{1}{4} ( \| \bfx^k - \bfx^{k-1} \| + \| \bfx^{k - 1} - \bfx^{k-2} \| ) + \frac{8\sqrt{c_6}}{\mu} \delta_k, ~\forall k \geq k_\nu. \label{recu1}
\end{align}
Now let us sum up (\ref{recu1}) from  $k+2$ to $\bar{k}$, where $ k_\nu \leq k \leq \bar{k} - 2$, to get
\begin{align}
	\sum_{l = k+2}^{\bar{k}} \| \bfx^{l+1} - \bfx^l \| 
	&\leq \frac{1}{4}\sum_{l = k+2}^{\bar{k}} \| \bfx^{l} - \bfx^{l-1} \| 
	+ \frac{1}{4} \sum_{l = k+2}^{\bar{k}} \| \bfx^{l-1} - \bfx^{l-2} \|  + \frac{8\sqrt{c_6}}{\mu} \sum_{l = k+2}^{\bar{k}} \delta_{l} \notag  \\
	&\leq \frac{1}{4}\sum_{l = k+2}^{\bar{k}} \| \bfx^{l+1} - \bfx^{l} \|
	+ \frac 14 \| \bfx^{k+2} - \bfx^{k+1} \| \notag \\
	& \ \ + \frac{1}{4}\sum_{l = k+2}^{\bar{k}} \| \bfx^{l+1} - \bfx^{l} \|
	+ \frac 14 \| \bfx^{k+2} - \bfx^{k+1} \| 
	+ \frac 14 \| \bfx^{k+1} - \bfx^{k}\| \notag \\
	& \ \ +  \frac{8\sqrt{c_6}}{\mu} \sum_{l = k+2}^{\bar{k}}  (\sqrt{\mathcal{V}_l - \mathcal{V}_*} - \sqrt{\mathcal{V}_{l + 1 } - \mathcal{V}_*}), \notag \\
	& = \frac{1}{2}\sum_{l = k+2}^{\bar{k}} \| \bfx^{l+1} - \bfx^{l} \|
	+  \frac 12 \| \bfx^{k+2} - \bfx^{k+1} \|
	+ \frac 14 \| \bfx^{k+1} - \bfx^{k} \| \notag \\
	& \ \ +  \frac{8\sqrt{c_6}}{\mu} \sum_{l = k+2}^{\bar{k}}  (\sqrt{\mathcal{V}_l - \mathcal{V}_*} - \sqrt{\mathcal{V}_{l + 1 } - \mathcal{V}_*}).
\end{align}
Thus, using
$
\sqrt{\mathcal{V}_{k+2} - \mathcal{V}_*} \ge \sqrt{\mathcal{V}_{\bar{k}+1} - \mathcal{V}_*},
$
we have
\begin{equation} \notag
	\frac{1}{2}\sum_{l = k+2}^{\bar{k}} \| \bfx^{l+1} - \bfx^l \| 
	\leq  \frac{1}{2} \| \bfx^{k+2} - \bfx^{k+1} \| 
	+ \frac{1}{4} \| \bfx^{k+1} - \bfx^{k} \| + \frac{8\sqrt{c_6}}{\mu} \Big(  \sqrt{\mathcal{V}_{k+2} - \mathcal{V}_*} 
	\Big) ,
\end{equation}
which further leads to
\begin{align}
	\sum_{l = k}^{\bar{k}} \| \bfx^{l+1} - \bfx^l \|	
	&\leq \frac{3}{2} \| \bfx^{k+1} - \bfx^{k} \| + 2 \| \bfx^{k+2} - \bfx^{k+1} \| + \frac{16\sqrt{c_6}}{\mu} \sqrt{\mathcal{V}_{k+2} - \mathcal{V}_*} \notag \\
	& \mathop{\leq}\limits^{(\ref{Decreasing-Vk})} \frac{3}{\sqrt{\mu}} \sqrt{\mathcal{V}_{k} - \mathcal{V}_{k+1}} 
	+ \frac{4}{\sqrt{\mu}} \sqrt{\mathcal{V}_{k+1} - \mathcal{V}_{k+2}} 
	+ \frac{16\sqrt{c_6}}{\mu} \sqrt{\mathcal{V}_{k+2} - \mathcal{V}_*} \notag \\
	& \leq \left(\frac{7}{\sqrt{\mu}} + \frac{16\sqrt{c_6}}{\mu} \right) \sqrt{\mathcal{V}_{k} - \mathcal{V}_*}. \notag
\end{align}
Since the right-hand of the inequality above is independent of $\bar{k}$, we 
let $\bar{k} \rightarrow \infty$ to get
\[ %\begin{equation}
\sum_{l = k}^{\infty} \| \bfx^{l+1} - \bfx^l \| 
\leq \left(\frac{7}{\sqrt{\mu}} + \frac{16\sqrt{c_6}}{\mu} \right) \sqrt{\mathcal{V}_{k} - \mathcal{V}_*}
, \quad\forall k \geq  k_\nu. %\label{ub_w_lyap}
\] %\end{equation}
The inequality above yields
\begin{align}
	\| \bfx^k - \bfx^* \| &\leq \| \bfx^{k+1} - \bfx^k \| + \| \bfx^{k+1} - \bfx^* \| \leq \cdots \leq \sum_{l = k}^{\infty} \| \bfx^{l+1} - \bfx^l \| \notag \\
	&\leq \left(\frac{7}{\sqrt{\mu}} + \frac{16\sqrt{c_6}}{\mu} \right) \sqrt{\mathcal{V}_{k} - \mathcal{V}_*} \label{rate_w1} \\
	& \mathop{\leq}\limits^{(\ref{Lyapunov-rate})} 
	\underbrace{\sqrt{c_\nu} \left({7}/{\sqrt{\mu}} + {16\sqrt{c_6}}/{\mu} \right)}_{:=c_x} q^{\frac{k}{2}} 
	= c_x \sqrt{q}^k, \quad \forall k \geq  k_\nu.\notag
\end{align}
This completes the R-linear convergence on the sequence $\{\bfx^k\}$. 
We now prove the remaining parts.
It follows from (\ref{gk-gradient}) that
\begin{align}
	\nabla_\bfx g_{k-1} (\bfx^k, \bfu^k) = \nabla f(\bfx^k) + \mu ( \bfx^k - \bfx^{k-1} ) + A^\top \bfy^k. \notag
\end{align}
We also note that 
$
\nabla f(\bfx^*) + A^\top \bfy^* = 0.
$
Subtracting the above two equations implies that
\begin{equation}
	\nabla f(\bfx^k)- \nabla f(\bfx^*) + \mu (\bfx^k - \bfx^{k-1}) 
	+ A^\top (\bfy^k - \bfy^*) = \nabla_\bfx g_{k-1} (\bfx^k, \bfu^k), \notag
\end{equation}
and thus we can estimate
\begin{align}
	\gamma\| \bfy^k - \bfy^*\| &\leq 
	\| A^\top ( \bfy^k - \bfy^* ) \| \notag \\ 
	& \leq \| \nabla f(\bfx^k)- \nabla f(\bfx^*) \| 
	+ \mu \| \bfx^k - \bfx^{k-1} \| + \| \nabla_\bfx g_{k-1} (\bfx^k, \bfu^k) \| \notag \\
	\mathop{\leq}\limits^{(\ref{Residual-error})}
	& \ell_f\| \bfx^k - \bfx^* \| + (\mu + c_1) \| \bfx^k - \bfx^{k-1} \|. \notag
\end{align}
Denote $k^*:=k_\nu + 1$. We have the further estimation for $k \geq k^*$,
\begin{align}
	\| \bfy^k - \bfy^*\| 
	&\leq \frac{\ell_f}{\gamma} \| \bfx^k - \bfx^* \| 
	+ \frac{\mu + c_1}{\gamma} \| \bfx^k - \bfx^{k-1} \| \notag \\
	& \mathop{\leq}\limits^{(\ref{rate_w1},\ref{Decreasing-Vk})} 
	\ell_f \left(\frac{7}{\sqrt{\mu}} + \frac{16\sqrt{c_6}}{\mu} \right)
	\sqrt{\mathcal{V}_{k} - \mathcal{V}_*} 
	+ \frac{2(\mu + c_1)}{\gamma\sqrt{\mu}} \sqrt{\mathcal{V}_{k-1} - \mathcal{V}_{k}} \notag\\
	& \leq \left(\frac{7\ell_f}{\sqrt{\mu}} + \frac{16\sqrt{c_6} \ell_f}{\mu} 
	+ \frac{2(\mu + c_1)}{\gamma\sqrt{\mu}} \right)
	\sqrt{\mathcal{V}_{k-1} - \mathcal{V}_{*}} \le c_y \sqrt{q}^k, \notag 
\end{align}
where the last inequality used (\ref{Lyapunov-rate}) and
\[
c_y :=  \left(\frac{7\ell_f}{\sqrt{\mu}} + \frac{16\sqrt{c_6} \ell_f}{\mu} 
+ \frac{2(\mu + c_1)}{\gamma\sqrt{\mu}} \right) \sqrt{c_\nu}/q.
\]
Similarly, for $k \geq k^*$, we can obtain
\begin{align}
	&\| \bfu^k - \bfu^*\|  \notag	 \\
	= & \| - (\bfy^k - \bfy^{k-1}) /\rho+ A(\bfx^k - \bfx^*) \| 
	\leq \frac{1}{\rho} \| \bfy^k - \bfy^{k-1} \| + \| A \| \| \bfx^k - \bfx^* \| \notag\\
	\mathop{\leq}\limits^{(\ref{eq3.5})}& \frac{c_3}{\rho} 
	\| \bfx^{k} - \bfx^{k-1} \| + \frac{c_4}{\rho}\| \bfx^{k-1} - \bfx^{k-2} \| 
	+ \| A \| \| \bfx^k - \bfx^* \| \notag\\
	\mathop{\leq}\limits^{(\ref{rate_w1},\ref{Decreasing-Vk})} & \frac{2c_3}{\rho\sqrt{\mu}} \sqrt{{\mathcal{V}_{k-1} - \mathcal{V}_{k}}} + \frac{2c_4}{\rho\sqrt{\mu}} \sqrt{{\mathcal{V}_{k-2} - \mathcal{V}_{k-1}}} 
	+ \|A\| \left(\frac{7}{\sqrt{\mu}} + \frac{16\sqrt{c_6}}{\mu} \right)
	\sqrt{\mathcal{V}_{k} - \mathcal{V}_*} \notag\\
	\leq & \left[ \|A\| \left(\frac{7}{\sqrt{\mu}} + \frac{16\sqrt{c_6}}{\mu} \right) 
	+ \frac{2(c_3+c_4)}{\rho\sqrt{\mu}} \right] 
	\sqrt{{\mathcal{V}_{k-2} - \mathcal{V}_{*}}} \le c_{u} \sqrt{q}^k, \notag
\end{align}
where the last inequality used the bound (\ref{Lyapunov-rate}) and
\[
c_u := \left[ \|A\| \left(\frac{7}{\sqrt{\mu}} + \frac{16\sqrt{c_6}}{\mu} \right) 
+ \frac{2(c_3+c_4)}{\rho\sqrt{\mu}} \right] \sqrt{c_\nu}/q.
\] 
Overall, for all $k \geq k^*$, we arrive at the desired R-linear convergence rates.
\hfill $\Box$

%%%%%%%%%%%%%%%%%%%%%%%%%%%%%%%%%%%%%%%%%%%%%%%%%%%%%%%%%%%
\section{Numerical Experiments} \label{Section-Numerical}

In this section, we demonstrate the numerical performance of our iNALM on two important applications in support vector machine (SVM) and multi-label classification (MLC) problems. All the experiments were coded in MATLAB2021 and
implemented on a laptop with 32 GB memory and an Intel CORE i7 2.6 GHz CPU.
We report the performance of iNALM on SVM and MLC in separate subsections.

%%%%%%%%%%%%%%%%%%%%%%%%%%%%%%%%%%%%%%%%%%%%%%%%%%%%%%%%%%%%%%
\subsection{Experiments for SVM.} \label{ex_svm} 

As a classical tool for solving binary classification problems, SVM has a large number of variants. Here, we aim at solving the one with 0/1 loss, abbreviated as 0/1-SVM. 
Given a training set of $m$ points $\{ (\bfx_i, z_i) : i \in [m] \}$ with the $i$th sample $\bfx_i \in \mathbb{R}^{n}$ and its last component $[\bfx_i]_n =1$, and the $i$th class label $z_i \in \{ -1,1 \}$, 
we denote $\bfz :=(z_1,\cdots, z_m)^\top \in \mathbb{R}^m$ and 
$X:=[\bfx_1, \bfx_2,\cdots, \bfx_m]^\top \in \mathbb{R}^{m \times n}$. 
Then 0/1-SVM can be seen as a special case of (\ref{COP}) with the following setting
%\textcolor{blue}
{\begin{align}
		f(\bfx) :=  \frac{1}{2} (\sum_{i= 1}^{n-1} x_i^2 + \vartheta x_p^2)  ,~ A := -( \bfz \textbf{1}^\top_n) \odot X, \ \bfb := \textbf{1}_m, \  \lambda > 0, \notag
\end{align}} 
where $\vartheta$ is a positive constant and $\odot$ is the Hadamard product of componentwise multiplication. 
Data comes from the following two examples.

%and in our numerical experiments, we set $c = 1$.
\begin{example}\label{examp1}
	Samples with positive (resp. negative) labels are drawn from 
	the normal distribution
	$N(\mu_1,\Sigma_1)$ (resp. $N(\mu_2,\Sigma_2)$), where the parameters
	$\mu_1, \mu_2, \Sigma_1$ and $\Sigma_2$ can be generated by the following matlab codes
	\begin{align} \notag
		&\texttt{$\mu_1$ = randn(n,1),$\mu_2$ = randn(n,1)}, \\
		&\texttt{$\Sigma_1$ = diag(randn(n,1)),$\Sigma_2$ = diag(randn(n,1))}. \notag
	\end{align}
	We then flip $r$ percentage (noise ratio) of those samples, making them be marked with reverse labels. 
\end{example}

\begin{example} \label{svm_real}
	%	\textcolor{blue}
	{We select the following binary classification datasets described in Table~\ref{svm_realdata}.}
	\begin{table}[htbp]
		%\TABLE
		\centering
		\captionsetup{justification=justified}
		\caption{Real binary classification datasets in different fields.} \label{svm_realdata}	
		\begin{tabular}{cccccc}
			\hline
			Abbreviation  & Dataset                & $m_{tr}$ & $m_{te}$ & $n$     & Domain            \\ \hline
			\texttt{cncr} & Colon$\_$cancer$\tablefootnote{https://jundongl.github.io/scikit-feature/\label{featureselection}}$           & 62       & 0        & 2000    & Biology           \\
			\texttt{alla} & ALLAML$\footref{featureselection}$                 & 72       & 0        & 7129    & Biology           \\
			\texttt{arce} & Arcene$\footref{featureselection}$                 & 200      & 0        & 10000   & Mass Spectrometry \\
			\texttt{glio} & Gliomas\_85$\footref{featureselection}$             & 85       & 0        & 22283   & Biology           \\
			\texttt{pege} & Prostate\_GE$\footref{featureselection}$            & 102      & 0        & 5966    & Biology           \\
			\texttt{news} & News20.binary$\tablefootnote{https://www.csie.ntu.edu.tw/~cjlin/libsvm/\label{libsvm}}$          & 19996    & 0        & 1355191 & Text              \\
			\texttt{dbw1} & Dbworld\_bodies$\tablefootnote{https://archive.ics.uci.edu/ml/datasets.php\label{uci}}$         & 64       & 0        & 4072    & Text              \\
			\texttt{dbw2} & Dbworld\_bodies\_stemmed$\footref{uci}$ & 64       & 0        & 3721    & Text              \\
			\texttt{dext} & Dexter$\tablefootnote{http://clopinet.com/isabelle/Projects/NIPS2003/\label{nips}}$                 & 600      & 0        & 20000   & Text              \\
			\texttt{doro} & Dorothea$\footref{nips}$               & 1150     & 0        & 100000  & Biology           \\
			\texttt{fmad} & Farm\_ads$\footref{uci}$               & 4143     & 0        & 54877   & Text              \\
			\texttt{covb} & Covtype.binary$\footref{libsvm}$         & 581012   & 0        & 54      & Biology           \\
			\texttt{mush} & Mushroom$\footref{libsvm}$               & 8124     & 0        & 112     & Biology           \\
			\texttt{phis} & Phishing$\footref{libsvm}$               & 11055    & 0        & 68      & Computer Security \\
			\texttt{rlsm} & Real-sim$\footref{libsvm}$               & 72309    & 0        & 20958   & Text              \\
			\texttt{leuk} & Leukemia$\footref{libsvm}$               & 38       & 34       & 7129    & Biology           \\
			\texttt{rcvb} & Rcv1.binary$\footref{libsvm}$            & 20242    & 677399   & 47236   & Text              \\
			\texttt{ijcn} & Ijcnn1$\footref{libsvm}$                 & 49990    & 91701    & 22      & Text              \\
			\texttt{a1a}  & a1a$\footref{libsvm}$                    & 1605     & 30956    & 123     & Social            \\
			\texttt{made} & Madelon$\footref{libsvm}$                & 2000     & 600      & 500     & Artificial        \\ \hline
		\end{tabular}{}
	\end{table}
	
	{In particular, the data \texttt{cncr}, \texttt{arce}, \texttt{glio}, \texttt{pege}, \texttt{dbw1}, \texttt{dbw2}, \texttt{phi} and \texttt{ijcn} are preprocessed by a sample-wise and then feature-wise normalization. \texttt{rlsm} and \texttt{made} are scaled to $[-1,1],$ while \texttt{fmad} is feature-wisely scaled to unit norm.}
\end{example} 

%%%%%%%%%%%%%%%%%%%%%%%%%%%%%%%%%%%%%%%%%%%%%%%%%%%%%%%%%%%%%%%
\subsubsection{Parameter setup and benchmark methods}

There are two types of parameters to set up in order to use iNALM.
The first type includes $\lambda$, $\vartheta$, $\rho$ and $\mu$, which are called model parameters that define the problem (\ref{COP-Constrained}). 
Note that $\mu$ is used to regularize the problem and defines the Lyapunov function
$\V_{\rho, \mu}$. 
Their best choice of values in general depends on problem data. We will test three groups of values to show their influence on the performance of iNALM.

The parameters of the second type include $c_1$, $c_2$, $\epsilon_k$, $\gamma$, which are 
used in Alg.~\ref{Alg-GSN} and (\ref{para_rhobeta}). 
They are the algorithmic parameters. We set
\[
c_1 = c_2 = 0.1, \quad \gamma = 0.1\min\{\|\bfa_i\|^2 \; | \ i\in [m]\},
\quad \epsilon_k = \frac{10\lambda\alpha}{k},
\]
where $\bfa_i$ denotes the $i$th row vector of $A$,
and  $\alpha$ and $\beta$ are  defined by  (\ref{para_rhobeta}).

We adopt a stopping criterion from \cite{li2015global} to terminate iNALM
if the sum of successive changes is small:
\begin{align} \notag
	\frac{\| \bfx^k - \bfx^{k-1} \| + \| \bfu^{k} - \bfu^{k-1} \| + \| \bfy^{k} - \bfy^{k-1} \|}{\| \bfx^k \|+\| \bfu^k \| + \| \bfy^k \| + 1} < 10^{-3}.
\end{align}
We mainly use three metrics for evaluating algorithms' performance: 
classification accuracy 
{$\texttt{Acc} := 1 - \| \mbox{sgn}(X\bfx) - \bfz\|_0/m$}, 
CPU time (\texttt{Time}), 
and the number of support vectors (\texttt{nSV}).
In addition, we monitor the violation of the first-order optimality condition
(\ref{P-stat}) by
\begin{align} \notag
	\texttt{ FOC}:=\max\{ \| \nabla f (\bfx^k) + A^\top \bfy^k \|, 
	{\rm dist}(\bfu^k,  {\rm Prox}_{\alpha\lambda h(\cdot)} ( \bfu^k + \alpha \bfy^k )), \| A\bfx^k + \bfb - \bfu^k \| \}.
\end{align}
%The approximated solution of subproblem (\ref{eq2.1}) is computed by MSN, whose parameters are set as  
%\begin{align} \notag
%	Q = \frac{1}{2} \| w \|^2_{tI - \rho A^\top A} - f(w),~ t = 1.
%\end{align}
For numerical comparison, we select three commonly used and efficient SVM solvers: HSVM from the library libsvm \cite{chang2011libsvm}, 
LSVM from \cite{pelckmans2002ls}, 
and RSVM from \cite{wu2007robust}. 
All these solvers adopt linear kernel and the involved parameters 
are set as their default values.

%%%%%%%%%%%%%%%%%%%%%%%%%%%%%%%%%%%%%%%%%%%%%%%%%%%%%%%%%
\subsubsection{Numerical comparison} 

In this part, we report our observation on the influence of the model parameters
$\lambda, \rho$ and $\mu$ on the performance of iNALM both graphically and
numerically by running many instances of the data from the two examples.
Our first study is to run iNALM against a 
simulated dataset with $m = 500, n = 2000$ and $r = 2\%$. 
Here, we set $\vartheta = 1$ and adopt initial point 
$\bfx^0 = \textbf{1}_n$, $\bfu^0 = \bfy^0 = \textbf{0}_n$. 
For each choice of values of $(\lambda, \rho, \mu)$, we plot four graphs:
$\V_{\alpha, \beta}$ vs \texttt{Iteration},
\texttt{FOC} vs \texttt{Iteration},
\texttt{Time} vs \texttt{Iteration},
and
\texttt{nSV} vs \texttt{Iteration}.
The observed trend in those graphs will give a good indication how we best
set the model parameters. We detail our observations below.

\begin{figure}[h]
	\subfigure{
		\begin{minipage}[t]{0.5\linewidth}
			\centering
			\includegraphics[width=2.6in]{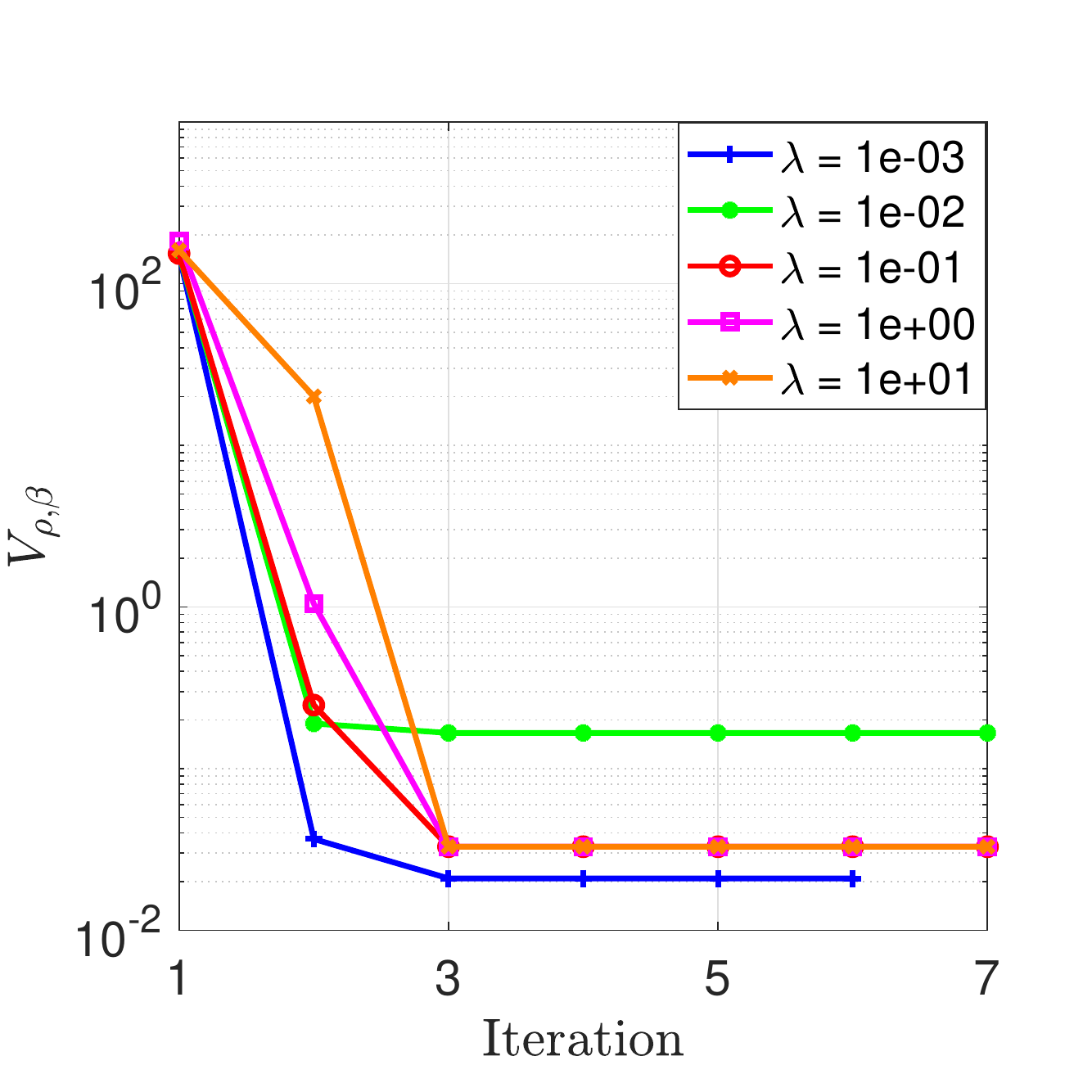}
			%\caption{fig1}
		\end{minipage}%
	}
	\subfigure{
		\begin{minipage}[t]{0.5\linewidth}
			\centering
			\includegraphics[width=2.6in]{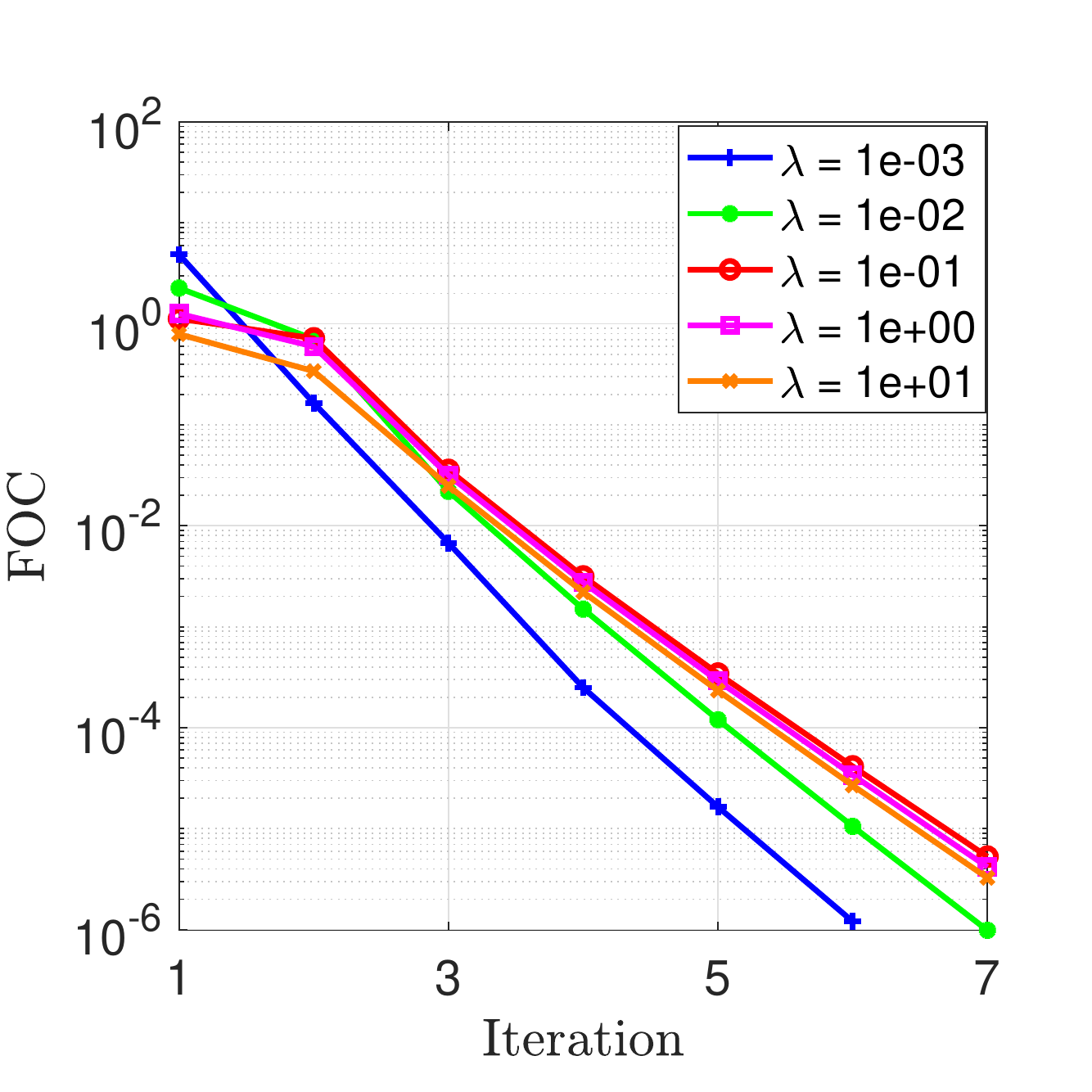}
			%\caption{fig2}
		\end{minipage}%
	}
	%	\subfigure{
	%	\begin{minipage}[t]{0.5\linewidth}
	%		\centering
	%		\includegraphics[width=1in]{lam_acc.eps}
	%		%\caption{fig2}
	%	\end{minipage}%
	%}
	
	\subfigure{
		\begin{minipage}[t]{0.5\linewidth}
			\centering
			\includegraphics[width=2.6in]{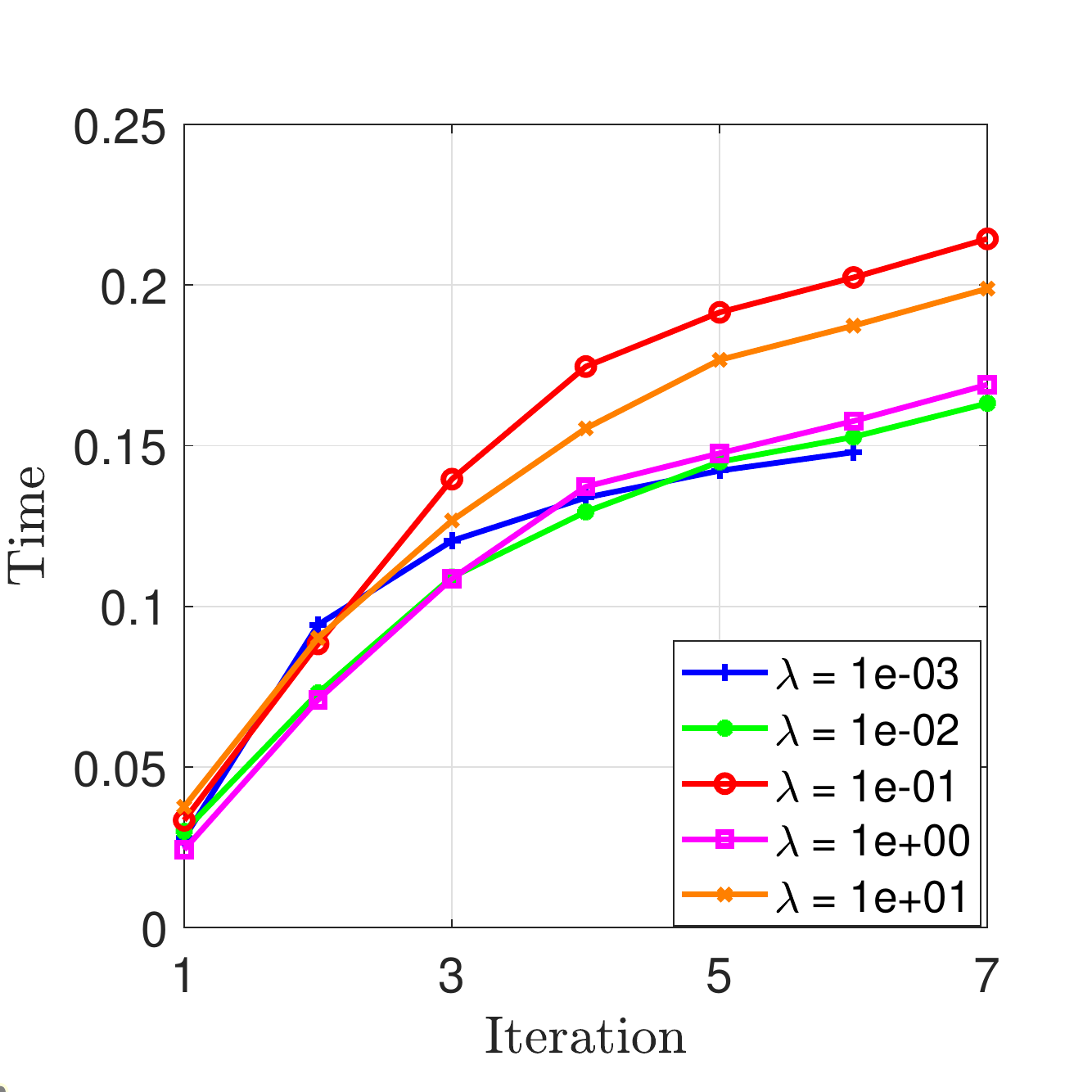}
			%\caption{fig2}
		\end{minipage}
	}%
	\subfigure{
		\begin{minipage}[t]{0.5\linewidth}
			\centering
			\includegraphics[width=2.6in]{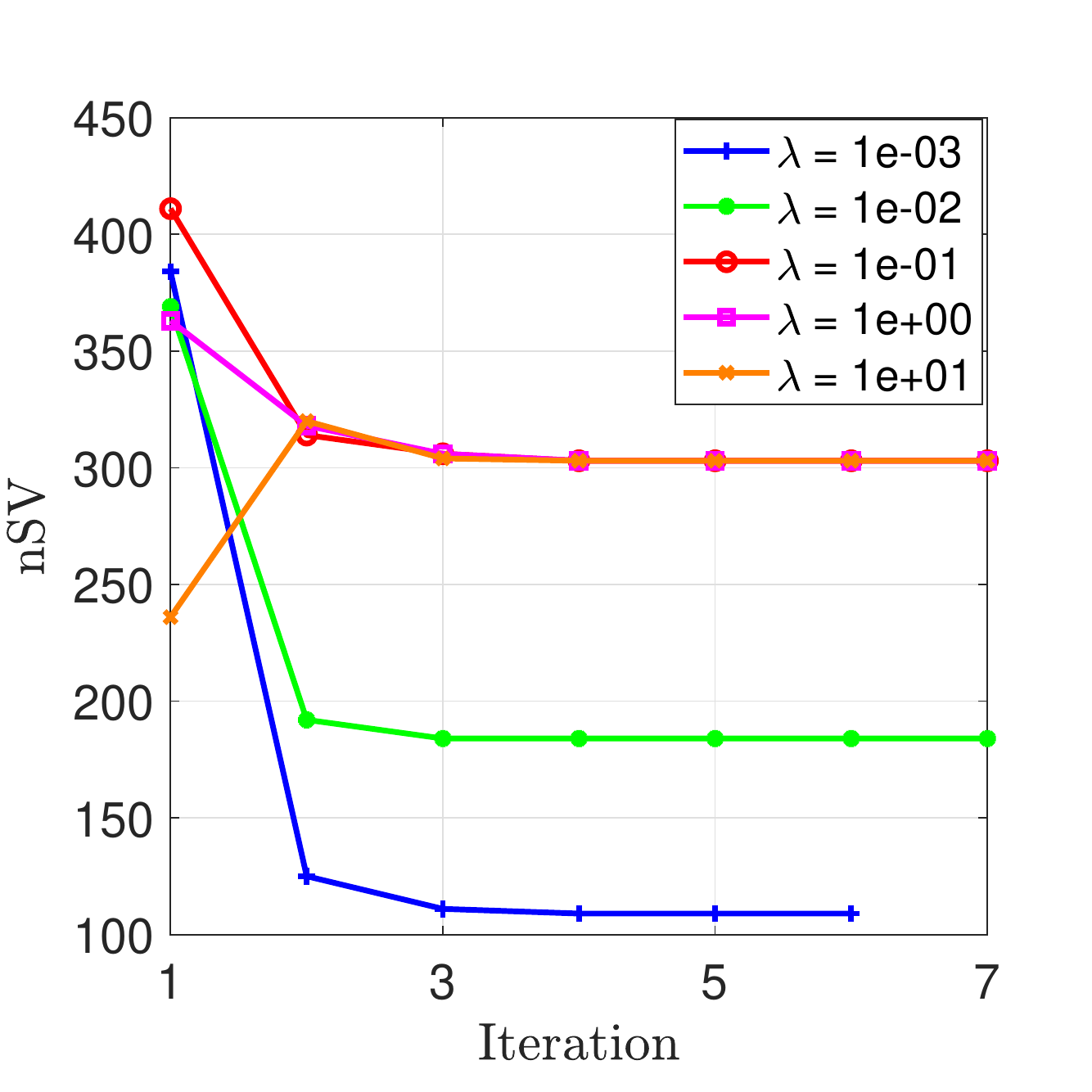}
			%\caption{fig2}
		\end{minipage}
	}
	\caption{Numerical results of iNALM when $\rho = 10^{-2}, \mu = 10^{-2}$ and $\lambda \in \{ 10^{-3}, 10^{-2}, \cdots, 10 \}$. 
		\label{fig1}}
	{}
\end{figure}

{\bf (i) First Test}.
We fix $\rho = 10^{-2}, \mu = 10^{-2}$ and vary $\lambda \in \{ 10^{-3}, 10^{-2}, \cdots, 10 \}$. 
The $\mathcal{V}_{\rho,\beta}$ vs Iteration in Figure \ref{fig1} shows that the Lyapunov function value sequence has a decreasing trend along with iteration.
This is consistent with our theoretical result of Lemma~\ref{Lemma-Vk}.
When $\lambda$ is small, iNALM tends to have a smaller \texttt{nSV} and slightly faster decrease rate on \texttt{FOC}. 

\begin{figure}[h]	
	\subfigure{
		\begin{minipage}[t]{0.5\linewidth}
			\centering
			\includegraphics[width=2.6in]{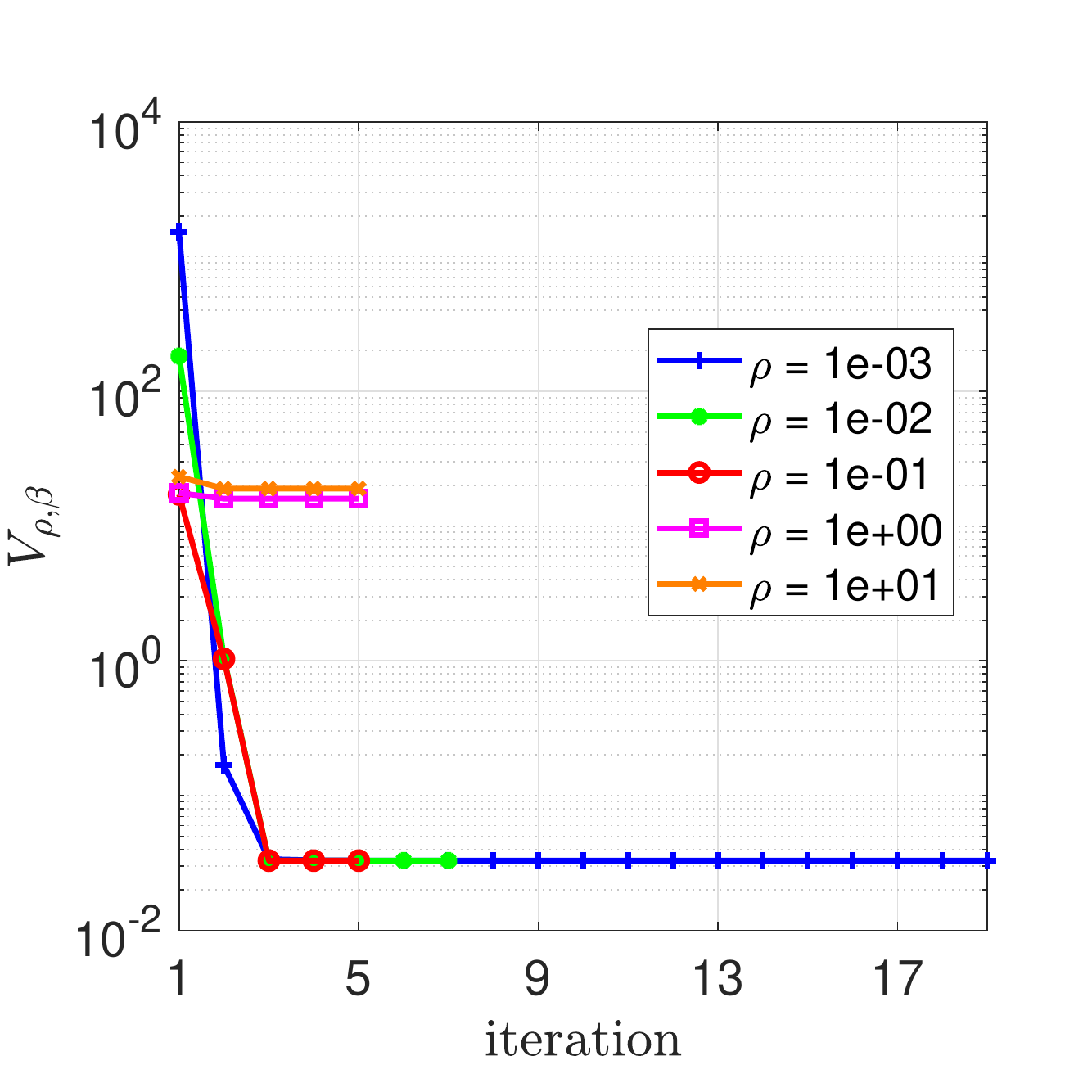}
			%\caption{fig1}
		\end{minipage}%
	}%
	\subfigure{
		\begin{minipage}[t]{0.5\linewidth}
			\centering
			\includegraphics[width=2.6in]{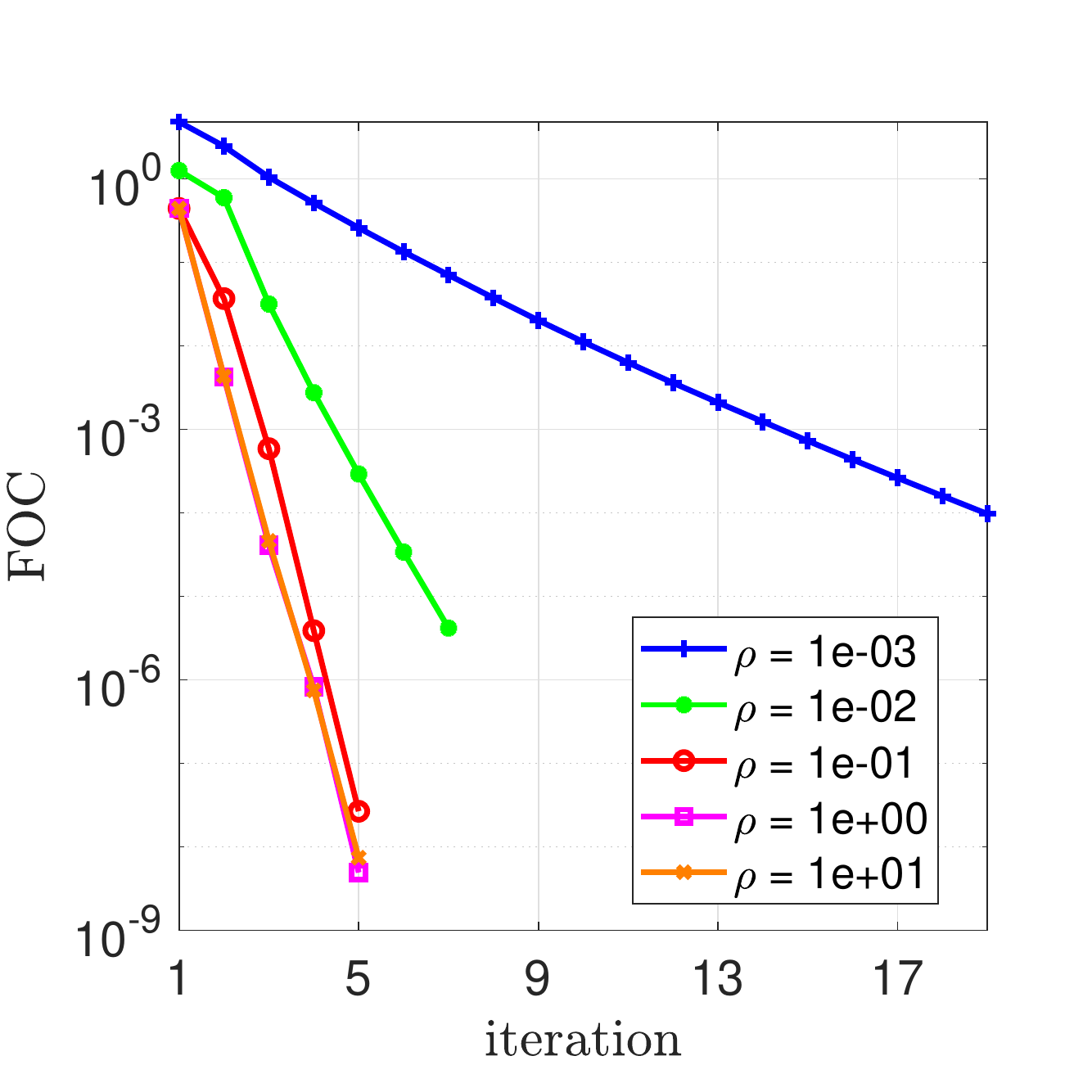}
			%\caption{fig2}
		\end{minipage}%
	}%
	%	\subfigure{
	%	\begin{minipage}[t]{0.334\linewidth}
	%		\centering
	%		\includegraphics[width=2.4in]{rho_acc.eps}
	%		%\caption{fig2}
	%	\end{minipage}%
	%}%
	
	\subfigure{
		\begin{minipage}[t]{0.5\linewidth}
			\centering
			\includegraphics[width=2.6in]{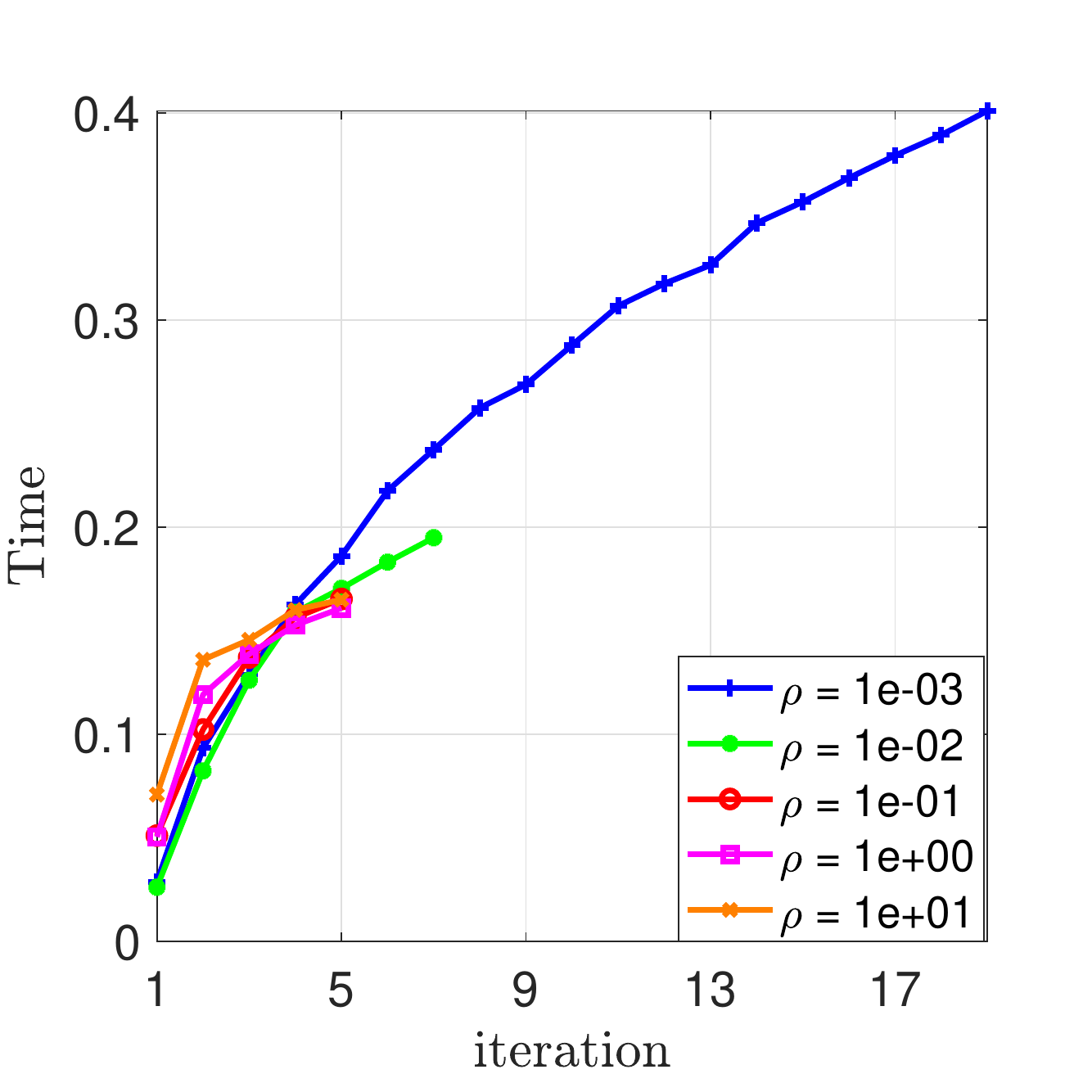}
			%\caption{fig2}
		\end{minipage}
	}%
	\subfigure{
		\begin{minipage}[t]{0.5\linewidth}
			\centering
			\includegraphics[width=2.6in]{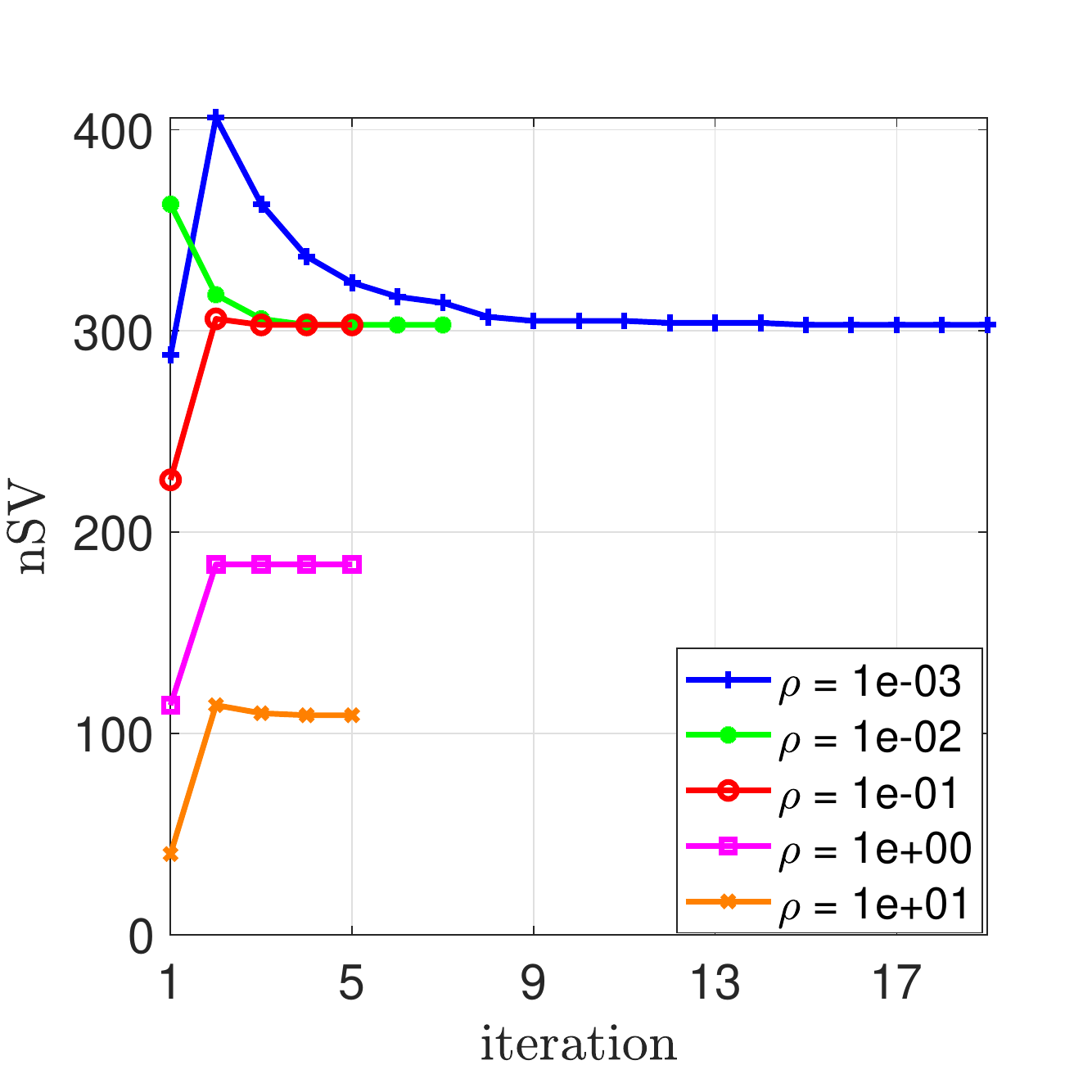}
			%\caption{fig2}
		\end{minipage}
	}%
	\centering
	\caption{Numerical results of iNALM on Example \ref{examp1} when $\lambda = 1, \mu = 10^{-2}$ and $\rho \in \{ 10^{-3}, 10^{-2}, \cdots, 10 \}$.}
	\label{fig2}
\end{figure}

{\bf (ii) Second Test}.
We set $\lambda = 1, \mu = 10^{-2}$ and vary $\rho \in \{ 10^{-3}, 10^{-2}, \cdots, 10 \}$. 
Figure \ref{fig2} demonstrates that $\rho$ has a significant influence on all metrics. For $\rho$ increasing from $10^{-3}$ to $10^{-1}$, 
the decreasing speed of \texttt{FOC} is significantly faster.
This phenomenon can be explained by Thms.~\ref{Thm-Lyapunov} and \ref{Thm-R-Iterates}  
because a larger $\rho$ implies smaller $q$, leading to a faster convergence rate of iNALM. 
The figure of \texttt{Time} is similar when $\rho$ is taken from $10^{-1}$ to $10^{1}$, and it then significantly rises when $\rho$ changes from $1$ to $10$. 
$V_{\rho,\beta}$ tends to get larger as $\rho$ increases, whereas \texttt{nSV} shows an opposite trend. 

\begin{figure}[h]	
	\subfigure{
		\begin{minipage}[t]{0.5\linewidth}
			\centering
			\includegraphics[width=2.6in]{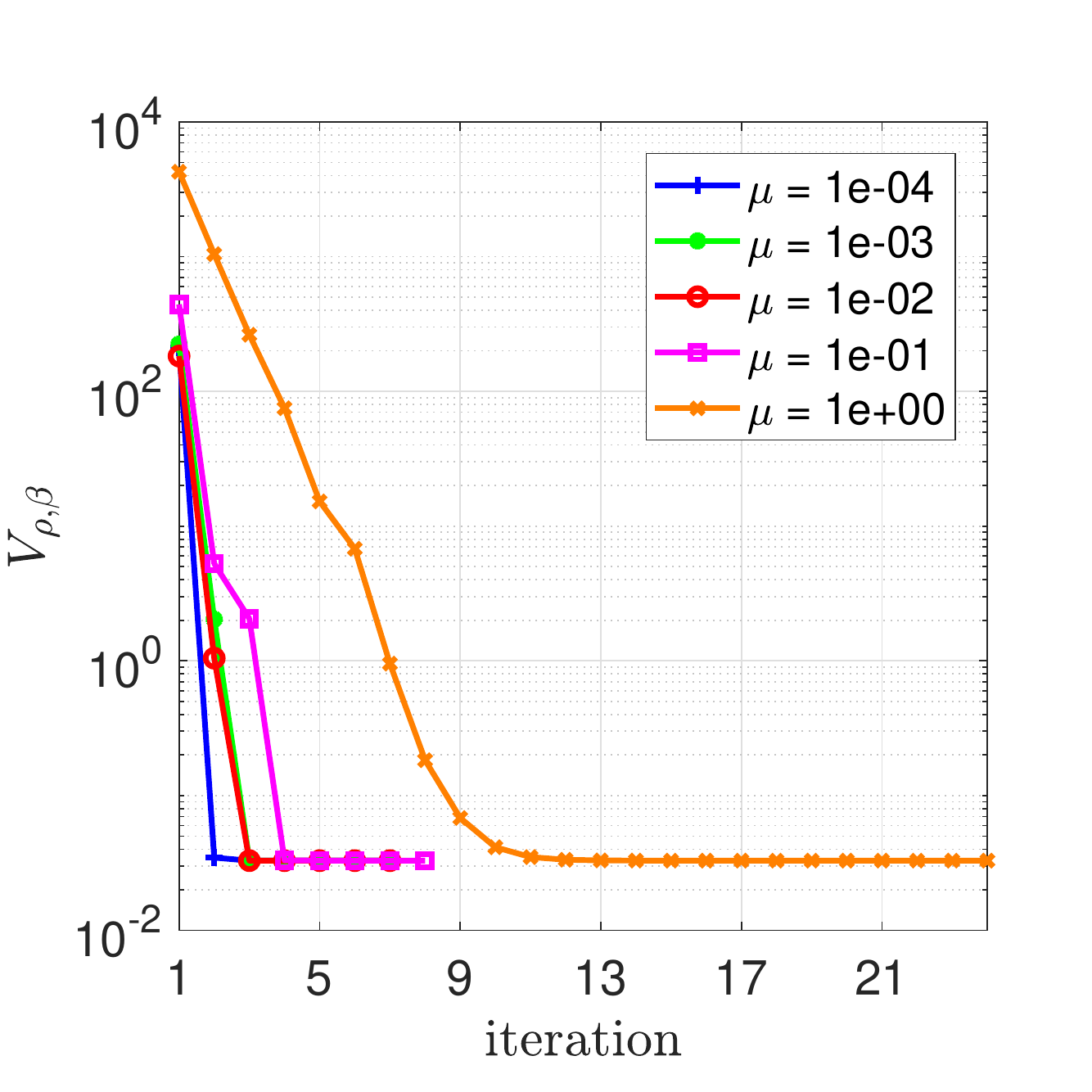}
			%\caption{fig1}
		\end{minipage}%
	}%
	\subfigure{
		\begin{minipage}[t]{0.5\linewidth}
			\centering
			\includegraphics[width=2.6in]{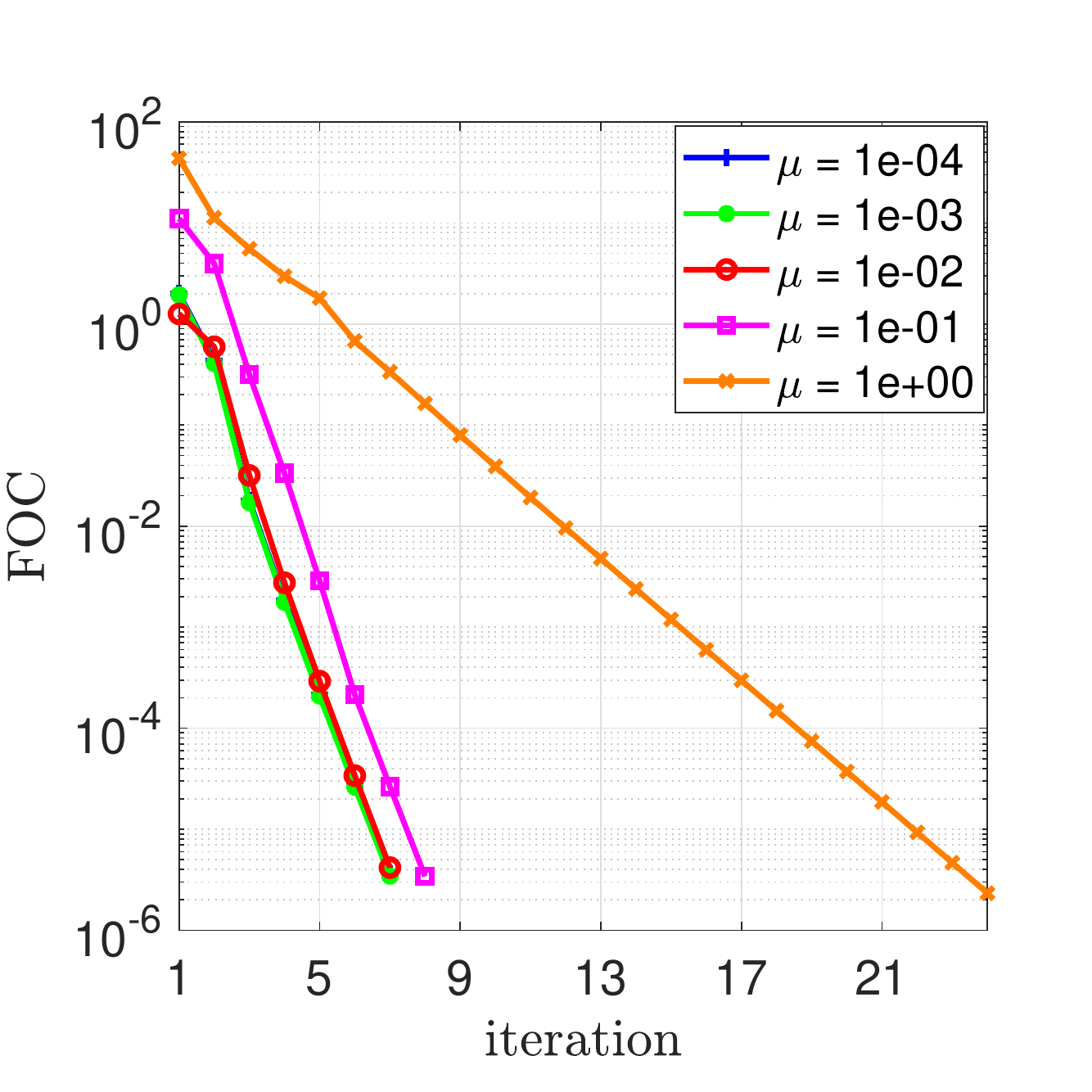}
			%\caption{fig2}
		\end{minipage}%
	}%
	\centering
	
	%		\subfigure{
	%		\begin{minipage}[t]{0.5\linewidth}
	%			\centering
	%			\includegraphics[width=2.5in]{mu_acc.eps}
	%			%\caption{fig1}
	%		\end{minipage}%
	%	}%
	
	\subfigure{
		\begin{minipage}[t]{0.5\linewidth}
			\centering
			\includegraphics[width=2.6in]{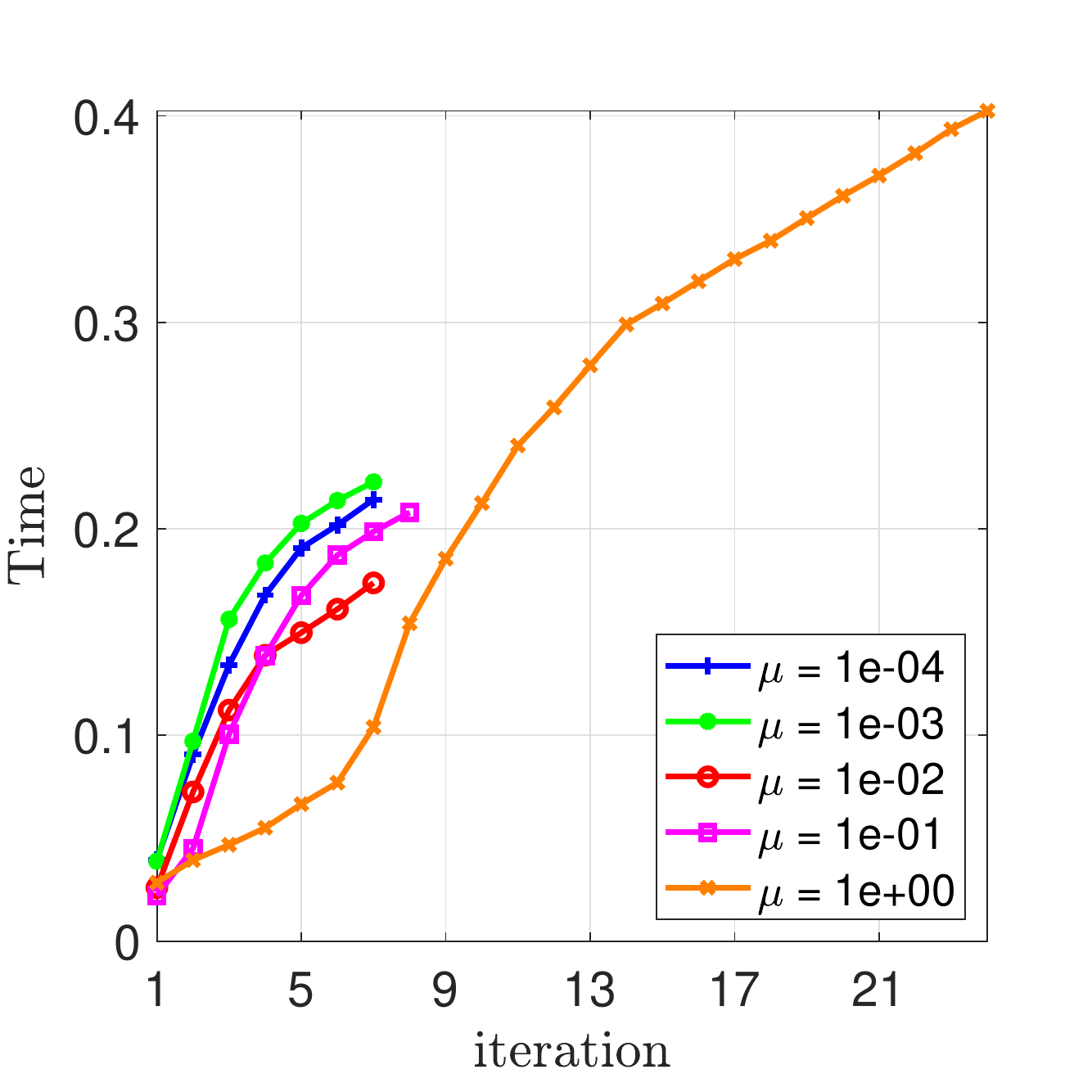}
			%\caption{fig2}
		\end{minipage}
	}%
	\subfigure{
		\begin{minipage}[t]{0.5\linewidth}
			\centering
			\includegraphics[width=2.6in]{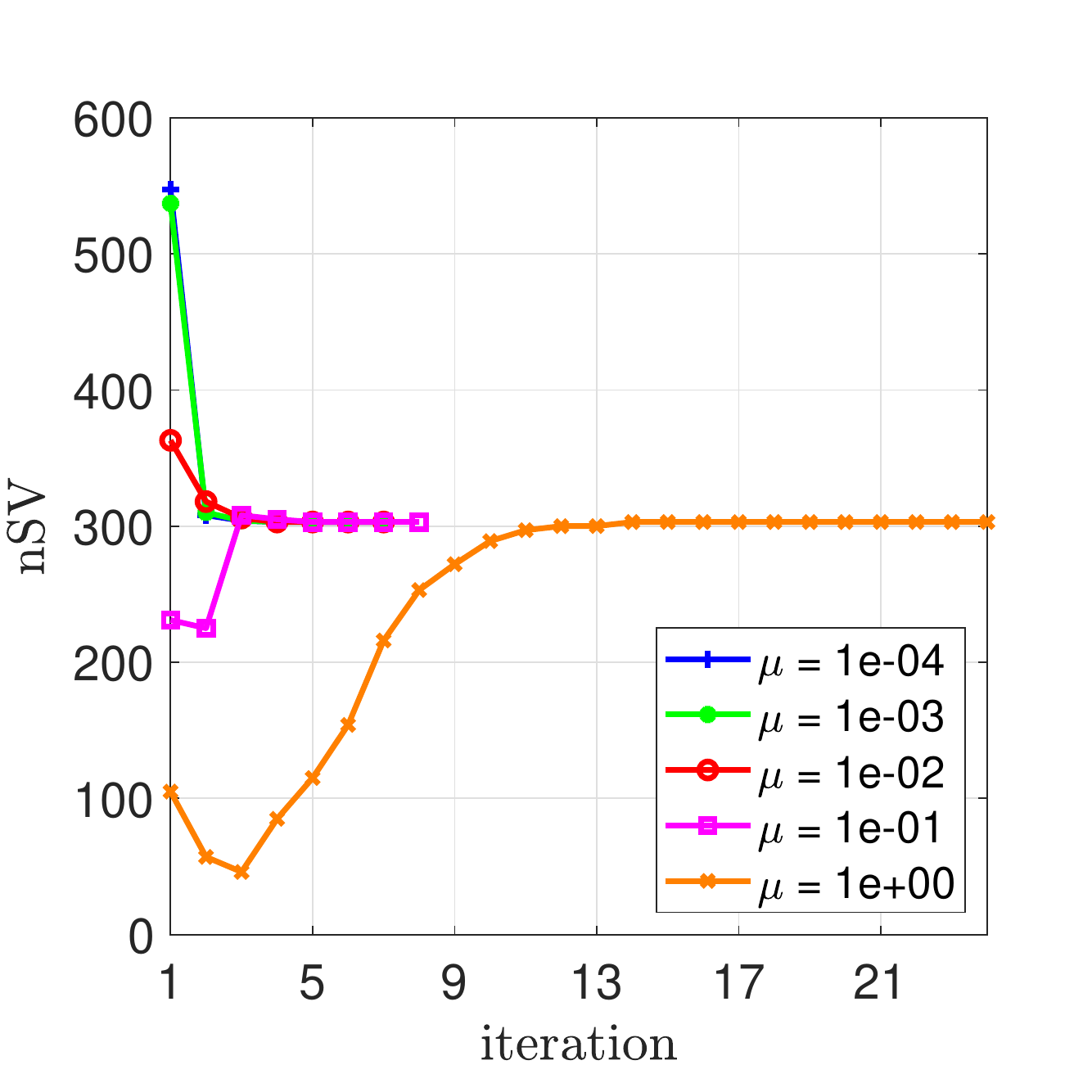}
			%\caption{fig2}
		\end{minipage}
	}%
	\centering
	\caption{Numerical results of iNALM on Example \ref{examp1} when $\lambda = 1, \rho = 10^{-2}$ and $\mu \in \{ 10^{-4}, 10^{-3}, \cdots, 1 \}$.}
	\label{fig3}
\end{figure}

{\bf (iii) Third Test.}
We set $\lambda = 1, \rho = 10^{-2}$ and change $\mu \in \{ 10^{-4}, 10^{-3}, \cdots, 1 \}$. 
The results in Figure \ref{fig3} show that when setting larger $\mu$, decreasing speed of \texttt{FOC} tends to be slower, and thus iNALM is likely to take more iterations to meet the stopping criteria. This is because a larger $\mu$ yields larger $q$, and thus iNALM tends to converge slower according to Thms.~\ref{Thm-Lyapunov} and \ref{Thm-R-Iterates}. Moreover, the parameter $\mu$ has little influence on the ultimate value of the Lyapunov function and \texttt{nSV}. 
Overall, from Figure \ref{fig1}, \ref{fig2} and \ref{fig3}, we can see the descent property of Lyapunov function and linear convergence rate for iNALM. The \texttt{nSV} fluctuates in the first few iterations and then stabilizes  during the last few iterations. 
%Lower nSV helps to decrease Time for 0/1-IALM, but it may lead to a low Acc. Therefore, $\lambda$ and $\rho$ needs to be properly selected to strike a balance between Time and Acc. 

To observe how the algorithms' performance is influenced by the dimension of data, we 
test Example \ref{examp1} with various $m, n$ and $r$. 
Half of the samples are drawn as a training set and the rest of the samples constitute a testing set. 
The number of training and testing sets are denoted as $m_{tr}$ and $m_{te}$ respectively. 
We set $\lambda = 1$, $\mu = 10^{-2}$, $\rho=1$.

\begin{table}[htbp]
	%\TABLE
	\captionsetup{justification=justified}
	\caption{Numerical results of four algorithms on Example \ref{examp1} with different $m_{tr},m_{te},p$ and $r$.} \label{tab1}
	\resizebox{\textwidth}{50mm}
	%{\textcolor{blue}
	{\begin{tabular}{c|cccccccccccc}
			\hline
			& \multicolumn{4}{c|}{\texttt{Acc} $\uparrow$}                                        & \multicolumn{4}{c|}{\texttt{Time} $\downarrow$}                                      & \multicolumn{4}{c}{\texttt{nSV} $\downarrow$}                   \\ \cline{2-13} 
			& \texttt{iNALM} & \texttt{HSVM} & \texttt{LSVM} & \multicolumn{1}{c|}{\texttt{RSVM}} & \texttt{iNALM} & \texttt{HSVM} & \texttt{LSVM} & \multicolumn{1}{c|}{\texttt{RSVM}} & \texttt{iNALM} & \texttt{HSVM} & \texttt{LSVM} & \texttt{RSVM} \\ \cline{2-13} 
			$m_{tr}=m_{te}$ & \multicolumn{12}{c}{$n = 5000, r = 0$}                                 \\ \hline
			1000            & 1.000          & 1.000         & 1.000         & \multicolumn{1}{c|}{1.000}         & 0.125          & 0.997         & 2.956         & \multicolumn{1}{c|}{0.171}         & 157            & 155           & 1000          & 273           \\
			2000            & 1.000          & 1.000         & 1.000         & \multicolumn{1}{c|}{1.000}         & 0.238          & 2.263         & 4.824         & \multicolumn{1}{c|}{0.597}         & 185            & 181           & 2000          & 223           \\
			3000            & 1.000          & 1.000         & 1.000         & \multicolumn{1}{c|}{1.000}         & 0.422          & 3.467         & 8.613         & \multicolumn{1}{c|}{1.690}         & 202            & 194           & 3000          & 285           \\
			4000            & 1.000          & 1.000         & 1.000         & \multicolumn{1}{c|}{1.000}         & 0.625          & 4.731         & 15.33         & \multicolumn{1}{c|}{3.909}         & 210            & 201           & 4000          & 209           \\
			5000            & 1.000          & 1.000         & 1.000         & \multicolumn{1}{c|}{1.000}         & 0.935          & 6.752         & 25.23         & \multicolumn{1}{c|}{7.157}         & 236            & 230           & 5000          & 258           \\ \hline
			$n$             & \multicolumn{12}{c}{$m_{tr} = m_{te} = 2500, r = 0$}                                                                                                                                                                                       \\ \hline
			2000            & 1.000          & 1.000         & 1.000         & \multicolumn{1}{c|}{1.000}         & 0.172          & 0.848         & 4.437         & \multicolumn{1}{c|}{1.006}         & 126            & 124           & 2500          & 209           \\
			4000            & 1.000          & 1.000         & 1.000         & \multicolumn{1}{c|}{1.000}         & 0.235          & 2.109         & 5.780         & \multicolumn{1}{c|}{1.028}         & 174            & 173           & 2500          & 295           \\
			6000            & 1.000          & 1.000         & 1.000         & \multicolumn{1}{c|}{1.000}         & 0.431          & 3.669         & 7.216         & \multicolumn{1}{c|}{0.883}         & 220            & 207           & 2500          & 1175          \\
			8000            & 1.000          & 1.000         & 1.000         & \multicolumn{1}{c|}{1.000}         & 0.898          & 5.208         & 8.426         & \multicolumn{1}{c|}{0.933}         & 236            & 223           & 2500          & 1094          \\
			10000           & 1.000          & 1.000         & 1.000         & \multicolumn{1}{c|}{1.000}         & 1.271          & 6.814         & 9.992         &  \multicolumn{1}{c|}{1.076}              & 242            & 233           & 2500          & 1184           \\ \hline
			$r$             & \multicolumn{12}{c}{$m_{tr} = m_{te} = 5000, n = 100$}                                                                                                                                                                                     \\ \hline
			2\%             & 0.980          & 0.980         & 0.980         & \multicolumn{1}{c|}{0.980}         & 0.086          & 2.557         & 200.3         & \multicolumn{1}{c|}{17.88}         & 40             & 288           & 5000          & 42            \\
			4\%             & 0.960          & 0.960         & 0.960         & \multicolumn{1}{c|}{0.960}         & 0.061          & 3.838         & 200.1         & \multicolumn{1}{c|}{17.74}         & 40             & 514           & 5000          & 42            \\
			6\%             & 0.940          & 0.940         & 0.940         & \multicolumn{1}{c|}{0.940}         & 0.049          & 4.645         & 199.8         & \multicolumn{1}{c|}{17.51}         & 46             & 730           & 5000          & 48            \\
			8\%             & 0.920          & 0.920         & 0.920         & \multicolumn{1}{c|}{0.920}         & 0.079          & 6.572         & 212.3         & \multicolumn{1}{c|}{21.82}         & 35             & 939           & 5000          & 39            \\
			10\%            & 0.900          & 0.900         & 0.900         & \multicolumn{1}{c|}{0.900}         & 0.060          & 7.074         & 234.0         & \multicolumn{1}{c|}{18.03}         & 36             & 1191          & 5000          & 37            \\ \hline
	\end{tabular}}
	
	{``$\uparrow$" means that the larger metric value, the better performance the algortihm has, while ``$\downarrow$" is the opposite situation.}
\end{table}

{\bf (iv) Fourth test.}
We fix $n = 5000, r = 0$ and then change $m_{tr} = m_{te} \in \{ 1000, 2000, \cdots, 5000 \}.$ As shown in Table \ref{tab1}, the \texttt{Time} of iNALM is the shortest, followed by that of RSVM. 
Also,  iNALM has the second smallest \texttt{nSV}. 
In particular, for the case of $m_{tr} = m_{te} = 4000$, iNALM spends about 15\% \texttt{Time} as that of RSVM, and obtains a solution with smaller \texttt{nSV}. 
When we set $m_{tr} = m_{te} = 2500$ and vary $n \in \{ 2000, 4000, \cdots, 10000 \}$, 
iNALM still has advantage on \texttt{Time} spent. 
In the case of $n = 10000$, the \texttt{Time} by iNALM is nearly $1/5$ as much as that of HSVM. 
Although RSVM is also competitive in \texttt{Time}, its \texttt{nSV} is much larger than that of iNALM. 
Finally, we fix $ m_{tr} = m_{te} = 5000$, $n =100$ and 
alter $r \in \{2\%, 4\%, \cdots, 10\%\}$. 
The results show that iNALM has a great advantage on \texttt{Time}, which is almost less than $1\%$ \texttt{Time} of other algorithms. iNALM also performs best on \texttt{nSV}, whereas HSVM fails on this metric. In the case of $r = 10 \%$, \texttt{nSV} of iNALM is approximately $1/30$ as much as that of HSVM. Through those experiments, we observed that iNALM has a fast running speed, and also shows high efficiency and stability for noisy data.

Finally, let us illustrate the numerical experiments on real binary classification datasets. 
%\textcolor{blue}
{For iNALM, we set $\mu = 10^{-2}, \lambda = 1, \rho = 1$ on all the datasets, while we set $\vartheta = 1$ on datasets \texttt{dbw1}, \texttt{dbw2}, \texttt{dext}, \texttt{doro}, \texttt{farm}, \texttt{covb}, \texttt{a1a} and \texttt{made}, and $\vartheta = 10^{-2}$ on other datasets. Particularly, for the data without testing sets, we conduct 5-fold cross validation and record the average results.} Numerical results in Table \ref{tab11} show that iNALM has excellent performance on all metrics. For most datasets, iNALM has the highest \texttt{Acc}. It also spends the shortest time on the datasets with $n \geq 10000$. For the \texttt{rcvb} dataset, LSVM ran for more than 10000s without giving a solution and RSVM ran out of the memory of our laptop, and thus their corresponding results are indicated by ``--". In this case, our iNALM takes less than 4\% \texttt{Time} of HSVM, while finding a solution with higher \texttt{Acc}.

\begin{table}[htbp]
	%\TABLE
	\captionsetup{justification=justified}
	\caption{Numerical results of four algorithms on real binary classification datasets.} \label{tab11}
	\resizebox{\textwidth}{50mm}
	%{\textcolor{blue}
	{\begin{tabular}{c|cccc|cccc|cccc}
			\hline
			& \multicolumn{4}{c|}{\texttt{Acc} $\uparrow$}                   & \multicolumn{4}{c|}{\texttt{Time} $\downarrow$}                & \multicolumn{4}{c}{\texttt{nSV} $\downarrow$}                  \\ \cline{2-13} 
			& \texttt{iNALM} & \texttt{HSVM} & \texttt{LSVM} & \texttt{RSVM} & \texttt{iNALM} & \texttt{HSVM} & \texttt{LSVM} & \texttt{RSVM} & \texttt{iNALM} & \texttt{HSVM} & \texttt{LSVM} & \texttt{RSVM} \\ \hline
			\texttt{cncr} & 0.871          & 0.855         & 0.855         & 0.838         & 0.036          & 0.008         & 1.157         & 0.063         & 36             & 34            & 50            & 43            \\
			\texttt{alla} & 0.961          & 0.923         & 0.936         & 0.923         & 0.031          & 0.035         & 3.328         & 0.036         & 47             & 46            & 58            & 46            \\
			\texttt{arce} & 0.905          & 0.890         & 0.820         & 0.875         & 0.027          & 0.789         & 4.809         & 0.042         & 128            & 112           & 160           & 130           \\
			\texttt{glio} & 0.882          & 0.858         & 0.906         & 0.647         & 0.025          & 0.119         & 9.759         & 0.064         & 48             & 36            & 68            & 68            \\
			\texttt{pege} & 0.940          & 0.930         & 0.921         & 0.930         & 0.029          & 0.052         & 2.889         & 0.038         & 64             & 64            & 82            & 65            \\
			\texttt{news} & 0.962          & 0.968         & --            & --            & 12.05          & 617.6         & --            & --            & 8899           & 7423          & --            & --            \\
			\texttt{dbw1} & 0.871          & 0.671         & 0.538         & 0.671         & 0.028          & 0.020         & 2.234         & 0.037         & 48             & 50            & 51            & 50            \\
			\texttt{dbw2} & 0.767          & 0.675         & 0.5538        & 0.675         & 0.025          & 0.016         & 1.824         & 0.036         & 50             & 50            & 51            & 50            \\
			\texttt{dext} & 0.953          & 0.947         & 0.880         & 0.923         & 0.031          & 0.152         & 69.35         & 0.208         & 418            & 319           & 480           & 480           \\
			\texttt{doro} & 0.928          & 0.927         & --            & 0.928         & 0.443          & 5.288         & --            & 0.824         & 803            & 787           & --            & 920           \\
			\texttt{fmad} & 0.888          & 0.892         & --            & 0.885         & 1.160          & 7.619         & --            & 42.86         & 2611           & 2079          & --            & 3314          \\
			\texttt{covb} & 0.769          & --            & --            & --            & 11.10          & --            & --            & --            & 15465          & --            & --            & --            \\
			\texttt{mush} & 1.000          & 1.000         & --            & --            & 0.084          & 0.130         & --            & --            & 348            & 222           & --            & --            \\
			\texttt{phis} & 0.942          & 0.938         & --            & --            & 0.132          & 2.379         & --            & --            & 194            & 1288          & --            & --            \\
			\texttt{rlsm} & 0.972          & 0.974         & --            & --            & 4.808          & 402.2         & --            & --            & 9305           & 9431          & --            & --            \\
			\texttt{leuk} & 0.882          & 0.824         & 0.853         & 0.824         & 0.076          & 0.020         & 3.740         & 0.098         & 33             & 29            & 38            & 37            \\
			\texttt{rcvb} & 0.962          & 0.963         & --            & --            & 1.950          & 113.3         & --            & --            & 5755           & 4695          & --            & --            \\
			\texttt{ijcn} & 0.951          & 0.926         & --            & --            & 2.586          & 45.74         & --            & --            & 291            & 8389          & --            & --            \\
			\texttt{a1a}  & 0.838          & 0.838         & 0.839         & 0.748         & 0.056          & 0.094         & 12.47         & 304.5         & 170            & 591           & 1605          & 1605          \\
			\texttt{made} & 0.593          & 0.575         & 0.570         & 0.578         & 1.220          & 2.075         & 22.81         & 837.3         & 638            & 1665          & 2000          & 1999          \\ \hline
	\end{tabular}}
	
	{``$\uparrow$" means that the larger metric value, the better performance the algorithm has, while ``$\downarrow$" is the opposite situation. ``--" means that the corresponding algorithm fails to compute a solution within 3 hours or runs out of memory.}
\end{table}

%%%%%%%%%%%%%%%%%%%%%%%%%%%%%%%%%%%%%%%%%%%%%%%%%%%%%%%%%%%%%%%%
\subsection{Experiments on MLC}

Let us consider an MLC problem with $\ell$ classes \cite{zhang2013review}. 
Given $m$ training instances $\bfx_i \in \mathbb{R}^n$ with its last element being
$[\bfx_i]_n = 1$ for all $ i \in [m]$ and 
the class labels $\bfz^{(i)}  \in \{ 1,-1 \}^\ell$ for all $i \in [m]$ with the relevant (resp. irrelevant) classes being $1$ (resp. $-1$), 
we denote $X := [ \bfx_1,\cdots, \bfx_m ]^\top \in \mathbb{R}^{m \times n}$ 
and $Z := [\bfz^{(1)}, \cdots, \bfz^{(\ell)}]^\top \in \mathbb{R}^{m \times \ell} $. Hamming loss, which can be precisely represented by 0/1 loss, is an important metric for evaluating the performance of a multi-label classifier. Meanwhile, sparse regularizers are often used for mitigating the risk of overfitting and selecting relevant features in practice. Here, we consider a sparse regularized MLC model as a special case of 0/1-COP in the following setting
\begin{align}
	&	f(\bfx) :=  \sum_{i=1}^{n\ell} \vartheta_i(x_i^2 + \vartheta_0)^{1/2}, \ \vartheta_i > 0 \ {\rm for~all} \ i \in\{0,1\cdots,\ell m\},  \ \bfb := \textbf{1}_{m\ell}, \ \lambda > 0, \notag \\
	& A := \left[\begin{array}{cccc}
		-( \bfz^{(1)} \textbf{1}^\top_n) \odot X & & &\\
		& -( \bfz^{(2)} \textbf{1}^\top_n) \odot X & & \\
		&  & \ddots &  \\
		&  &  & -( \bfz^{(\ell)} \textbf{1}^\top_n) \odot X 
	\end{array} \right] \in \mathbb{R}^{m\ell \times n \ell}, \notag
\end{align}   
where $f: \mathbb{R}^{n\ell} \to \mathbb{R}$ is a weighted smooth $\ell_1$ regularizer, referring to the smooth $\ell_1$ regularizer proposed by \cite{fountoulakis2016second}, but equipped with weighted constants $\vartheta_i, \forall i \in [m\ell]$. Particularly, we adopt $\vartheta_0 = 10^{-3}$ here.

\begin{example} \label{examp2}
	We generate the feature matrix $X$ and label matrix $Z$ by the following Matlab codes:
	\begin{align}
		&\texttt{X = [randn(m, n-1); ones(m,1)]; W = 2*rand(n, $\ell$) - 1};
		\texttt{Z = sign(X*W);} \notag\\
		& \texttt{J = randperm(m); X = X(J,:); Z = Z(J,:);} \notag
	\end{align}
	In the subsequent numerical comparison for this example, we draw $90 \%$ of samples as a training set and the rest composes a testing set. We denote the number of samples in training (resp. testing) set as $m_{tr}$ ($m_{te}$) and the total number of samples $m = m_{tr} + m_{te}$.
\end{example}

\begin{example} \label{example_mlc_real}
	%\textcolor{blue}
	{We select the following real datasets with higher dimensions from the multi-label classification dataset repository\footnote{https://www.uco.es/kdis/mllresources/}. Particularly, datasets \texttt{sbbc}, \texttt{guar} and \texttt{reut} are preprocessed by sample-wise and then feature-wise normalization, while \texttt{plgo} is feature-wisely scaled to $[-1,1]$.}
	\begin{table}[htbp]
		%\TABLE
		\centering
		\caption{Real multi-label classification datasets with higher dimension.} \label{mlc_realdata}	
		%{\textcolor{blue}
		{\begin{tabular}{ccccccc}
				\hline
				Abbreviation  & Dataset             & $m_{tr}$ & $m_{te}$ & $n$   & $\ell$ & Domain  \\ \hline
				\texttt{sbbc} & 3s-bbc1000          & 239      & 112      & 1185  & 27     & Text    \\
				\texttt{guar} & 3s-guardian1000     & 204      & 98       & 1000  & 6      & Text    \\
				\texttt{genb} & Genbase             & 444      & 218      & 1185  & 27     & Biology \\
				\texttt{reut} & 3s-reuters1000      & 201      & 93       & 1000  & 6      & Text    \\
				\texttt{eugo} & EukaryoteGO         & 5168     & 2584     & 12689 & 22     & Biology \\
				\texttt{gpgo} & GpositiveGO         & 347      & 172      & 912   & 4      & Biology \\
				\texttt{hugo} & HumanGO             & 2050     & 1053     & 9844  & 14     & Biology \\
				\texttt{medi} & Medical             & 659      & 319      & 1449  & 45     & Biology \\
				\texttt{plgo} & PlantGO             & 638      & 331      & 3091  & 12     & Biology \\
				\texttt{rcs1} & Rcv1subset1         & 4045     & 1955     & 47236 & 101    & Text    \\
				\texttt{rcs2} & Rcv1subset2         & 4045     & 1955     & 47236 & 101    & Text    \\
				\texttt{rcs3} & Rcv1subset3         & 4021     & 1979     & 47236 & 101    & Text    \\
				\texttt{rcs4} & Rcv1subset4         & 3997     & 2003     & 47229 & 101    & Text    \\
				\texttt{rcs5} & Rcv1subset5         & 3964     & 2036     & 47235 & 101    & Text    \\
				\texttt{yhbs} & Yahoo\_Business      & 7523     & 3691     & 21924 & 30     & Text    \\
				\texttt{yhet} & Yahoo\_Entertainment & 8569     & 4161     & 32001 & 21     & Text    \\
				\texttt{yhhl} & Yahoo\_Health        & 6158     & 3047     & 30605 & 32     & Text    \\
				\texttt{scen} & Scene               & 1618     & 789      & 294   & 6      & Image   \\
				\texttt{yelp} & Yelp                & 7231     & 3569     & 668   & 8      & Text    \\
				\texttt{hpaa} & HumanPseAAC         & 2053     & 1053     & 440   & 14     & Biology \\
				\texttt{20ng} & 20NG                & 6158     & 3047     & 1006  & 20     & Text    \\
				\texttt{ohsu} & Ohsumed             & 9301     & 4628     & 1002  & 23     & Text    \\
				\texttt{reuk} & Reuters-K500        & 3925     & 1958     & 500   & 103    & Text    \\
				\texttt{slas} & Slashdot            & 2527     & 1222     & 1079  & 22     & Text    \\ \hline
		\end{tabular}}{}
	\end{table}
\end{example}

\subsubsection{Benchmark Methods and Experimental Setup.} 

We select three competitive algorithms in the MLC field. They are ranking support vector machine \cite{elisseeff2001kernel} (Rank-SVM), random $k$-labelsets \cite{tsoumakas2010random} (RAKEL) and multi-label twin support machine \cite{chen2016mltsvm} (MLTSVM). The regularization parameter $\lambda$ of these three algorithms depends on the type of datasets. 
For RAKEL, when $\ell = 3$ or $4$, the cardinality of sublabel set and the number of sublabel set are set as $3$ and $\ell$ respectively, and these two parameters are set as $3$ and $2\ell$ in other cases. 
For iNALM, we mainly tune $\lambda, \rho$ and $\mu$, and other parameters are set 
as in Subsection \ref{ex_svm}. 
More details about parameter setting will be illustrated in the numerical comparison part. 
%and other parameters are set as they had originally been given in those codes. For 0/1-IALM, we select $\lambda = 1, \rho = 10^{-2}, \mu = 10^{-2}$, and other parameters are set the same as those in Subsection \ref{ex_svm}.

Here we select three representative metrics (see e.g. \cite{zhang2013review}) for evaluating multi-label classifier, and they are Hamming loss (\texttt{HL}), ranking loss (\texttt{RL}) and average precision (\texttt{AP}).

%%%%%%%%%%%%%%%%%%%%%%%%%%%%%%%%%%%%%%%%%%%%%%%%%%%%%%%
\subsubsection{Numerical Comparison}

The comparison on Example \ref{examp2} aims to observe the performance of four MLC algorithms with various $m,n$ and $\ell$. For algorithms Rank-SVM, RAKEL and MLTSVM, regularization parameters are set as $\lambda = 1$. For our iNALM, we set $\lambda = 1$, $\mu = 10^{-2}$ and $\rho=1$. 
For the weighted constants in $f$, we set $\vartheta_i = 10^{-2}$ for $i \in \{n,2n,\cdots,\ell n\}$, and $\vartheta_i = 1$ otherwise.

{\bf (i) First Test.}
We generate datasets by Example \ref{examp2} with 
$m \in \{1000,2000,\cdots,5000  \}$ and $n = 5000, \ell = 3$. 
As shown in Table \ref{tab2.1}, iNALM performs best on \texttt{HL}. 
Meanwhile, Rank-SVM has the best \texttt{RL}, but its \texttt{Time} is the largest, more than 10 times that of iNALM. 
Although MLTSVM spends shorter \texttt{Time}, 
it fails on the other three metrics. 
We can also observe that \texttt{Time} of Rank-SVM grows from $9.159$s to $236.1$s when $m$ changes from $1000$ to $5000$, which means that it is sensitive to the magnitude of $m$. 

{\bf (ii) Second Test.}
We fix $m = 1000, \ell = 3$ 
and change $n \in \{6000,7000,\cdots,10000  \}$. 
Table \ref{tab2.2} shows that iNALM has the best performance on all the metrics. In particular, 
compared with other algorithms,
its \texttt{Time} is not sensitive to the change of $n$.

\begin{table}[htbp]
	%\TABLE
	\captionsetup{justification=justified}
	\caption{Numerical results of four algorithms on Example \ref{examp2} with $n = 5000, \ell = 3$ and $m \in \{1000,2000,\cdots,5000  \}$.} \label{tab2.1}
	%{\textcolor{blue}
	{\begin{tabular}{c|cccc|cccc}
			\hline
			\begin{tabular}[c]{@{}c@{}}$n=5000$\\ $\ell=3$\end{tabular} & \multicolumn{4}{c|}{\texttt{HL} $\downarrow$}                         & \multicolumn{4}{c}{\texttt{Time} $\downarrow$}                        \\ \hline
			$m$                                                         & \texttt{iNALM} & \texttt{Rank-SVM} & \texttt{RAKEL} & \texttt{MLTSVM} & \texttt{iNALM} & \texttt{Rank-SVM} & \texttt{RAKEL} & \texttt{MLTSVM} \\ \hline
			1000                                                        & 0.378          & 0.389             & 0.381          & 0.454           & 0.728          & 9.159             & 2.834          & 2.390           \\
			2000                                                        & 0.342          & 0.346             & 0.341          & 0.404           & 2.751          & 35.94             & 7.183          & 3.219           \\
			3000                                                        & 0.310          & 0.331             & 0.325          & 0.402           & 4.453          & 81.27             & 12.21          & 4.169           \\
			4000                                                        & 0.289          & 0.309             & 0.311          & 0.384           & 9.993          & 150.4             & 21.20          & 5.465           \\
			5000                                                        & 0.272          & 0.290             & 0.299          & 0.366           & 21.71          & 236.1             & 31.08          & 6.681           \\ \hline
			& \multicolumn{4}{c|}{\texttt{RL} $\downarrow$}                         & \multicolumn{4}{c}{\texttt{AP} $\uparrow$}                            \\ \hline
			$m$                                                         & \texttt{iNALM} & \texttt{Rank-SVM} & \texttt{RAKEL} & \texttt{MLTSVM} & \texttt{iNALM} & \texttt{Rank-SVM} & \texttt{RAKEL} & \texttt{MLTSVM} \\ \hline
			1000                                                        & 0.348          & 0.343             & 0.505          & 0.479           & 0.785          & 0.786             & 0.780          & 0.712           \\
			2000                                                        & 0.284          & 0.269             & 0.450          & 0.397           & 0.832          & 0.841             & 0.812          & 0.762           \\
			3000                                                        & 0.238          & 0.247             & 0.415          & 0.388           & 0.857          & 0.849             & 0.831          & 0.772           \\
			4000                                                        & 0.208          & 0.209             & 0.392          & 0.360           & 0.873          & 0.871             & 0.839          & 0.781           \\
			5000                                                        & 0.193          & 0.177             & 0.372          & 0.343           & 0.881          & 0.892             & 0.847          & 0.791           \\ \hline
	\end{tabular}} 
	
	{``$\uparrow$" means that the larger metric value, the better performance the algorithm has, while ``$\downarrow$" is the opposite situation.}
\end{table}

\begin{table}[htbp]
	%\TABLE
	\captionsetup{justification=justified} 
	\caption{Numerical results of four algorithms on Example \ref{examp2} with $m = 1000, \ell = 3$ and $n \in \{6000,7000,\cdots,10000  \}$}. \label{tab2.2}
	%{\textcolor{blue}
	{\begin{tabular}{c|cccc|cccc}
			\hline
			\begin{tabular}[c]{@{}c@{}}$m=1000$\\ $\ell=3$\end{tabular} & \multicolumn{4}{c|}{\texttt{HL} $\downarrow$}                         & \multicolumn{4}{c}{\texttt{Time} $\downarrow$}                        \\ \hline
			$n$                                                             & \texttt{iNALM} & \texttt{Rank-SVM} & \texttt{RAKEL} & \texttt{MLTSVM} & \texttt{iNALM} & \texttt{Rank-SVM} & \texttt{RAKEL} & \texttt{MLTSVM} \\ \hline
			6000                                                            & 0.375          & 0.413             & 0.402          & 0.427           & 0.788          & 7.559             & 3.349          & 3.344           \\
			7000                                                            & 0.386          & 0.412             & 0.402          & 0.435           & 0.803          & 8.832             & 3.802          & 4.998           \\
			8000                                                            & 0.398          & 0.419             & 0.412          & 0.421           & 0.817          & 8.957             & 4.350          & 6.725           \\
			9000                                                            & 0.402          & 0.432             & 0.416          & 0.423           & 0.845          & 10.03             & 4.890          & 8.881           \\
			10000                                                           & 0.404          & 0.420             & 0.416          & 0.431           & 0.924          & 10.36             & 5.397          & 10.97           \\ \hline
			& \multicolumn{4}{c|}{\texttt{RL} $\downarrow$}                         & \multicolumn{4}{c}{\texttt{AP} $\uparrow$}                            \\ \hline
			$n$                                                             & \texttt{iNALM} & \texttt{Rank-SVM} & \texttt{RAKEL} & \texttt{MLTSVM} & \texttt{iNALM} & \texttt{Rank-SVM} & \texttt{RAKEL} & \texttt{MLTSVM} \\ \hline
			6000                                                            & 0.351          & 0.368             & 0.529          & 0.460           & 0.789          & 0.778             & 0.770          & 0.727           \\
			7000                                                            & 0.350          & 0.370             & 0.544          & 0.468           & 0.797          & 0.775             & 0.763          & 0.720           \\
			8000                                                            & 0.369          & 0.374             & 0.544          & 0.478           & 0.783          & 0.780             & 0.762          & 0.722           \\
			9000                                                            & 0.377          & 0.393             & 0.562          & 0.472           & 0.780          & 0.765             & 0.757          & 0.723           \\
			10000                                                           & 0.382          & 0.387             & 0.551          & 0.485           & 0.770          & 0.769             & 0.762          & 0.712           \\ \hline
	\end{tabular}}
	{``$\uparrow$" means that the larger metric value, the better performance the algorithm has, while ``$\downarrow$" is the opposite situation.}
\end{table}

{\bf (iii) Third Test.} Finally we fix $m = 1000, n = 5000$ and change $\ell \in \{6,8,\cdots,14\}$. In most cases of this test, iNALM has the best results on all four metrics. Again, it has a significant advantage on \texttt{Time}.

\begin{table}[htbp]
	%\TABLE
	\centering 
	\caption{Numerical results of four algorithms on Example \ref{examp2} with $m = 1000, n = 5000$ and $\ell \in \{6,8,\cdots,14  \}$.} \label{tab2.3}
	%{\textcolor{blue}
	{\begin{tabular}{c|cccc|cccc}
			\hline
			\begin{tabular}[c]{@{}c@{}}$m = 1000$\\ $n=5000$\end{tabular} & \multicolumn{4}{c|}{\texttt{HL} $\downarrow$}                         & \multicolumn{4}{c}{\texttt{Time} $\downarrow$}                        \\ \hline
			$\ell$                                                      & \texttt{iNALM} & \texttt{Rank-SVM} & \texttt{RAKEL} & \texttt{MLTSVM} & \texttt{iNALM} & \texttt{Rank-SVM} & \texttt{RAKEL} & \texttt{MLTSVM} \\ \hline
			6                                                           & 0.392          & 0.432             & 0.432          & 0.515           & 1.585          & 17.55             & 20.41          & 4.430           \\
			8                                                           & 0.391          & 0.404             & 0.400          & 0.466           & 2.055          & 22.10             & 27.24          & 5.823           \\
			10                                                          & 0.397          & 0.411             & 0.415          & 0.474           & 2.568          & 37.27             & 34.24          & 7.294           \\
			12                                                          & 0.389          & 0.386             & 0.386          & 0.446           & 3.053          & 63.38             & 40.84          & 8.868           \\
			14                                                          & 0.393          & 0.429             & 0.431          & 0.484           & 3.780          & 102.5             & 47.55          & 10.27           \\ \hline
			& \multicolumn{4}{c|}{\texttt{RL} $\downarrow$}                         & \multicolumn{4}{c}{\texttt{AP} $\uparrow$}                            \\ \hline
			$\ell$                                                      & \texttt{iNALM} & \texttt{Rank-SVM} & \texttt{RAKEL} & \texttt{MLTSVM} & \texttt{iNALM} & \texttt{Rank-SVM} & \texttt{RAKEL} & \texttt{MLTSVM} \\ \hline
			6                                                           & 0.350          & 0.401             & 0.435          & 0.459           & 0.746          & 0.724             & 0.698          & 0.657           \\
			8                                                           & 0.346          & 0.362             & 0.430          & 0.441           & 0.727          & 0.718             & 0.702          & 0.676           \\
			10                                                          & 0.340          & 0.361             & 0.427          & 0.486           & 0.721          & 0.715             & 0.692          & 0.642           \\
			12                                                          & 0.346          & 0.343             & 0.413          & 0.443           & 0.706          & 0.711             & 0.701          & 0.640           \\
			14                                                          & 0.354          & 0.386             & 0.455          & 0.489           & 0.696          & 0.676             & 0.650          & 0.601           \\ \hline
	\end{tabular}}
	
	{``$\uparrow$" means that the larger metric value, the better performance the algorithm has, while ``$\downarrow$" is the opposite situation.}
\end{table}

\begin{table}[H]
	%\TABLE
	\centering 
	\caption{The \texttt{HL} and \texttt{Time} of four algorithms on real multi-label classification datasets.} \label{mlc_real1}
	%{\textcolor{blue}
	{\begin{tabular}{c|cccc|cccc}
			\hline
			& \multicolumn{4}{c|}{\texttt{HL} $\downarrow$}                         & \multicolumn{4}{c}{\texttt{Time} $\downarrow$}                        \\ \cline{2-9} 
			& \texttt{iNALM} & \texttt{Rank-SVM} & \texttt{RAKEL} & \texttt{MLTSVM} & \texttt{iNALM} & \texttt{Rank-SVM} & \texttt{RAKEL} & \texttt{MLTSVM} \\ \hline
			\texttt{sbbc} & 2.04e-1        & 1.93e-1           & 2.29e-1        & 2.80e-1         & 0.178          & 6.235             & 0.065          & 0.205           \\
			\texttt{guar} & 2.18e-1        & 2.31e-1           & 2.45e-1        & 2.91e-1         & 0.163          & 3.500             & 0.066          & 0.180           \\
			\texttt{genb} & 1.02e-3        & 1.87e-3           & 2.89e-3        & 1.02e-2         & 0.273          & 214.2             & 0.085          & 0.192           \\
			\texttt{reut} & 2.28e-1        & 2.08e-1           & 2.56e-1        & 2.89e-1         & 0.167          & 5.575             & 0.080          & 0.192           \\
			\texttt{eugo} & 2.21e-2        & --                & 2.32e-2        & 3.11e-2         & 7.038          & --                & 6.273          & 3028            \\
			\texttt{gpgo} & 3.05e-2        & 4.07e-2           & 2.76e-2        & 5.67e-2         & 0.147          & 4.851             & 0.022          & 0.150           \\
			\texttt{hugo} & 3.87e-2        & 4.89e-2           & 4.17e-2        & 5.76e-2         & 2.215          & 485.5             & 3.526          & 54.45           \\
			\texttt{medi} & 1.02e-2        & 1.14e-2           & 1.11e-2        & 1.37e-2         & 1.070          & 889.8             & 0.624          & 5.314           \\
			\texttt{plgo} & 3.60e-2        & 4.91e-2           & 4.25e-2        & 5.24e-2         & 0.533          & 61.13             & 0.621          & 3.653           \\
			\texttt{rcs1} & 2.66e-2        & --                & 2.74e-2        & 3.17e-2         & 30.45          & --                & 54.10          & 5792            \\
			\texttt{rcs2} & 2.28e-2        & --                & 2.36e-2        & 3.15e-2         & 36.16          & --                & 61.19          & 5696            \\
			\texttt{rcs3} & 2.29e-2        & --                & 2.39e-2        & 3.37e-2         & 35.40          & --                & 62.66          & 5694            \\
			\texttt{rcs4} & 1.94e-2        & --                & 1.99e-2        & 2.58e-2         & 37.70          & --                & 62.46          & 4847            \\
			\texttt{rcs5} & 2.23e-2        & --                & 2.32e-2        & 2.77e-2         & 35.54          & --                & 60.42          & 5393            \\
			\texttt{yhbs} & 2.48e-2        & --                & 2.56e-2        & --              & 116.4          & --                & 206.0          & ---             \\
			\texttt{yhet} & 4.58e-2        & --                & 4.46e-2        & --              & 234.3          & --                & 363.7          & ---             \\
			\texttt{yhhl} & 3.11e-2        & --                & 3.34e-2        & --              & 129.4          & --                & 197.2          & --              \\
			\texttt{scen} & 1.08e-1        & 1.89e-1           & 1.21e-1        & 1.65e-1         & 0.631          & 62.51             & 6.163          & 44.77           \\
			\texttt{yelp} & 1.65e-1        & 2.20e-1           & 1.69e-1        & 2.44e-1         & 8.995          & 2968              & 9.594          & 3139            \\
			\texttt{hpaa} & 8.80e-2        & 1.16e-1           & 1.16e-1        & 1.28e-1         & 10.79          & 588.3             & 26.18          & 1639            \\
			\texttt{20ng} & 3.19e-2        & --                & 3.84e-2        & --              & 3.030          & --                & 29.81          & --              \\
			\texttt{ohsu} & 5.88e-2        & --                & 7.44e-2        & --              & 7.276          & --                & 29.14          & --              \\
			\texttt{reuk} & 1.25e-2        & --                & 1.55e-2        & --              & 1.745          & --                & 6.005          & --              \\
			\texttt{slas} & 4.05e-2        & 4.85e-2           & 6.29e-2        & 7.23e-2         & 0.253          & 1946              & 0.953          & 212.7           \\ \hline
	\end{tabular}}
	
	{``$\uparrow$" means that the larger metric value, the better performance the algorithm has, while ``$\downarrow$" is the opposite situation. ``--" means that the corresponding algorithm fails to compute a solution within 3 hours or runs out of memory.}
\end{table}

Now let us show the numerical results of the four algorithms on the real data sets in Example~\ref{example_mlc_real}. 
Regularization parameters for Rank-SVM, RAKEL and MLTSVM are set as $\lambda = 10^3$, and other parameters are set as default values. 
For iNALM, we set $ \rho = 10^2$ and $\lambda = 10^{3}$, while the weighted constants in $f$ are selected as $\vartheta_i = 1$ for $i \in \{ 1,\cdots,p\ell \}$. For parameter $\mu$ of iNALM, we set $\mu = 10^{-2}$ on \texttt{sbbc}, \texttt{genb}, \texttt{guar} and \texttt{reut}, $\mu = 10^{3}$ on \texttt{yhbs}, \texttt{yhet} and \texttt{hpaa}, while $\mu = 10^2$ on the rest of datasets. We can see from Table \ref{mlc_real1} and \ref{mlc_real2} that iNALM can achieve the best \texttt{HL}, \texttt{RL} and \texttt{AP} on most of the datasets in Example \ref{example_mlc_real}. Particularly, RAKEL spends the shortest \texttt{Time} on datasets with $m_{tr} \leq 1000$ and $n \leq 2000$, whereas iNALM has the shortest \texttt{Time} on other datasets.

\begin{table}[H]
	%\TABLE
	\centering 
	\caption{The \texttt{RL} and \texttt{AP} of four algorithms on real multi-label classification datasets.} \label{mlc_real2}
	%{\textcolor{blue}
	{\begin{tabular}{c|cccc|cccc}
			\hline
			& \multicolumn{4}{c|}{\texttt{RL} $\downarrow$}                         & \multicolumn{4}{c}{\texttt{AP} $\uparrow$}                           \\ \cline{2-9} 
			& \texttt{iNALM} & \texttt{Rank-SVM} & \texttt{RAKEL} & \texttt{MLTSVM} & \texttt{iNALM} & \texttt{Rank-SVM} & \texttt{RAKEL} & \texttt{MLTSVM} \\ \hline
			\texttt{sbbc} & 3.08e-1        & 3.57e-1           & 6.82e-1        & 3.16e-1         & 0.595          & 0.571             & 0.542          & 0.573           \\
			\texttt{guar} & 3.39e-1        & 4.05e-1           & 7.43e-1        & 3.89e-1         & 0.564          & 0.499             & 0.494          & 0.519           \\
			\texttt{genb} & 2.15e-3        & 8.68e-4           & 7.65e-3        & 1.96e-3         & 0.993          & 0.995             & 0.515          & 0.991           \\
			\texttt{reut} & 3.35e-1        & 3.95e-1           & 7.10e-1        & 3.41e-1         & 0.563          & 0.538             & 0.494          & 0.565           \\
			\texttt{eugo} & 2.48e-2        & --                & 1.45e-1        & 3.38e-2         & 0.881          & --                & 0.841          & 0.815           \\
			\texttt{gpgo} & 2.71e-2        & 3.49e-2           & 5.14e-2        & 3.88e-2         & 0.969          & 0.960             & 0.966          & 0.954           \\
			\texttt{hugo} & 3.63e-2        & 3.87e-2           & 1.60e-1        & 5.26e-2         & 0.877          & 0.858             & 0.836          & 0.795           \\
			\texttt{medi} & 2.12e-2        & 1.49e-2           & 1.72e-1        & 2.91e-2         & 0.900          & 0.888             & 0.829          & 0.836           \\
			\texttt{plgo} & 3.14e-2        & 3.23e-2           & 2.02e-1        & 5.06e-2         & 0.889          & 0.881             & 0.815          & 0.824           \\
			\texttt{rcs1} & 4.27e-2        & --                & 5.52e-1        & 8.67e-2         & 0.631          & --                & 0.435          & 0.499           \\
			\texttt{rcs2} & 3.99e-2        & --                & 5.17e-1        & 8.04e-2         & 0.660          & --                & 0.483          & 0.550           \\
			\texttt{rcs3} & 3.82e-2        & --                & 5.20e-1        & 8.40e-2         & 0.660          & --                & 0.490          & 0.548           \\
			\texttt{rcs4} & 3.13e-2        & --                & 4.41e-1        & 7.01e-2         & 0.714          & --                & 0.562          & 0.605           \\
			\texttt{rcs5} & 3.90e-2        & --                & 4.94e-1        & 7.98e-2         & 0.657          & --                & 0.301          & 0.547           \\
			\texttt{yhbs} & 4.70e-2        & --                & 1.59e-1        & --              & 0.888          & --                & 0.487          & --              \\
			\texttt{yhet} & 9.66e-2        & --                & 2.96e-1        & --              & 0.764          & --                & 0.741          & --              \\
			\texttt{yhhl} & 7.11e-2        & --                & 2.77e-1        & --              & 0.803          & --                & 0.758          & --              \\
			\texttt{scen} & 1.19e-1        & 1.05e-1           & 2.38e-1        & 1.44e-1         & 0.823          & 0.816             & 0.782          & 0.768           \\
			\texttt{yelp} & 1.95e-1        & 2.02e-1           & 2.00e-1        & 2.54e-1         & 0.777          & 0.750             & 0.814          & 0.686           \\
			\texttt{hpaa} & 3.69e-1        & 3.56e-1           & 6.19e-1        & 2.22e-1         & 0.469          & 0.408             & 0.442          & 0.450           \\
			\texttt{20ng} & 5.64e-2        & --                & 2.64e-1        & --              & 0.814          & --                & 0.734          & --              \\
			\texttt{ohsu} & 2.31e-1        & --                & 3.83e-1        & --              & 0.615          & --                & 0.580          & --              \\
			\texttt{reuk} & 1.01e-1        & --                & 4.75e-1        & --              & 0.635          & --                & 0.501          & --              \\
			\texttt{slas} & 1.22e-1        & 1.63e-1           & 4.07e-1        & 1.27e-1         & 0.678          & 0.625             & 0.555          & 0.553           \\ \hline
	\end{tabular}}
	
	{``$\uparrow$" means that the larger metric value, the better performance the algorithm has, while ``$\downarrow$" is the opposite situation. ``--" means that the corresponding algorithm fails to compute a solution within 3 hours or runs out of memory.}
\end{table}

\begin{remark} \label{Remark-Numerical}
	One major reason why iNALM is fast for both SVM and MLC problems is that
	the objective function $f(\bfx)$ has a separable form: $f(\bfx) = \sum_{i=1}^n w_i f_i(x_i)$, where $w_i >0$ and $f_i$ is twice continuously differentiable and convex.
	Therefore, the Hessian matrix $\nabla^2 f(\bfx)$ is diagonal and the Newton equation
	of the type (\ref{Reduced-Newton-Eq})
	encountered in iNALM is cheap to solve. 
	Moreover, the number of support vectors is small. 
	This in turn reduces the computational complexity when it comes to matrix-vector 
	multiplications (e.g., computing $A_{\Gamma} \bfx$).
	Finally, we make an interesting observation about the performance of iNALM on the test
	problems that have $m >n$. 
	In theory, iNALM requires $m \le n$ so that the full row-rank assumption on the data 
	matrix $A$ may be satisfied. 
	We note that iNALM also worked quite well when this assumption was violated.
	This brought us to think whether the full row-rank assumption may be weakened or replaced by other conditions that would allow $m >n$. 
	It is certainly an important research question. 
	%Another reason is that the considered problems all have $m < n$, making it possible that
	%$A$ has full row-rank (Assumption~\ref{Assumption-Fullrank}).
	%This is important to the convergence rate results of iNALM.
\end{remark}

%%%%%%%%%%%%%%%%%%%%%%%%%%%%%%%%%%%%%%%%%%%%%%%%%%%%%%%%%%%
\section{Conclusion} \label{Section-Conclusion}

This paper aims to answer an open question how to extend the classical theory
of augmented Lagrangian method (ALM) involving second-order conditions from smooth optimization 
to nonsmooth, nonconvex optimization.
For the 0/1 composite optimization (\ref{COP}),
we demonstrated that it can be achieved through 
successfully identifying the active set by
making use of the proximal operator of the 0/1 loss function $h(\bfu)$.
Consequently, we are able to define a second-order necessary/sufficient condition that
is essential for the convergence of ALM.
Since the active-set defines a subspace, it is natural to use Newton's method within
this subspace to solve the subproblem arising from ALM. 
The resulting subspace Newton's method is proved to be globally and locally quadratically convergent under reasonable conditions.

For the Newton method to be implementable, we designed a set of stopping criteria that
are closely related to a local error bound for the KKT system of (\ref{COP}).
This local error bound allows us to prove the R-linear convergence rate of the resulting inexact Newton ALM under the proposed second-order sufficient condition.
Therefore, we satisfactorily extended the classical ALM theory to (\ref{COP}) via an
implementable inexact Newton ALM. 

This research raises a hope that similar theory may be obtained for a wider class of
nonsmooth, nonconvex optimization problems considered in \cite{bolte2018nonconvex,bot2019proximal}. A key question to answer here is what form a second-order sufficient 
condition (SOSC) may take and how their ALMs are related to such SOSC properties.
It seems that there is no easy answer to those questions, which are certainly worth for
further investigation.

$\newline$
\section*{Acknowledgment}
This work was supported by Fundamental Research Funds for the Central Universities (2022YJS099), the National Natural Science Foundation of China (12131004, 11971052) and Beijing Natural Science Foundation (Z190002).

%%%%%%%%%%%%%%%%%%%%%%%%%%%%%%%%%%%%%%%%%%%%%%%%%%%%%%%%%%%%%%
\bibliographystyle{siamplain}
\bibliography{NALM_Refs} % if more than one, comma separated

\end{document}